\documentclass[reqno,english]{amsart}
\usepackage{amsfonts,amsmath,latexsym,verbatim,amscd,mathrsfs,color,array}
\usepackage[colorlinks=true]{hyperref}
\usepackage{amsmath,amssymb,amsthm,amsfonts,graphicx,color}
\usepackage[hmargin=3cm, vmargin=3cm]{geometry}
\usepackage{bbm}
\usepackage{amssymb}
\usepackage{pdfsync}
\usepackage{epstopdf}
\usepackage{cite}
\usepackage{graphicx}
\usepackage{multirow}
\usepackage[font=bf,aboveskip=15pt]{caption}
\usepackage[toc,page]{appendix}
\usepackage[colorlinks=true]{hyperref}
\newcommand\la{\lambda}

\newcommand{\inn}{{\quad\mbox{in } }}

\newcommand\De{\Delta}

\newcommand{\R}{\mathbb{R}}

\renewcommand{\Re} {\mathop{\mathrm{Re}}}
\renewcommand{\Im} {\mathop{\mathrm{Im}}}
\renewcommand{\div}{\mathop{\rm div}}
\newcommand{\curl} {\mathop{\rm curl}}

\newcommand{\pd}{\partial}

\newcommand{\na}{\nabla}
\newcommand{\ttt}{\tilde}

\newcommand{\RR}{\mathcal{R}}
\newcommand{\EE}{\mathcal{E}}
\newcommand{\KK}{\mathcal{K}}
\newcommand{\EQ}[1]{\begin{equation}\begin{split} #1 \end{split}\end{equation}}

\renewcommand{\div}{\mathop{\rm div}}

\newtheorem{theorem}{Theorem}[section]
\newtheorem{lemma}{Lemma}[section]
\newtheorem{cor}{Corollary}[section]
\newtheorem{prop}{Proposition}[section]
\newtheorem{remark}{Remark}[section]

\numberwithin{equation}{section}

\begin{document}
\title[Nematic liquid crystal flow with partially free boundary]{Nematic liquid crystal flow with partially free boundary}

\author[F. Lin]{Fanghua Lin}
\address{\noindent
Courant Institute of Mathematical Sciences,
New York University, NY 10012, USA}
\email{linf@cims.nyu.edu}

\author[Y. Sire]{Yannick Sire}
\address{\noindent
Department of Mathematics,
Johns Hopkins University, 3400 N. Charles Street, Baltimore, MD 21218, USA}
\email{ysire1@jhu.edu}

\author[J. Wei]{Juncheng Wei}
\address{\noindent
Department of Mathematics,
University of British Columbia, Vancouver, B.C., V6T 1Z2, Canada}
\email{jcwei@math.ubc.ca}

\author[Y. Zhou]{Yifu Zhou}
\address{\noindent
Department of Mathematics,
Johns Hopkins University, 3400 N. Charles Street, Baltimore, MD 21218, USA}
\email{yzhou173@jhu.edu}

\begin{abstract}
We study a simplified Ericksen-Leslie system modeling the flow of nematic liquid crystals with partially free boundary conditions. It is a coupled system between the Navier-Stokes equation for the fluid velocity with a transported heat flow of harmonic maps, and both of these parabolic equations are critical for analysis in  two dimensions. The boundary conditions are physically natural and they correspond to the Navier slip boundary condition with zero friction for the velocity field and a Plateau-Neumann type boundary condition for the map. In this paper we construct smooth solutions of this coupled system that blow up in  finite time at any finitely many given points on the boundary or in the interior of the domain.
\end{abstract}
\maketitle

{
  \hypersetup{linkcolor=black}
  \tableofcontents
}


\bigskip

\section{Introduction}

The aim of the present work is to investigate liquid crystal flows with partially free boundary conditions. Let $\Omega \subset \mathbb{R}^d$ $(d\leq 3)$ be a smooth domain.  We consider the following system
\begin{equation}\label{LCF}
\begin{cases}
\partial_t v + v\cdot \nabla v +\nabla P = \Delta v - \varepsilon_0\nabla\cdot \left(\nabla u \odot \nabla u-\frac12 |\nabla u|^2 \mathbb{I}_d\right)~&\mbox{ in }~\Omega\times (0,T),\\
\nabla\cdot v =0~&\mbox{ in }~\Omega\times (0,T),\\
\partial_t u+v\cdot\nabla u=\Delta u +|\nabla u|^2 u~&\mbox{ in }~\Omega\times (0,T),\\
\end{cases}
\end{equation}
with the partially free boundary conditions
\begin{equation}\label{FB}
\begin{cases}
v\cdot \nu=0~&\mbox{ on }~\partial\Omega\times(0,T),\\
(S v \cdot \nu)_{\tau} = 0~&\mbox{ on }~\partial\Omega\times(0,T),\\
u(x,t) \in \Sigma~&\mbox{ on }~\partial\Omega\times(0,T),\\
\frac{\partial u}{\partial \nu}(x,t) \perp T_{u(x,t)} \Sigma~&\mbox{ on }~\partial\Omega\times(0,T),\\
\end{cases}
\end{equation}
where $v:\Omega\times[0,T)\to\R^d$ is the fluid velocity field, $P:\Omega\times[0,T)\to\R$ is the fluid pressure function, $u:\Omega\times[0,T)\to\mathbb{S}^2$ stands for the orientation unit vector field of nematic liquid crystals, $(\na u\odot\na u)_{ij}:=\na_i u\cdot\na_j u$, and $\mathbb I_d$ is the identity matrix on $\R^d$, $\varepsilon_0>0$ (coupling constant) represents a competition between kinetic energy and elastic energy, $\nu$ is the unit outer normal of $\partial\Omega$, $S$ is the strain tensor 
$$
S v=\frac{1}{2}(\nabla v + (\nabla v)^T),
$$
 and $\Sigma\subset \mathbb S^2$ is a simple, closed and smooth curve. Let us assume that $\Sigma$ is the equator for simplicity. The case that $\Sigma$ is a circle in $\mathbb S^2$ is physically relevant, see for example \cite{Ericksen1976,DeGennes1974}. And the proofs here work equally well with some simple modifications.  Boundary conditions for \eqref{FB}$_1$--\eqref{FB}$_2$ are the usual {\it Navier boundary conditions} for the Navier--Stokes equation with a zero friction, and \eqref{FB}$_3$--\eqref{FB}$_4$ are referred as {\it partially free boundary conditions} for the harmonic map heat flow, see for examples \cite{Hamilton1975LNM,MaLiCMH1991,StruweManMath1991,ChenLin98JGA} and the references therein.

\medskip

This system \eqref{LCF}, which was first introduced by the first author in \cite{L1989CPAM} as a simplified version of Ericksen--Leslie system established by Ericksen \cite{Ericksen1962ARMA} and Leslie
\cite{Leslie1968ARMA}, enjoys the same type energy law, coupling structure and dissipative properties. The system under consideration is a nonlinearly coupled system between the incompressible Navier-Stokes equations and the heat flow of harmonic maps with a (partially) free boundary condition. The latter is a geometric flow with the Plateau and Neumann type boundary conditions. Let us first describe briefly the latter system
in a more geometric set up.

\medskip

Let $(M, g)$ be an $m$-dimensional smooth Riemannian manifold with boundary $\partial M$ and $N$ be another  smooth  Riemannian manifold without boundary. Suppose $\Sigma$ is a $k$-dimensional submanifold of $N$ without boundary. Any continuous map $u_0: M\to N$ satisfying $u_0(\partial M)\subset \Sigma$ defines a relative homotopy class in maps from $(M, \partial M)$ to $(N, \Sigma)$. A map $u: M\to N$ with $u(\partial M)\subset \Sigma$ is called homotopic to $u_0$ if there exists a continuous homotopy $h:[0, 1]\times M \to N$ satisfying $h([0,1]\times \partial M)\subset \Sigma$, $h(0) = u_0$ and $h(1) = u$. An interesting problem is that whether or not each relative homotopy class of maps has a representation by harmonic maps. The latter must be solutions to the following problem:
\begin{equation}\label{e:harmonicmapwithfreeboundary}
\begin{cases}
-\Delta u = \Gamma(u)(\nabla u, \nabla u),\\
u(\partial M)\subset \Sigma,\\
\frac{\partial u}{\partial \nu}\perp T_u\Sigma.
  \end{cases}
\end{equation}
Here $\nu$ is the unit normal vector of $M$ along the boundary $\partial M$, $\Delta \equiv \Delta_M$ is the Laplace-Beltrami operator of $(M,g)$, $\Gamma$ is the second fundamental form of $N$ (viewed as a submanifold in $\mathbb{R}^\ell$ via Nash's isometric embedding), $T_pN$ is the tangent space in $\mathbb{R}^\ell$ of $N$ at $p$ and $\perp$ means orthogonal to in $\mathbb{R}^\ell$. (\ref{e:harmonicmapwithfreeboundary}) is the Euler--Lagrange equation for critical points of the Dirichlet energy functional
\begin{equation*}
E(u) = \int_{M}|\nabla u|^2\,dv_g
\end{equation*}
defined over the space of maps
\begin{equation*}
H^1_\Sigma(M, N) = \{u\in H^1(M, N): u(x)\subset \Sigma \ {\rm{a.e.}}\ x\in\partial M \}.
\end{equation*}
Here $H^1(M, N)=\big\{u\in H^1(M,\mathbb R^\ell): u(x)\in N \ {\rm{a.e.}}\ x\in M\big\}$. Both the existence and partial regularity of energy minimizing harmonic maps in $H^1_\Sigma(M, N)$ have been established for examples, in \cite{Hamilton1975LNM,BaldesMM1982,GulliverJostJRAM1987} under special assumptions, and in \cite{DuzaarSteffen1989JRAM,DuzaarSteffenAA1989,HardtLinCPAM1989} in general cases. Another standard approach to investigate (\ref{e:harmonicmapwithfreeboundary}) is to study the following parabolic problem
\begin{equation}\label{e:harmonicmapflowwithfreeboundary}
\begin{cases}
\partial_t u -\Delta u = \Gamma(u)(\nabla u, \nabla u)&\text{ on }M\times [0, \infty),\\
u(x, t)\in \Sigma &\text{ on } \partial M\times [0,\infty),\\
\frac{\partial u}{\partial \nu}(x, t)\perp T_{u(x, t)}\Sigma &\text{ on}\ \partial M\times [0,\infty),\\
u(\cdot, 0) = u_0 &\text{ on }M.
  \end{cases}
\end{equation}
This is the so-called harmonic map heat flow with a partially free boundary. (\ref{e:harmonicmapflowwithfreeboundary}) was first studied in \cite{Hamilton1975LNM}, Hamilton considered the case where $\Sigma$ is totally geodesic and the sectional curvature $K_N\leq 0$.
He proved the existence of a unique global smooth solution for (\ref{e:harmonicmapflowwithfreeboundary}). The global existence of weak solutions of (\ref{e:harmonicmapflowwithfreeboundary}) was established by Struwe in \cite{StruweManMath1991} for $m\ge 3$, see also \cite{ChenLin98JGA}.
 In \cite{MaLiCMH1991}, the case $m = {\rm dim} M = 2$ was considered, where a global existence and uniqueness result for finite energy weak solutions was obtained under some {suitable} geometric hypotheses on $N$ and $\Sigma$.  When $N$ is an Euclidean space, the first equation in (\ref{e:harmonicmapflowwithfreeboundary}) becomes the standard heat equation
\begin{equation*}
u_t - \Delta u = 0 \text{ on }M \times [0, \infty).
\end{equation*}
Even in this special case, as pointed out in \cite{ChenLin98JGA} and \cite{StruweManMath1991}, estimates near the boundary for (\ref{e:harmonicmapflowwithfreeboundary}) are difficult due to this highly nonlinear boundary condition. As far as the heat flow is concerned, Struwe in \cite{StruweManMath1991} studied the problem using the intrinsic version of harmonic maps with a free boundary condition. In particular, he used a Ginzburg-Landau approximation in the interior, hence keeping the same nonlinear boundary condition. Another approach was considered in \cite{HSSW}, where the approximation is on the boundary.

\medskip

The finite time singularity (as conjectured in \cite{ChenLin98JGA}) for \eqref{e:harmonicmapflowwithfreeboundary} was proven only recently in \cite{sire2019singularity} with  $N=\R^2$, $M=\R^2_+$ and $\Sigma=\mathbb S^1\subset \R^2$. The analysis there cannot be generalized directly to the current situation as the target is no longer flat and the standard heat equation has to be replaced by the heat flow of maps into the sphere. Despite the nonlinearity of the system and the nonlinear coupling, the Navier slip boundary condition turns out to be consistent with the partially free boundary condition of the map. The latter is important for our analysis.

\medskip

In the aspect of the incompressible Navier--Stokes equations, we refer the readers to, for example, the books \cite{Constantin-Foias,Tartar-Book} and the references therein for comprehensive theories. Of particular interest here is the incompressible Navier-Stokes equation with Navier boundary conditions since the system \eqref{LCF} turns out to be more compatible and physically natural with the Navier boundary conditions \eqref{FB}$_1$--\eqref{FB}$_2$ compared to the no-slip boundary condition. In the setting of Navier boundary conditions, without being exhaustive, we refer to \cite{daVeiga2006CPAA} for Constantin-Fefferman type regularity results in the case of $\R^3_+$ and in \cite{daVeiga2009,SiranLi2019} for the case of general domains in $\R^3$. See also \cite{Neustupa2010,Chen2010Indiana} for the local existence of strong solutions and \cite{Chen2010Indiana} for the global existence of weak solutions in dimension three. It is worth mentioning that in order to treat boundaries, the authors in \cite{daVeiga2009,SiranLi2019} used the Solonnikov's theory developed in \cite{Solonnikov1971II,Solonnikov1970I} on Green's matrices for elliptic systems of Petrovsky type, which is a subclass of Agmon--Douglis--Nirenberg (ADN) elliptic systems (see the seminal works \cite{ADN1,ADN2}). Our case, especially for the dealing with the boundary bubbling, is closely related to the aforementioned theories.

\medskip

For the study of the nematic liquid crystal flows, there have been growing interests concerning the global existence of weak solutions, partial regularity results, singularity formation and others. We refer to \cite{L1989CPAM,LL1995CPAM,LL1996DCDS,LLW2010ARMA,LW2014RSTA,LW2016CPAM,HLLW2016ARMA,LCF2D} and the references therein.

\medskip

\subsection*{Main results}
In the following, we consider the nematic liquid crystal flow with partially free boundary \eqref{LCF}--\eqref{FB} in the half space case $\Omega=\R^2_+$. Our method of construction could be adapted to the case of general domain, but it would involve more technical computations and we refrain considering such a generality.

\medskip

We first construct finite time blow-up solutions to the partially free boundary system \eqref{LCF}--\eqref{FB}, where the singularities can actually take place both in the interior and on the boundary, and as a direct consequence of the construction, we give an example (different from the one constructed in \cite{sire2019singularity}) of finite time singularities for the harmonic map heat flow with partially free boundary \eqref{e:harmonicmapflowwithfreeboundary} as conjectured in \cite{ChenLin98JGA}. 

\medskip

Our first theorem is stated as follows
\begin{theorem}\label{thm1}
For $T,\varepsilon_0>0$ sufficiently small and any given points $\{q^{(j)}_{\mathcal B}\}_{j=1}^{k_{\mathcal B}}\cup\{q^{(j)}_{\mathcal I}\}_{j=1}^{k_{\mathcal I}}\subset\overline{\R^2_+}$ with $q^{(j)}_{\mathcal B}\in \partial \R^2_+$ and $q^{(j)}_{\mathcal I}\in\mathring{\R}^2_+$, there exists initial data $(u_0,v_0)$ such that the solution $(u,v)$ to problem \eqref{LCF} with partially free boundary conditions \eqref{FB} blows up at finite time $t=T$ exactly at these given points. More precisely,
\begin{equation*}
u(x,t)-u_*(x)-\sum_{j=1}^{k_{\mathcal B}}  \left[W_1\left(\frac{x-q^{(j)}_{\mathcal B}}{\la^{(j)}(t)}\right)-W_1(\infty)\right]-\sum_{j=1}^{k_{\mathcal I}} Q_{\omega_{\mathcal I}^{(j)}} \left[W_2\left(\frac{x-q^{(j)}_{\mathcal I}}{\la^{(j)}(t)}\right)-W_2(\infty)\right]\to 0 ~~\mbox{ as }~t\to T
\end{equation*}
in $H^1(\R^2_+)\cap L^\infty(\R^2_+)$, where $u_*\in H^1(\R^2_+)\cap C(\bar\R^2_+)$, profiles $W_1$ and $W_2$ are defined in \eqref{def-W1} and \eqref{def-W2}, respectively, the rotation $Q_{\omega_{\mathcal I}^{(j)}}$ is defined in \eqref{def-Qz}, and the blow-up rate and angles satisfy, for some $\kappa_j^*>0$ and $\omega_j^*$,
\begin{equation*}
\begin{aligned}
&\la^{(j)}(t)\sim\kappa_j^*\frac{T-t}{|\log(T-t)|^2}~\mbox{ as }~t\to T,\\
&\omega^{(j)}_{\mathcal I} \to \omega_j^{*} ~\ \mbox{ as } \ t \to T.\\
\end{aligned}
\end{equation*}
In particular, it holds that
\begin{equation*}
|\nabla u(\cdot,t)|^2\,dx\rightharpoonup |\nabla u_*|^2\,dx +4\pi \sum_{j=1}^{k_{\mathcal B}} \delta_{q^{(j)}_{\mathcal B}}+8\pi \sum_{j=1}^{k_{\mathcal I}} \delta_{q^{(j)}_{\mathcal I}}~\mbox{ as }~t\to T
\end{equation*}
as convergence of Radon measures.
Furthermore, the velocity field satisfies
\begin{equation*}
|v(x,t)|\leq c\sum_{j=1}^k  \frac{\la_j^{\nu_{j}-1}(t)}{1+\left|\frac{x-q_j}{\la_j(t)}\right|}, ~~0<t<T,
\end{equation*}
for some $c>0$ and $0<\nu_j<1$, $j=1,\cdots,k$. Here $k=k_{\mathcal B}+k_{\mathcal I}$ and $\{q_j\}_{j=1}^k=\{q^{(j)}_{\mathcal B}\}_{j=1}^{k_{\mathcal B}}\cup\{q^{(j)}_{\mathcal I}\}_{j=1}^{k_{\mathcal I}}$.
\end{theorem}

\medskip

\begin{remark}
\noindent
\begin{itemize}
    \item Each bubble on the boundary might be viewed as a ``half'' bubble.
    \item The absence of the rotations for the boundary bubbles is in fact a consequence of the partially free boundary conditions \eqref{FB}. See more detailed discussions in next subsection.
\end{itemize}
\end{remark}

\medskip

For the harmonic map heat flow with free boundary, the question whether finite time singularity exists or not was originally raised by Chen and Lin \cite{ChenLin98JGA}. The first example was constructed recently by Sire, Wei and Zheng \cite{sire2019singularity} using a caloric extension. We would like to point out that the proof of Theorem \ref{thm1} actually gives another different example of finite time singularity. In fact, as a consequence of the construction of Theorem \ref{thm1}, we have
\begin{cor}
 Assume $M=\R^2_+$, $N=\mathbb S^2$, and $\Sigma=\{(x_1,x_2,x_3)\in\mathbb S^2:x_3=0\}$ in \eqref{e:harmonicmapflowwithfreeboundary}.
Given any finitely many distinct points $q_k$ in $\R^2_+$ or on $\partial \R^2_+$, for $T>0$ sufficiently small, there exists initial data $u_0$ such that the solution to \eqref{e:harmonicmapflowwithfreeboundary} blows up exactly at these prescribed points at time $t=T$. Moreover, the blow-up profile takes the form of sharply scaled 1-corotational profile around each point $q_k$ with type II blow-up rate
$$
\la_k(t)\sim \frac{T-t}{|\log(T-t)|^2}~\mbox{ as }~t\to T.
$$
\end{cor}

\medskip

In order to deal with the Navier--Stokes equation with Navier boundary conditions, we consider the following Stokes system with Navier boundary conditions
\begin{equation}\label{eqn-Stokes-Navier}
\left\{
\begin{aligned}
&\partial_t v  +\nabla P = \Delta v +F~&\mbox{ in }~\R^2_+\times (0,\infty),\\
&\nabla\cdot v =0~&\mbox{ in }~\R^2_+\times (0,\infty),\\
&\partial_{x_2} v_1\Big|_{x_2=0}=0,\quad v_2\Big|_{x_2=0}=0,\\
&v\big|_{t=0}=0,\\
\end{aligned}
\right.
\end{equation}
where $F$ is solenoidal:
$$
\nabla\cdot F=0,\quad F_2\big|_{x_2=0}=0.
$$
For the system \eqref{eqn-Stokes-Navier}, we derive the Green's tensor and its associated pressure tensor and further obtain the following pointwise upper bounds:
\begin{theorem}\label{thm2}
The solution to \eqref{eqn-Stokes-Navier} with solenoidal forcing can be expressed in the form
\begin{equation*}
\begin{aligned}
&v(x,t)=\int_0^t \int_{\R^2_+} \mathcal G^0(x,y,t-\tau)F(y,\tau)  dyd\tau+\int_0^t\int_{\R^2_+} \mathcal G^*(x,y,t-\tau)\int_0^{\tau}F(y,s)ds dyd\tau\\
&P(x,t)=\int_0^t \int_{\R^2_+} \mathcal P(x,y,t-\tau)\cdot  F(y,\tau) dyd\tau\\
\end{aligned}
\end{equation*}
with $\mathcal G^0=(G_{ij}^0)_{i,j=1,2}$, $\mathcal G^*=(G_{ij}^*)_{i,j=1,2}$, $\mathcal P=(P_j)_{j=1,2}$
\begin{equation*}
\begin{aligned}
G^0_{ij}(x,y,t)=&~\delta_{ij}(\Gamma(x-y,t)-\Gamma(x-y^*,t))\\
G^*_{ij}(x,y,t)=&~(1-\delta_{ij})\frac{\partial}{\partial x_1}\Bigg[-2\frac{\partial\Gamma(x-y^*,t)}{\partial x_2}-4\int_{\R}\frac{\partial E(x_1-z_1,x_2)}{\partial x_2}\frac{\partial\Gamma(z_1-y_1,y_2,t)}{\partial y_2}dz_1\Bigg]\\
&~+\delta_{ij}\frac{\partial}{\partial x_2}\Bigg[-2\frac{\partial\Gamma(x-y^*,t)}{\partial x_2}-4\int_{\R}\frac{\partial E(x_1-z_1,x_2)}{\partial x_2}\frac{\partial\Gamma(z_1-y_1,y_2,t)}{\partial y_2}dz_1\Bigg],\\
P_j(x,y,t)= &~4(1-\delta_{j2}) \frac{\partial}{\partial x_j}\Bigg[\int_{\R} E(x_1-z_1,x_2)\frac{\partial \Gamma(z_1-y_1,y_2,t)}{\partial y_2}dz_1\Bigg].
\end{aligned}
\end{equation*}
Moreover, the following pointwise upper bounds hold
\begin{equation*}
\begin{aligned}
&|\partial_t^sD_x^kD_y^mP_j(x,y,t)|\lesssim t^{-1-s-\frac{m_2}{2}}(|x-y^*|^2+t)^{-\frac{1+|k|+|m'|}{2}}e^{-\frac{cy_2^2}{t}},\\
&|\partial_t^sD_x^kD_y^mG^*_{ij}(x,y,t)|\lesssim t^{-1-s-\frac{m_2}{2}}(|x-y^*|^2+t)^{-\frac{2+|k|+|m'|}{2}}e^{-\frac{cy_2^2}{t}}.
\end{aligned}
\end{equation*}
\end{theorem}

\medskip

As far as we know, above explicit representation formulae and pointwise estimates are not present anywhere and our construction requires rather precise pointwise estimates of the velocity field, so we include those here for self-containedness. The proof of Theorem \ref{thm2} is in a similar spirit as the works by Solonnikov, see for example \cite{Solonnikov1976}. In fact, another way to deal with the forced Navier-Stokes equation is to use the symmetry encoded in the partially free boundary conditions, which simplifies the analysis. More precisely, under certain reflections thanks to the partially free boundary conditions \eqref{FB}, the structure of the full system \eqref{LCF} is preserved, and thus the partially free boundary problem can be regarded as an ``interior'' problem across the boundary. Essentially, this reflection technique shares similarities with the classical Agmon-Douglis-Nirenberg theory (see \cite{ADN1,ADN2}).

\medskip

The key strategy in the proof of Theorem \ref{thm1} is in similar spirit as that of \cite{LCF2D}, namely, one starts first from the harmonic map heat flow and regards the transported term $v\cdot\nabla u$ as a perturbation. Then the leading part of the solution $u$ to the harmonic map heat flow provides external forcing to the Navier-Stokes equation. The velocity, which carries information of $u$, enters harmonic map heat flow in the form of a transported term. Finally, this loop argument is closed once one shows that the transported term is indeed a perturbation. The transported harmonic map heat flow and the incompressible Navier-Stokes equation with forcing are in fact strongly coupled as one can see from the natural scaling invariance, and the smallness of the universal coupling constant $\varepsilon_0$ is to ensure that the full system is less coupled and the loop argument can be implemented. One interesting feature is the natural symmetry/reflection encoded in the partially free boundary conditions \eqref{FB}, which has also been used in \cite{ChenLin98JGA}, and is crucial in our construction, especially for the behaviors for both the orientation field and the velocity field near the boundary. In fact, since the inner concentration zone of each boundary bubble touches the boundary, the use of reflection greatly simplifies the analysis of the linearization, and to be more precise, one can regard both the linearized harmonic map heat flow and the Stokes system near boundary bubble as interior problems. 

\medskip

The construction is based on recently developed {\it parabolic gluing method}, which has been successfully applied in the studies of singularity formation in parabolic equations and systems, geometric flows, fluid equations and others. See for example \cite{Green16JEMS,17HMF,173D,17halfHMF,18Euler} and the references therein.

\medskip

The rest of the paper is devoted to the proofs of the above mentioned  theorems. For simplicity, we construct the most representative case  of one interior bubble and one boundary bubble. The construction of any finitely linear combination of bubbles either in the interior or on the boundary is similar. See Appendix \ref{k-bubbles} for detailed discussions. Before carrying out the rigorous constructions, we first give a brief roadmap and introduce the key ingredients in next subsection.

\bigskip

\subsection*{Roadmap to the construction}

\bigskip

The starting point of the construction is the symmetry (see Section \ref{symmetry} for details) encoded in the partially free boundary conditions \eqref{FB}. The free boundary conditions not only guarantee the energy dissipation but also suggest the correct ansatz for the bubbling, which is crucial for the linearization around the boundary bubble. We should first note that the linearization for the interior bubbling is completely ``localized'' because of our inner--outer construction. In other words, the linearization near interior bubbles does not touch the boundary. There are two crucial aspects needed to be analyzed carefully in the construction:
\begin{itemize}
\item how the partially free boundary conditions \eqref{FB} affect the boundary linearization;
\item how the interior bubble(s) interact with the boundary bubble(s).
\end{itemize}

\medskip

In the case of the half space $\R^2_+$, the partially free boundary conditions \eqref{FB} can be expressed as
\begin{equation}\label{FB1}
\begin{cases}
\partial_{x_2} u_1=0,\\
\partial_{x_2} u_2=0,\\
u_3=0,\qquad\qquad\mbox{ on } \partial\R^2_+,\\
\partial_{x_2} v_1=0,\\
v_2=0,\\
\end{cases}
\end{equation}
so one expects certain symmetry across $\partial\R^2_+$. See the discussions in Section \ref{symmetry} for the preservation of the structure for the full system under the reflection class \eqref{symmetryclass}, and this in turn suggests  to take as a first approximation (for the map)
\begin{equation*}
\begin{aligned}
U_*(x,t)
=U^{(1)}(x,t)+U^{(2)}(x,t)-U^{(2)*}(x,t)
\end{aligned}
\end{equation*}
as in \eqref{def-U*}, where $U^{(1)}$ and $U^{(2)}$ are bubbles placed on the boundary and in the interior, respectively. Here the purpose of the ``reflected'' bubble $U^{(2)*}$ is to make the error in the right symmetry class \eqref{symmetryclass} across the boundary $\partial\R^2_+$; which crucially simplifies our analysis for the linearization around boundary bubble. Indeed, careful analysis for the linearization in Fourier expansion suggests that the inner problem touching boundary can in fact be regarded as an interior problem after proper reflections across $\partial\R^2_+$, provided the right hand side in the linearization has no projection onto certain direction on the tangent plane of boundary bubble (see Section \ref{sec-lineartheories} for more details). One role that $U^{(2)*}$ plays is to enforce such symmetry. In the gluing construction, another important role that the partially free boundary conditions \eqref{FB1} play is to also ensure that the outer ``noises'' coupled into the linearization on the boundary do not destroy the structure discussed above. 

\medskip

The next step is to find perturbation consisting of inner and outer profiles such that a real solution with desired asymptotics can be found. We look for solution in the form
$$
u\sim U_*+u^{(1)}_{\rm inner}+u^{(2)}_{\rm inner}+u_{\rm outer},
$$
where the inner parts $u^{(1)}_{\rm inner}$, $u^{(2)}_{\rm inner}$, expanding on corresponding tangent plane, solve the linearization around the interior bubble and boundary bubble, respectively, and the outer part $u_{\rm outer}$, solving essentially a non-homogeneous heat equation, handles the external noises. Above ansatz then leads to a inner--outer gluing system for $(u^{(1)}_{\rm inner},u^{(2)}_{\rm inner},u_{\rm outer})$. Here $u_{\rm outer}$ is relatively straightforward to solve, while in order to find well-behaved $(u^{(1)}_{\rm inner},u^{(2)}_{\rm inner})$, careful adjustment of modulation parameters is required so that certain orthogonalities are satisfied. The adjustment determines the right dynamics, in particular for the scaling parameters.

\medskip

Dealing with the velocity field $v$ requires the analysis of the Stokes operator with Navier boundary conditions. A direct way is to use its associated Green's tensor derived in Appendix \ref{sec-SO} to capture precise pointwise control. In a similar spirit as the ADN theory as well as Solonnikov's theory, the use of reflections thanks to the partially free boundary \eqref{FB1} in fact reduces the problem into an interior one. 
 
 \medskip
 
 With all the analysis for linearized harmonic map heat flow and Stokes system set up, one can start the loop argument described in the picture below, and the construction is done in appropriate weighted topologies by fixed point arguments.
 
 \medskip

In conclusion, the partially free boundary conditions \eqref{FB} not only imply a natural and physical model of hydrodynamics of nematic liquid crystals, but also encode good structure that triggers the new boundary bubbling phenomenon. With such structure, we then carry out the iterative scheme as in \cite{LCF2D} and close the loop by using $0<\varepsilon_0\ll 1$ in the refined perturbation argument depicted as follows
\catcode`\@=11
\newdimen\cdsep
\cdsep=3em

\def\cdstrut{\vrule height .5\cdsep width 0pt depth .3\cdsep}
\def\@cdstrut{{\advance\cdsep by 2em\cdstrut}}

\def\arrow#1#2{
  \ifx d#1
    \llap{$\scriptstyle#2$}\left\downarrow\cdstrut\right.\@cdstrut\fi
  \ifx u#1
    \llap{$\scriptstyle#2$}\left\uparrow\cdstrut\right.\@cdstrut\fi
  \ifx r#1
    \mathop{\hbox to \cdsep{\rightarrowfill}}\limits^{#2}\fi
  \ifx l#1
    \mathop{\hbox to \cdsep{\leftarrowfill}}\limits^{#2}\fi
}
\catcode`\@=12

\cdsep=3em
$$
\begin{matrix}
\mbox{\eqref{LCF}$_1$:}  -\varepsilon_0 \nabla\cdot \left(\nabla u\odot \nabla u-\frac12 |\nabla u|^2 \mathbb I_2\right)                   & \arrow{r}{{\text{Green's tensor}}}_{\text{Reflection}}  & v ~\mbox{ in \eqref{LCF}$_3$}                    \cr
  \arrow{u}{{\text{Strongly coupled}}} &                      & \arrow{d}{} \cr
\mbox{ Mode } k \mbox{ in the inner problem of $u$: } \phi_k                  & \arrow{l}{{\text{Same asymptotics as the RHS}}}_{0<\varepsilon_0\ll 1} &  v\cdot \nabla u \mbox{ in \eqref{LCF}$_3$ }                  \cr
\end{matrix}
$$

\medskip

\bigskip

\section{Reflection and symmetry of the full system}\label{symmetry}

\bigskip

Recall that the partially free boundary conditions \eqref{FB} in the case $\Omega=\R^2_+$ can be written as 
$$
\partial_{x_2} u_1=\partial_{x_2} u_2=u_3=\partial_{x_2} v_1=v_2=0 ~\mbox{ on }~\partial\R^2_+.\\
$$
So we perform even reflection for $u_1$, $u_2$, $v_1$ and odd reflection for $u_3$, $v_2$:
\begin{equation*}
\tilde u(x_1,x_2,t)=\begin{bmatrix}
u_1(x_1,-x_2,t)\\
u_2(x_1,-x_2,t)\\
-u_3(x_1,-x_2,t)\\
\end{bmatrix}, \quad
\tilde v(x_1,x_2,t)=\begin{bmatrix}
v_1(x_1,-x_2,t)\\
-v_2(x_1,-x_2,t)\\
\end{bmatrix}, \quad
x_2<0
\end{equation*}
such that the partially free boundary conditions are automatically satisfied. Note that this reflection technique has already  been used in the harmonic map heat flow with partially free boundary (see \cite{ChenLin98JGA} for example). By
\begin{equation*}
\nabla\cdot( \nabla u \odot \nabla u)=\begin{bmatrix}
2\partial_1 u_k \partial_{11}u_k+\partial_{22}u_k\partial_1 u_k+\partial_{21}u_k\partial_2 u_k\\
2\partial_2 u_k \partial_{22}u_k+\partial_{11}u_k\partial_2 u_k+\partial_{12}u_k\partial_1 u_k\\
\end{bmatrix}
\end{equation*}
and
\begin{equation*}
v\cdot\nabla u=\begin{bmatrix}
v_1\partial_1 u_1+v_2\partial_2 u_1\\
v_1\partial_1 u_2+v_2\partial_2 u_2\\
v_1\partial_1 u_3+v_2\partial_2 u_3\\
\end{bmatrix},
\end{equation*}
it is direct to check that
$$
\partial_t \tilde u+\tilde v\cdot\nabla \tilde u=\Delta \tilde u +|\nabla \tilde u|^2 \tilde u,
$$
and
$$
\partial_t \tilde v +\tilde v\cdot \nabla \tilde v +\nabla \tilde P = \Delta \tilde v - \varepsilon_0\nabla\cdot \left(\nabla \tilde u \odot \nabla \tilde u-\frac12 |\nabla \tilde u|^2 \mathbb{I}_2\right)
$$
if $\tilde P(x_1,x_2,t)=P(x_1,-x_2,t)$ for $x_2<0$, which means $\frac{\partial P}{\partial \nu}=0$. In fact, this is true with our partially free boundary conditions \eqref{FB}. Indeed, applying the outer normal $\nu$ to \eqref{LCF}$_1$, we get
$$
\partial_t v_2 + v_1\partial_1 v_2 + v_2 \partial_2 v_2 + \partial_2 P=\Delta v_2- \varepsilon_0(\partial_2 u_k \partial_{22}u_k+\partial_{11}u_k\partial_2 u_k).
$$
By the divergence free condition and the Navier slip boundary condition, we have
$$
\partial_{22} v_1=-\partial_1(\partial_2 v_1)=0 ~\mbox{ on }~ \partial \R^2_+.
$$
Similarly, on $\partial \R^2_+$, by the partially free boundary condition
\begin{equation*}
\begin{aligned}
\partial_2 u_k \partial_{22}u_k+\partial_{11}u_k\partial_2 u_k=&~\partial_2 u_3 \partial_{22} u_3 + \partial_{11}u_3 \partial_2 u_3\\
=&~\partial_2 u_3 \Delta u_3\\
=&~\partial_2 u_3 (\partial_t u_3 + v_1\partial_1 u_3+v_2 \partial_2 u_3-|\nabla u|^2 u_3)\\
=&~0.
\end{aligned}
\end{equation*}
Therefore, we have $\frac{\partial P}{\partial \nu}=0$ on $\partial \R^2_+$.

\medskip

In conclusion, with the following reflections
\begin{equation*}
\tilde u(x_1,x_2,t)=\begin{bmatrix}
u_1(x_1,-x_2,t)\\
u_2(x_1,-x_2,t)\\
-u_3(x_1,-x_2,t)\\
\end{bmatrix}, \quad
\tilde v(x_1,x_2,t)=\begin{bmatrix}
v_1(x_1,-x_2,t)\\
-v_2(x_1,-x_2,t)\\
\end{bmatrix}, \quad
\tilde P(x_1,x_2,t)=P(x_1,-x_2,t),
\end{equation*}
the structure of the equation \eqref{LCF} is preserved, i.e.,
\begin{equation*}
\begin{cases}
\partial_t \tilde v +\tilde v\cdot \nabla \tilde v +\nabla \tilde P = \Delta \tilde v - \varepsilon_0\nabla\cdot \left(\nabla \tilde u \odot \nabla \tilde u-\frac12 |\nabla \tilde u|^2 \mathbb{I}_2\right),\\
\nabla\cdot \tilde v=0,\\
\partial_t \tilde u+\tilde v\cdot\nabla \tilde u=\Delta \tilde u +|\nabla \tilde u|^2 \tilde u.\\
\end{cases}
\end{equation*}
We shall look for a solution to \eqref{LCF} with partially free boundary condition \eqref{FB} in the symmetry class across the boundary $\partial \R^2_+$ given by:
\begin{equation}\label{symmetryclass}
u_1,~u_2,~v_1,~P \mbox{ are even in } x_2 \mbox{ and } u_3,~v_2 \mbox{ are odd in } x_2
\end{equation}
with $u=(u_1,u_2,u_3)^T$, $v=(v_1,v_2)^T$.

\begin{remark}
It is worth mentioning that if one imposes the no-slip boundary condition for the velocity
\begin{equation*}
v\big|_{\partial \R^2_+}=0
\end{equation*}
instead of the Navier boundary conditions \eqref{FB}$_1$--\eqref{FB}$_2$, the natural energy dissipation is also preserved. However, this artificial boundary may destroy the structure of the coupled system in nature, as one can see in the reflections, i.e., there is no reflection preserving the structure of the entire system.
\end{remark}

\medskip

Next, we try to gain some information from the symmetry at the linearized level. We consider the infinitesimal generator of rigid motions: dilation, translations and rotations. More precisely, the invariance from scaling and rotation around $z$-axis corresponds to mode $0$, the invariance from translations corresponds to mode $1$, and the invariance from rotations around $x$ and $y$ axes corresponds to mode $-1$. 

\medskip

Let $W_2$ be the least energy (degree 1) harmonic map
\begin{equation}\label{def-W2}
W_2(x)=\begin{bmatrix}
\frac{2x}{1+|x|^2}\\
\frac{|x|^2-1}{1+|x|^2}\\
\end{bmatrix},~x\in\R^2,
\end{equation}
namely, $\int_{\R^2}|\nabla W_2|^2=8\pi$. Our first approximation of the boundary bubble will be based on the degree 1 profile
\begin{equation}\label{def-W1}
W_1:=Q_*W_2.
\end{equation}
Here we introduce
\begin{equation}\label{def-Q_*}
Q_*=
\begin{bmatrix}
1&0&0\\
0&0&1\\
0&1&0\\
\end{bmatrix}
\end{equation}
because of the reflection \eqref{symmetryclass}.

\medskip

In fact, there is some ``rigidity'' produced by the partially free boundary conditions, especially for the boundary bubble(s). To see this, we formally compute below the first variations with respect to different parameters. For the rotations around $x$ and $y$ axes (mode $-1$), we consider

$$
W_{1,\alpha}:=\begin{bmatrix}
       1 & 0 & 0 \\[0.3em]
       0 & \cos\alpha & -\sin\alpha \\[0.3em]
       0 & \sin\alpha & \cos\alpha\\
     \end{bmatrix}
     \begin{bmatrix}
     \frac{2x_1}{x_1^2+x_2^2+1}\\
     \frac{x_1^2+x_2^2-1}{x_1^2+x_2^2+1}\\
     \frac{2x_2}{x_1^2+x_2^2+1}\\
     \end{bmatrix},\quad
     W_{1,\beta}:=\begin{bmatrix}
       \cos\beta & 0 & \sin\beta \\[0.3em]
       0 & 1 & 0 \\[0.3em]
       -\sin\beta & 0 & \cos\beta\\
     \end{bmatrix}
     \begin{bmatrix}
     \frac{2x_1}{x_1^2+x_2^2+1}\\
     \frac{x_1^2+x_2^2-1}{x_1^2+x_2^2+1}\\
     \frac{2x_2}{x_1^2+x_2^2+1}\\
     \end{bmatrix}.
$$
The first variations
$$
\partial_{\alpha} W_{1,\alpha}\big|_{\alpha=0}=\begin{bmatrix}
0\\
-\frac{2x_2}{x_1^2+x_2^2+1}\\
\frac{x_1^2+x_2^2-1}{x_1^2+x_2^2+1}\\
\end{bmatrix}, \quad \partial_{\beta} W_{1,\beta}\big|_{\beta=0}=\begin{bmatrix}
\frac{2x_2}{x_1^2+x_2^2+1}\\
0\\
-\frac{2x_1}{x_1^2+x_2^2+1}\\
\end{bmatrix}
$$
are not in the symmetry class \eqref{symmetryclass}.

\medskip

Similarly, for the scaling and rotation around $z$ axis (mode $0$)
$$
W_{1,\lambda}:=\begin{bmatrix}
\frac{2\la x_1}{x_1^2+x_2^2+\la^2}\\
\frac{x_1^2+x_2^2-\la^2}{x_1^2+x_2^2+\la^2}\\
\frac{2\la x_2}{x_1^2+x_2^2+\la^2}\\
\end{bmatrix}, \quad
W_{1,\omega}:=\begin{bmatrix}
       \cos\omega & -\sin\omega & 0 \\[0.3em]
       \sin \omega & \cos\omega  & 0 \\[0.3em]
       0 & 0 & 1\\
     \end{bmatrix} \begin{bmatrix}
     \frac{2x_1}{x_1^2+x_2^2+1}\\
     \frac{x_1^2+x_2^2-1}{x_1^2+x_2^2+1}\\
     \frac{2x_2}{x_1^2+x_2^2+1}\\
     \end{bmatrix},
$$
one has
$$
\partial_{\la} W_{1,\la}\big|_{\la=1}=\begin{bmatrix}
\frac{2x_1(x_1^2+x_2^2-1)}{(x_1^2+x_2^2+1)^2}\\
-\frac{4(x_1^2+x_2^2)}{(x_1^2+x_2^2+1)^2}\\
\frac{2x_2(x_1^2+x_2^2-1)}{(x_1^2+x_2^2+1)^2}\\
\end{bmatrix}, \quad
\partial_{\omega} W_{1,\omega}\big|_{\omega=0}=\begin{bmatrix}
     -\frac{x_1^2+x_2^2-1}{x_1^2+x_2^2+1}\\
     \frac{2x_1}{x_1^2+x_2^2+1}\\
     0\\
     \end{bmatrix},
$$
which are in the symmetry class.

\medskip

For the translations (mode $1$)
$$
W_{1,\xi_1}:=\begin{bmatrix}
     \frac{2(x_1-\xi_1)}{(x_1-\xi_1)^2+x_2^2+1}\\
     \frac{(x_1-\xi_1)^2+x_2^2-1}{(x_1-\xi_1)^2+x_2^2+1}\\
     \frac{2x_2}{(x_1-\xi_1)^2+x_2^2+1}\\
     \end{bmatrix}, \quad
W_{1,\xi_2}:=\begin{bmatrix}
     \frac{2x_1}{x_1^2+(x_2-\xi_2)^2+1}\\
     \frac{x_1^2+(x_2-\xi_2)^2-1}{x_1^2+(x_2-\xi_2)^2+1}\\
     \frac{2(x_2-\xi_2)}{x_1^2+(x_2-\xi_2)^2+1}\\
     \end{bmatrix},
$$
we have
$$
\partial_{\xi_1} W_{1,\xi_1}\big|_{\xi_1=0}=\begin{bmatrix}
-\frac{2(x_2^2-x_1^2+1)}{(x_1^2+x_2^2+1)^2}\\
-\frac{4x_1}{(x_1^2+x_2^2+1)^2}\\
\frac{4x_1x_2}{(x_1^2+x_2^2+1)^2}\\
\end{bmatrix}, \quad
\partial_{\xi_2} W_{1,\xi_2}\big|_{\xi_2=0}=\begin{bmatrix}
\frac{4x_1x_2}{(x_1^2+x_2^2+1)^2}\\
-\frac{4x_2}{(x_1^2+x_2^2+1)^2}\\
-\frac{2(x_1^2-x_2^2+1)}{(x_1^2+x_2^2+1)^2}\\
\end{bmatrix}.
$$
Only the first one is in the symmetry class. In other words, the symmetry only allows us to translate the bubble along the boundary $\partial \R^2_+$.

\medskip

Heuristically, the above argument suggests certain rigidity produced by the partially free boundary conditions, which reflects in the presence of modulation parameters.

\bigskip

\medskip

\section{Notations and preliminaries}

\medskip

In order to analyze the transported harmonic map heat flow
\begin{equation*}
u_t+v\cdot\nabla u=\Delta u+|\nabla u|^2 u ~\mbox{ in }\R^2_+
\end{equation*}
in the symmetry class \eqref{symmetryclass} across $\partial\R^2_+$, we first regard the transported term $v\cdot\nabla u$ as a perturbation and introduce our first approximation and its correction to the harmonic map heat flow. Later, after the analysis of the forced incompressible Navier-Stokes equation, we will show that the transported term is indeed a perturbation. We first give some useful notations and formulae.

\medskip

Recall that we will construct a bubbling solution which blows up at a given boundary point and interior point. In the following, the superscripts ``$(1)$'', ``$(2)$'' refer to the bubble placed on the boundary and in the interior respectively; and we will repeatedly adopt this notation to distinguish these two bubbles and their associated tangent planes.

Introduce polar coordinates near two concentration zones
$$
y_j=\frac{x-\xi^{(j)}}{\la^{(j)}},\quad x=\xi^{(j)}+\la^{(j)}\rho_j e^{i\theta_j},\quad r_j=\la^{(j)}\rho_j,\quad j=1,2.
$$
In polar coordinates $y_j=\rho_j e^{i\theta_j}$, $j=1,2$, the least energy harmonic maps $W_1(y_1)$, $W_2(y_2)$, defined in \eqref{def-W1} and \eqref{def-W2} respectively, can be represented as

\begin{equation*}
W_1(y_1)=\begin{bmatrix}
\cos\theta_1 \sin w(\rho_1)\\
\cos w(\rho_1)\\
\sin\theta_1 \sin w(\rho_1)\\
\end{bmatrix}
,\quad W_2(y_2)=\begin{bmatrix}
\cos\theta_2 \sin w(\rho_2)\\
\sin\theta_2 \sin w(\rho_2)\\
\cos w(\rho_2)\\
\end{bmatrix}:=\begin{bmatrix}
e^{i\theta_2} \sin w(\rho_2)\\
\cos w(\rho_2)\\
\end{bmatrix}
\end{equation*}
with
$$
w(\rho_j)=\pi-2 \arctan (\rho_j),
$$
and we have
\begin{equation*}
w_{\rho_j}=-\frac{2}{\rho_j^2+1},~~\sin w(\rho_j)=-\rho_j w_{\rho_j}=\frac{2\rho_j}{\rho_j^2+1},~~\cos w(\rho_j)=\frac{\rho_j^2-1}{\rho_j^2+1}.
\end{equation*}
The linearization of the harmonic map operator around $W_j$ is the elliptic operator
\begin{equation}\label{def-linearization}
L_{W}^{(j)}[\phi]:=\Delta_{y_j} \phi + |\nabla W_j(y_j)|^2 \phi +2(\nabla W_j(y_j)\cdot\nabla\phi) W_j(y_j),
\end{equation}
whose kernel functions are given by
\begin{equation}\label{kernels}
\left\{
\begin{aligned}
&Z_{0,1}^{(j)}(y_j)=\rho_j w_{\rho_j}(\rho_j) E_1^{(j)}(y_j),\\
&Z_{0,2}^{(j)}(y_j)=\rho_j w_{\rho_j}(\rho_j) E_2^{(j)}(y_j),\\
&Z_{1,1}^{(j)}(y_j)=w_{\rho_j}(\rho_j)[\cos \theta_j E_1^{(j)}(y_j)+\sin\theta_j E_2^{(j)}(y_j)],\\
&Z_{1,2}^{(j)}(y_j)=w_{\rho_j}(\rho_j)[\sin \theta_j E_1^{(j)}(y_j)-\cos\theta_j E_2^{(j)}(y_j)],\\
&Z_{-1,1}^{(j)}(y_j)=\rho_j^2 w_{\rho_j}(\rho_j)[\cos \theta_j E_1^{(j)}(y_j)-\sin \theta_j E_2^{(j)}(y_j)],\\
&Z_{-1,2}^{(j)}(y_j)=\rho_j^2 w_{\rho_j}(\rho_j)[\sin \theta_j E_1^{(j)}(y_j)+\cos \theta_j E_2^{(j)}(y_j)],\\
\end{aligned}
\right.
\end{equation}
where the vectors
\begin{equation}\label{def-E1E2}
\begin{aligned}
E_1^{(1)}(y_1)=\begin{bmatrix}\cos\theta_1\cos w(\rho_1)\\ -\sin w(\rho_1)\\\sin\theta_1 \cos w(\rho_1)\\ \end{bmatrix},~~E_2^{(1)}(y_1)=\begin{bmatrix} -\sin\theta_1\\ 0\\ \cos\theta_1\\ \end{bmatrix},\\
E_1^{(2)}(y_2)=\begin{bmatrix}\cos\theta_2\cos w(\rho_2)\\ \sin\theta_2 \cos w(\rho_2)\\-\sin w(\rho_2)\\ \end{bmatrix},~~E_2^{(2)}(y_2)=\begin{bmatrix} -\sin\theta_2\\  \cos\theta_2\\ 0\\ \end{bmatrix}\\
\end{aligned}
\end{equation}
form an orthonormal basis of the tangent space $T_{W_j(y_j)} \mathbb{S}^2$, i.e., Frenet basis associated to $W_j$. We see that
$$L^{(j)}_{W}[Z^{(j)}_{p,q}]=0 ~~ {~\rm{ for }} ~ p=\pm 1,~0, ~j,~q=1,2.$$
Because of the scaling, rotation and translation symmetries,
\begin{equation}\label{def-U1U2}
U^{(1)}(x,t):=W_1\left(\frac{x-\xi^{(1)}}{\la^{(1)}}\right),\quad U^{(2)}(x,t):=Q_{\omega}W_2\left(\frac{x-\xi^{(2)}}{\la^{(2)}}\right)
\end{equation}
solve the harmonic map equation, where $Q_{\omega}$ is the rotation of angle $\omega$ matrix around $z$-axis (viewing the target $\mathbb S^2$ embedded into $\R^3$) 
\begin{equation}\label{def-Qz}
Q_{\omega}:=
\begin{bmatrix}
       \cos\omega & -\sin\omega & 0 \\[0.3em]
       \sin \omega & \cos\omega  & 0 \\[0.3em]
       0 & 0 & 1\\
     \end{bmatrix}.
\end{equation}
Notice that
\begin{align*}
&~L^{(1)}_U[\varphi^{(1)}]=(\la^{(1)})^{-2} L^{(1)}_W[\phi^{(1)}],\quad L^{(2)}_U[\varphi^{(2)}]=(\la^{(2)})^{-2} Q_{\omega} L^{(2)}_W[\phi^{(2)}],\,\,\text{where}\\
&~\varphi^{(1)}(x)=\phi^{(1)}(y_1),\quad \varphi^{(2)}(x)=Q_{\omega}\phi^{(2)}(y_2),\quad y_j=\frac{x-\xi^{(j)}}{\la^{(j)}},
\end{align*}
and $L_U^{(j)}$ stands for the linearization around $U^{(j)}$, $j=1,\,2$.
In the sequel, it is of significance to compute the action of $L^{(j)}_U$ on functions whose values are orthogonal to $U^{(j)}$ pointwise. Define
\begin{equation}\label{def-proj}
\Pi_{U^{\perp}}^{(j)}\varphi:=\varphi-(\varphi\cdot U^{(j)})U^{(j)}.
\end{equation}
We now give several useful formulae whose proof is similar to that of \cite[Section 3]{17HMF}:
\begin{equation*}
L^{(j)}_U[\Pi^{(j)}_{U^{\perp}}\Phi]=\Pi^{(j)}_{U^{\perp}}\Delta \Phi +\tilde L^{(j)}_{U}[\Phi],
\end{equation*}
where we denote
\begin{equation}\label{def-tildeL}
\tilde L^{(j)}_{U}[\Phi]:=|\nabla U^{(j)}|^2\Pi^{(j)}_{U^{\perp}}\Phi-2\nabla(\Phi\cdot U^{(j)})\nabla U^{(j)},
\end{equation}
with
$$\nabla(\Phi\cdot U^{(j)})\nabla U^{(j)}=\sum_{k=1}^2\partial_{x_k} (\Phi\cdot U^{(j)}) \partial_{x_k} U^{(j)}, \quad x=(x_1,x_2).$$

We give several useful expressions of the operator $\tilde L^{(j)}_{U}$ acting on $\Phi$ in different forms:
\begin{itemize}
\item In the polar coordinates
$$\Phi(x)=\Phi(r_j,\theta_j),~~x=\xi^{(j)}+r_j e^{i\theta_j},$$
the operator \eqref{def-tildeL} can be expressed as
\begin{equation*}
\begin{aligned}
\tilde L_U^{(1)}[\Phi]&=-\frac{2}{\la^{(1)}} w_{\rho_1}(\rho_1)\left[(\partial_{r_1}\Phi\cdot U^{(1)})  E^{(1)}_1 -\frac1{r_1} (\partial_{\theta_1}\Phi\cdot U^{(1)})  E_2^{(1)}\right],~~r_1=\la^{(1)} \rho_1,\\
\tilde L_U^{(2)}[\Phi]&=-\frac{2}{\la^{(2)}} w_{\rho_2}(\rho_2)\left[(\partial_{r_2}\Phi\cdot U^{(2)})  Q_{\omega}E^{(2)}_1 -\frac1{r_2} (\partial_{\theta_2}\Phi\cdot U^{(2)}) Q_{\omega} E_2^{(2)}\right],~~r_2=\la^{(2)} \rho_2.\\
\end{aligned}
\end{equation*}
\medskip
\item For the operator $\tilde L_U^{(1)}$ acting on $\Phi(x)=(\varphi_1(x),\varphi_2(x),\varphi_3(x))^T$, we can compute
\begin{equation*}
\begin{aligned}
\partial_{r_1}\Phi\cdot U^{(1)}=&~\left((\cos\theta_1 \partial_{x_1}+\sin\theta_1 \partial_{x_2})\begin{bmatrix}
\varphi_1\\
\varphi_2\\
\varphi_3\\
\end{bmatrix}\right)\cdot \begin{bmatrix}
\cos\theta_1 \sin w(\rho_1)\\
\cos w(\rho_1)\\
\sin\theta_1 \sin w(\rho_1)\\
\end{bmatrix}\\
=&~\frac{\sin w(\rho_1)}{2}\bigg[(\partial_{x_1}\varphi_1+\partial_{x_2}\varphi_3)+\sin2\theta_1 (\partial_{x_2}\varphi_1+\partial_{x_1}\varphi_3)\\
&~+\cos2\theta_1 (\partial_{x_1}\varphi_1-\partial_{x_2}\varphi_3)
\bigg]+\cos w(\rho_1)[\cos\theta_1 \partial_{x_1}\varphi_2+\sin\theta_1 \partial_{x_2}\varphi_2],
\end{aligned}
\end{equation*}
\begin{equation*}
\begin{aligned}
\frac1{r_1}\partial_{\theta_1}\Phi\cdot U^{(1)}=&~\left((\cos\theta_1 \partial_{x_2}-\sin\theta_1 \partial_{x_1})\begin{bmatrix}
\varphi_1\\
\varphi_2\\
\varphi_3\\
\end{bmatrix}\right)\cdot \begin{bmatrix}
\cos\theta_1 \sin w(\rho_1)\\
\cos w(\rho_1)\\
\sin\theta_1 \sin w(\rho_1)\\
\end{bmatrix}\\
=&~\frac{\sin w(\rho_1)}{2}\bigg[(\partial_{x_2}\varphi_1-\partial_{x_1}\varphi_3)+\sin2\theta_1 (\partial_{x_2}\varphi_3-\partial_{x_1}\varphi_1)\\
&~+\cos2\theta_1 (\partial_{x_2}\varphi_1+\partial_{x_1}\varphi_3)
\bigg]+\cos w(\rho_1)[\cos\theta_1 \partial_{x_2}\varphi_2-\sin\theta_1 \partial_{x_1}\varphi_2],
\end{aligned}
\end{equation*}
and thus
\begin{equation}\label{def-tildeL(1)-1}
\begin{aligned}
\tilde L_U^{(1)}[\Phi]&:=[\tilde L_U]^{(1)}_0[\Phi]+[\tilde L_U]^{(1)}_1[\Phi]+[\tilde L_U]^{(1)}_2[\Phi]
\end{aligned}
\end{equation}
with
\begin{equation}\label{def-tildeL(1)-2}
\left\{
\begin{aligned}
~ [\tilde L_U]^{(1)}_0[\Phi]:=&~ (\la^{(1)})^{-1} \rho_1 w^2_{\rho_1}(\rho_1) \bigg[(\partial_{x_1}\varphi_1+\partial_{x_2}\varphi_3)  E^{(1)}_1+  (\partial_{x_1}\varphi_3-\partial_{x_2}\varphi_1)E^{(1)}_2\bigg],\\
[\tilde L_U]^{(1)}_1[\Phi]:=&~ -2(\la^{(1)})^{-1} w_{\rho_1}(\rho_1)\cos w(\rho_1) \bigg[(\partial_{x_1}\varphi_2)\cos\theta_1+(\partial_{x_2}\varphi_2)\sin\theta_1\bigg]  E_1^{(1)}\\
&~+2(\la^{(1)})^{-1} w_{\rho_1}(\rho_1)\cos w(\rho_1) \bigg[(\partial_{x_2}\varphi_2)\cos\theta_1-(\partial_{x_1}\varphi_2)\sin\theta_1\bigg]  E^{(1)}_2,\\
[\tilde L_U]^{(1)}_2[\Phi]:=&~ (\la^{(1)})^{-1} \rho_1 w^2_{\rho_1}(\rho_1)  \bigg[(\partial_{x_2}\varphi_1+\partial_{x_1}\varphi_3)\sin2\theta_1+(\partial_{x_1}\varphi_1-\partial_{x_2}\varphi_3)\cos2\theta_1\bigg]  E^{(1)}_1\\
&~+(\la^{(1)})^{-1}\rho_1 w^2_{\rho_1}(\rho_1)  \bigg[(\partial_{x_1}\varphi_1-\partial_{x_2}\varphi_3)\sin2\theta_1-(\partial_{x_2}\varphi_1+\partial_{x_1}\varphi_3)\cos2\theta_1\bigg] E_2^{(1)},
\end{aligned}
\right.
\end{equation}
where we have used $\sin w(\rho_1)=-\rho_1 w_{\rho_1}(\rho_1)$.

\medskip

\item Another convenient form is the following: for a $C^1$ function $\Phi(x): \R^2\to \mathbb C\times \R$ written in the complex form
\begin{equation*}
\Phi(x)=(\varphi_1(x),\varphi_2(x),\varphi_3(x))^T
:=\begin{bmatrix}
\varphi_1(x)+i\varphi_2(x)\\
\varphi_3(x)\\
\end{bmatrix},
\end{equation*}
if we write
\begin{equation}\label{def-complexnotations}
\begin{aligned}
&\varphi=\varphi_1 +i\varphi_2,\quad\bar\varphi=\varphi_1 -i\varphi_2,\\
&{\rm div}\varphi=\partial_{x_1} \varphi_1 +\partial_{x_2} \varphi_2,\quad{\rm curl} \varphi=\partial_{x_1} \varphi_2 -\partial_{x_2} \varphi_1,
\end{aligned}
\end{equation}
then we can express
\EQ{\label{def-tildeL(2)-1}
\tilde L^{(2)}_U[\Phi]:=[\tilde L_U]^{(2)}_0[\Phi]+[\tilde L_U]^{(2)}_1[\Phi]+[\tilde L_U]^{(2)}_2[\Phi]
}
as
\begin{equation}\label{def-tildeL(2)-2}
\left\{
\begin{aligned}
~ [\tilde L_U]^{(2)}_0[\Phi]:=&~ (\la^{(2)})^{-1} \rho_2 w^2_{\rho_2}(\rho_2) [{\rm div}(e^{-i\omega}\varphi) Q_{\omega} E^{(2)}_1+{\rm curl}(e^{-i\omega}\varphi) Q_{\omega} E^{(2)}_2],\\
[\tilde L_U]^{(2)}_1[\Phi]:=&~ -2(\la^{(2)})^{-1} w_{\rho_2}(\rho_2)\cos w(\rho_2) [(\partial_{x_1}\varphi_3)\cos\theta_2+(\partial_{x_2}\varphi_3)\sin\theta_2] Q_{\omega} E_1^{(2)}\\
&~-2(\la^{(2)})^{-1} w_{\rho_2}(\rho_2)\cos w(\rho_2) [(\partial_{x_1}\varphi_3)\sin\theta_2-(\partial_{x_2}\varphi_3)\cos\theta_2] Q_{\omega} E^{(2)}_2,\\
[\tilde L_U]^{(2)}_2[\Phi]:=&~ (\la^{(2)})^{-1} \rho_2 w^2_{\rho_2}(\rho_2)  [{\rm div}(e^{i\omega}\bar\varphi)\cos 2\theta_2-{\rm curl}(e^{i\omega}\bar\varphi)\sin 2\theta_2] Q_{\omega} E^{(2)}_1\\
&~+(\la^{(2)})^{-1}\rho_2 w^2_{\rho_2}(\rho_2)  [{\rm div}(e^{i\omega}\bar\varphi)\sin 2\theta_2+{\rm curl}(e^{i\omega}\bar\varphi)\cos 2\theta_2] Q_{\omega} E_2^{(2)}.
\end{aligned}
\right.
\end{equation}
The proof is similar to that of $\tilde L^{(1)}_U$.
\medskip
\item If we assume
\begin{equation*}
\Phi(x)=\begin{bmatrix}
\phi(r_2) e^{i\theta_2}\\
0\\
\end{bmatrix}
,\quad x=\xi^{(2)}+r_2 e^{i\theta_2},\quad r_2=\la^{(2)} \rho_2,
\end{equation*}
where $\phi(r)$ is complex-valued, then we have the following formula
\begin{equation*}
\tilde L^{(2)}_U[\Phi]=\frac{2}{\la^{(2)}} w^2_{\rho_2}(\rho_2) \left[{\rm Re}(e^{-i\omega} \partial_{r_2} \phi(r_2)) Q_{\omega} E^{(2)}_1 + \frac1{r_2} {\rm Im}(e^{-i\omega}\phi(r_2)) Q_{\omega} E^{(2)}_2\right].
\end{equation*}
\end{itemize}

\medskip

\section{Approximation and improvement}\label{Sec-approx}

\medskip

We consider the case of two bubbles with one placed in the interior and the other placed on the boundary. In this section, we introduce the approximate solution and its improvement. We consider the approximation
\begin{equation}\label{def-U*}
\begin{aligned}
U_*(x,t)
:=&~U^{(1)}(x,t)+U^{(2)}(x,t)-U^{(2)*}(x,t), \quad x\in\R^2_+,
\end{aligned}
\end{equation}
where $U^{(1)}$ and $U^{(2)}$ are given in \eqref{def-U1U2} and 
\begin{equation}
U^{(2)*}(x,t):=Q_{\omega}\begin{bmatrix}
-\frac{2\la^{(2)}(x_1-\xi^{(2)}_1)}{|x-\xi^{(2)*}|^2+(\la^{(2)})^2}\\
\frac{2\la^{(2)}(x_2+\xi^{(2)}_2)}{|x-\xi^{(2)*}|^2+(\la^{(2)})^2}\\
\frac{|x-\xi^{(2)*}|^2-(\la^{(2)})^2}{|x-\xi^{(2)*}|^2+(\la^{(2)})^2}\\
\end{bmatrix},
\end{equation}
$Q_{\omega}$ is the rotation matrix defined in \eqref{def-Qz},
and we take
\begin{equation}\label{def-xi}
\begin{aligned}
&\xi^{(1)}(t)=(\xi^{(1)}_{1}(t),0)\in\partial\R_+^2,\\
&\xi^{(2)}(t)=(\xi^{(2)}_{1}(t),\xi^{(2)}_{2}(t)),\quad \xi^{(2)*}(t)=(\xi^{(2)}_{1}(t),-\xi^{(2)}_{2}(t)),\quad \xi^{(2)}_2(t)>0.
\end{aligned}
\end{equation}
The reason for taking the above approximation is twofold: the reflection term is to preserve symmetry, and the leading profile should be approximately of length one. Clearly $U^{(2)*}$ is a least energy harmonic map. Denoting the error operator by 
$$
S(u) :=  -u_t+ \Delta u + |\nabla u|^2 u,
$$ we then have
\begin{equation*}
\begin{aligned}
S(U_*)=&~-\partial_t U^{(1)}-\partial_t (U^{(2)}-U^{(2)*})+|\nabla U_*|^2 U_*\\
&~-|\nabla U^{(1)}|^2U^{(1)}-|\nabla U^{(2)}|^2U^{(2)}+|\nabla U^{(2)*}|^2 U^{(2)*}\\
=&~-\bigg[\underbrace{\dot\la^{(1)} \pd_{\la^{(1)}} U^{(1)}}_{:=\mathcal E_0^{(1)}} +\underbrace{ \dot\xi^{(1)}_1  \pd_{\xi^{(1)}_1} U^{(1)}}_{:=\mathcal E_1^{(1)}} \bigg]\\
&~-\bigg[\underbrace{\dot\la^{(2)} \pd_{\la^{(2)}} (U^{(2)}-U^{(2)*}) + \dot \omega \pd_{\omega} (U^{(2)}-U^{(2)*})}_{:=\mathcal E_0^{(2)}} \\
&~+\underbrace{ \dot\xi^{(2)}_1  \pd_{\xi^{(2)}_1} (U^{(2)}-U^{(2)*})+\dot\xi^{(2)}_2  \pd_{\xi^{(2)}_2} (U^{(2)}-U^{(2)*})}_{:=\mathcal E_1^{(2)}} \bigg]\\
&~+\bigg|\nabla(U^{(1)}+U^{(2)}-U^{(2)*})\bigg|^2\big[U^{(1)}+U^{(2)}-U^{(2)*}\big]\\
&~-|\nabla U^{(1)}|^2U^{(1)}-|\nabla U^{(2)}|^2U^{(2)}+|\nabla U^{(2)*}|^2 U^{(2)*}\\
:=&~\mathcal E_0^{(1)}+\mathcal E_1^{(1)}+\mathcal E_0^{(2)}+\mathcal E_1^{(2)}+\tilde{\mathcal E},
\end{aligned}
\end{equation*}
where
\begin{equation*}
\left\{
\begin{aligned}
&\pd_{\la^{(1)}} U^{(1)}(x)=-(\la^{(1)})^{-1}  Z^{(1)}_{0,1}(y_1),\\
&\pd_{\xi_1^{(1)}} U^{(1)}(x) = (\la^{(1)})^{-1}  Z^{(1)}_{1,1}(y_1),\\
&\pd_{\xi_2^{(1)}} U^{(1)}(x) = (\la^{(1)})^{-1}   Z_{1,2}^{(1)}(y_1),\\
&\pd_{\la^{(2)}} U^{(2)}(x)=-(\la^{(2)})^{-1} Q_{\omega} Z^{(2)}_{0,1}(y_2),\\
&\pd_{\omega} U^{(2)}(x)=  -Q_{\omega} Z_{0,2}^{(2)}(y_2),\\
&\pd_{\xi_1^{(2)}} U^{(2)}(x) = (\la^{(2)})^{-1}  Q_{\omega} Z^{(2)}_{1,1}(y_2),\\
&\pd_{\xi_2^{(2)}} U^{(2)}(x) = (\la^{(2)})^{-1}  Q_{\omega} Z_{1,2}^{(2)}(y_2),\\
\end{aligned}
\right.
\end{equation*}
with
\begin{equation*}
\left\{
\begin{aligned}
&Z^{(j)}_{0,1}(y_j )=\rho_j  w_{\rho_j }(\rho_j ) E^{(j)}_1(y_j ),\\
&Z^{(2)}_{0,2}(y_2 )=\rho_2  w_{\rho_2 }(\rho_2 ) E^{(2)}_2(y_2 ),\\
&Z^{(j)}_{1,1}(y_j )=w_{\rho_j }(\rho_j )[\cos \theta_j  E^{(j)}_1(y_j )+\sin\theta_j  E^{(j)}_2(y_j )],\\
&Z^{(j)}_{1,2}(y_j )=w_{\rho_j }(\rho_j )[\sin \theta_j  E^{(j)}_1(y_j )-\cos\theta_j  E^{(j)}_2(y_j )],\\
\end{aligned}
\right.
\end{equation*}
for $j=1,2.$ Here the definitions of $Z_{p,q}^{(j)}$, $E_1^{(j)}$, $E_2^{(j)}$ can be found in \eqref{kernels}, \eqref{def-E1E2}.

We notice that the error $S(U_*)$ contains slow spatial decaying terms in $\mathcal E_0^{(j)}$, in other words, these terms are not in $L^{2}(\R^2)$, and $\mathcal E_0^{(j)}$ in fact corresponds to mode $0$ of each concentration. To improve the spatial decay, we add two global corrections $\Phi_0^{(j)}$ solving at leading order
\begin{equation*}
\partial_t \Phi_0^{(j)}\approx\Delta \Phi_0^{(j)}-\mathcal E_0^{(j)}.
\end{equation*}
More precisely, to improve $\mathcal E_0^{(1)}$, we take the corrections in the form
\begin{equation*}
\Phi_0^{(1)}=\begin{bmatrix}
\varphi_0^{(1)}(r_1,t)\cos\theta_1\\
0\\
\varphi_0^{(1)}(r_1,t)\sin\theta_1\\
\end{bmatrix}, \quad x-\xi^{(1)}=r_1 e^{i\theta_1},
\end{equation*}
where
\begin{equation*}
\begin{aligned}
&\varphi_0^{(1)}(r_1,t)=-\int_T^t\dot \la_1(s)r_1 K(z_1(r_1),t-s)ds,\\
&z_1(r_1)=(r_1^2+(\la^{(1)})^2)^{1/2},\quad K(z,t)=\frac{2(1-e^{-\frac{z^2}{4t}})}{z^2}.
\end{aligned}
\end{equation*}
The reason for the above correction is the following. The slow decaying error in $\mathcal E_0^{(1)}$ is
\begin{equation*}
\begin{aligned}
\dot\la^{(1)} \pd_{\la^{(1)}} U^{(1)}=&~\dot\la^{(1)}(\la^{(1)})^{-1}  Z^{(1)}_{0,1}(y_1)\\
\approx&~-\dot\la^{(1)}\frac{2r_1}{r_1^2+(\la^{(1)})^2}\begin{bmatrix}
\cos\theta_1\\
0\\
\sin\theta_1\\
\end{bmatrix}\approx -\frac{2\dot\la^{(1)}}{r_1}\begin{bmatrix}
\cos\theta_1\\
0\\
\sin\theta_1\\
\end{bmatrix}.
\end{aligned}
\end{equation*}
Then the scalar $\varphi_0^{(1)}$ roughly solves
$$
\partial_t \varphi_0^{(1)}=\partial_{r_1r_1}\varphi_0^{(1)}+ \frac{1}{r_1}\partial_{r_1}\varphi_0^{(1)}-\frac{1}{r_1^2}\varphi_0^{(1)}-\frac{2\dot\la^{(1)}}{r_1},
$$
whose explicit form was derived in \cite{17HMF}. Similarly, for the slowing decaying error in $\mathcal E_0^{(2)}$
\begin{equation*}
\begin{aligned}
\dot\la^{(2)} \pd_{\la^{(2)}} U^{(2)} + \dot \omega \pd_{\omega} U^{(2)}=&~\dot\la^{(2)}(\la^{(2)})^{-1}Q_{\omega} Z_{0,1}^{(2)}+\dot\omega Q_{\omega} Z_{0,2}^{(2)}\\
\approx &~ -\frac{2r_2}{r_2^2+(\la^{(2)})^2}Q_{\omega}\left(\dot\la^{(2)}\begin{bmatrix}
\cos\theta_2\\
\sin\theta_2\\
0\\
\end{bmatrix}+\la^{(2)}\dot\omega\begin{bmatrix}
-\sin\theta_2\\
\cos\theta_2\\
0\\
\end{bmatrix}\right)\\
\approx&~ -\frac{2\dot\la^{(2)}}{r_2}\begin{bmatrix}
\dot p(t)e^{i\theta_2}\\
0\\
\end{bmatrix},
\end{aligned}
\end{equation*}
we add a global correction of the form
\begin{equation*}
\Phi_0^{(2)}=\begin{bmatrix}
\varphi_0^{(2)}(r_2,t)e^{i\theta_2}\\
0\\
\end{bmatrix}, \quad x-\xi^{(2)}=r_2 e^{i\theta_2},
\end{equation*}
where
\begin{equation*}
\begin{aligned}
&\varphi_0^{(2)}(r_2,t)=-\int_T^t\dot p_2(s)r_2 K(z_2(r_2),t-s)ds,\\
&p_2(t)=\la^{(2)}(t)e^{i\omega(t)},\quad z_2(r_2)=(r_2^2+(\la^{(2)})^2)^{1/2}.
\end{aligned}
\end{equation*}
Here the complex notation
\begin{equation*}
v=\begin{bmatrix}
a+ib\\
c\\
\end{bmatrix}
\end{equation*}
means the 3-vector $v=(a,b,c)^T$. We then write
\begin{equation}\label{def-Phi0}
\Phi_0:=\Phi_0^{(1)}+\Phi_0^{(2)}.
\end{equation}
By direct computations, the new error produced by $\Phi_0$ is
$$
\partial_t\Phi_0 - \Delta_x \Phi_0  +  \EE_0^{(1)}+\EE_0^{(2)}    =
\sum_{j=1}^2(\ttt \RR_0^{(j)} +\ttt \RR_1^{(j)}) , $$
$$\ttt \RR_0^{(1)} =  Q_*\begin{bmatrix} \RR_0^{(1)}   \\ 0 \end{bmatrix}  ,\quad
\ttt\RR_1^{(1)} =  Q_*\begin{bmatrix} \RR_1^{(1)}   \\ 0 \end{bmatrix}, \quad \ttt \RR_0^{(2)} =  \begin{bmatrix} \RR_0^{(2)}   \\ 0 \end{bmatrix}  ,\quad
\ttt\RR_1^{(2)} =  \begin{bmatrix} \RR_1^{(2)}   \\ 0 \end{bmatrix},
$$
where $Q_*$ is defined in \eqref{def-Q_*},
\begin{equation}
\begin{aligned}
\RR_0^{(1)}&:=  - r_1e^{i\theta_1}   \frac {(\la^{(1)})^2}{z_1^4} \int_{-T}^t  \dot \la^{(1)}(s)  ( z_1{K_{z_1}} - z_1^2 K_{z_1 z_1}) (z_1(r_1),t-s) \, ds,\\
\RR_1^{(1)} & :=
- e^{ i\theta_1}  {\rm Re}\,( e^{-i\theta_1} \dot \xi^{(1)}(t))
 \int_{-T}^t  \dot \la^{(1)}(s) \, K(z_1(r_1),t-s) \, ds
\\
&\qquad
+  \frac{r_1}{z_1^2} e^{i\theta_1} \, (\la^{(1)}\dot\la^{(1)}(t)  -  {\rm Re}\,( re^{i\theta_1} \dot\xi^{(1)}(t)) )
\int_{-T}^t  \dot \la^{(1)}(s) \ {z_1K_{z_1}}(z_1(r_1),t-s)\,  ds,\\
\RR_0^{(2)}&:=  - r_2e^{i\theta_2}   \frac {(\la^{(2)})^2}{z_2^4} \int_{-T}^t  \dot p_2(s)  ( z_2{K_{z_2}} - z_2^2 K_{z_2 z_2}) (z_2(r_2),t-s) \, ds,\\
\RR_1^{(2)} & :=
- e^{ i\theta_2}  {\rm Re}\,( e^{-i\theta_2} \dot \xi^{(2)}(t))
 \int_{-T}^t  \dot p_2(s) \, K(z_2(r_2),t-s) \, ds
\\
&\qquad
+  \frac{r_2}{z_2^2} e^{i\theta_2} \, (\la^{(2)}\dot\la^{(2)}(t)  -  {\rm Re}\,( r_2e^{i\theta_2} \dot\xi^{(2)}(t)) )
\int_{-T}^t  \dot p_2(s) \ {z_2K_{z_2}}(z_2(r_2),t-s)\,  ds.\\
\end{aligned}
\end{equation}
Observe that $\RR_1^{(j)}$ is of smaller order. Moreover, we can evaluate
\begin{align*}
\sum_{j=1}^2\ttt L^{(j)}_{U}[\Phi_0^{(j)}]  -\partial_t U_* + \Delta \Phi_0 -\partial_t\Phi_0
=  \KK_{0}^{(j)} + \KK_1^{(j)}-\Pi_{U^\perp}^{(j)} [\ttt \RR_1^{(j)}]
\end{align*}
where for $j=1,2,$ operators $\ttt L^{(j)}_{U}$ and $\Pi_{U^\perp}^{(j)}$ are defined in \eqref{def-tildeL} and \eqref{def-proj}, respectively, and
\begin{align}\label{def-K0j}
\KK_{0}^{(j)} :=  \KK_{01}^{(j)}+ \KK_{02}^{(j)}
\end{align}
with
\begin{equation}\label{def-K1j}
\begin{aligned}
\KK_{01}^{(1)}
:=& - \frac {2}{\la^{(1)}} \rho_1 w_{\rho_1}^2
\int_{-T} ^t   \dot \la^{(1)}(s)   E^{(1)}_1
 K(z_1,t-s)  \, ds, \\
\KK_{02}^{(1)}
:= & \frac 1{\la^{(1)}} \rho_1 w_{\rho_1}^2  \left [  {\dot\la^{(1)}}
-
\int_{-T} ^t   \dot \la^{(1)}(s)  r_1 K_{z_1}(z_1,t-s) {(z_1)}_{r_1} \, ds\, \right]   E^{(1)}_1
\\
&~\quad
-   \frac{1}{4\lambda^{(1)}} \rho_1 w_{\rho_1}^2 \cos w(\rho_1)  \left [  \int_{-T}^t    \dot \la^{(1)}(s)
\, ( z_1{K_{z_1}} - z_1^2 K_{z_1 z_1}) (z_1,t-s)\, ds\, \right ] E_1^{(1)},
\\
\KK_{1}^{(1)}
 :=&
\frac 1{\la^{(1)}}  w_{\rho_1} \, \big [
\dot\xi_1^{(1)}\cos\theta_1  E^{(1)}_1
+ \dot\xi_1^{(1)}\sin\theta_1 E_2^{(1)}       \big ],\\
\KK_{01}^{(2)}
:=& - \frac {2}{\la^{(2)}} \rho_2 w_{\rho_2}^2
\int_{-T} ^t  \left [ {\rm Re  } \,( \dot p_2(s) e^{-i\omega(t)} )   Q_{\omega} E^{(2)}_1+
 {\rm Im  } \,( \dot p_2(s) e^{-i\omega(t)} ) Q_{\omega} E_2^{(2)}   \right ]
\cdot  K(z_2,t-s)  \, ds ,\\
\KK_{02}^{(2)}
:= & \frac 1{\la^{(2)}} \rho_2 w_{\rho_2}^2  \left [  {\dot\la^{(2)}}
-
\int_{-T} ^t  {\rm Re  } \,( \dot p_2(s) e^{-i\omega(t)} ) r_2 K_{z_2}(z_2,t-s) {(z_2)}_{r_2} \, ds\, \right]  Q_{\omega} E^{(2)}_1
\\
&~\quad
-   \frac{1}{4\lambda^{(2)}} \rho_2 w_{\rho_2}^2 \cos w(\rho_2)  \left [  \int_{-T}^t   {\rm Re}\, ( \dot p_2(s)e^{-i\omega(t) }  )
\, ( z_2{K_{z_2}} - z_2^2 K_{z_2 z_2}) (z_2,t-s)\, ds\, \right ] Q_{\omega} E_1^{(2)}
\\
&~\quad
-    \frac{1}{4\lambda^{(2)}} \rho_2 w_{\rho_2}^2  \left [  \int_{-T}^t   {\rm Im }\, ( \dot p_2(s)e^{-i\omega(t) }  )
\, ( z_2{K_{z_2}} - z_2^2 K_{z_2 z_2}) (z_2,t-s)\, ds\,  \right ]  Q_{\omega} E_2^{(2)}  ,\\
\KK_{1}^{(2)}
 :=&
\frac 1{\la^{(2)}}  w_{\rho_2} \, \big [
\Re \big (  (\dot  \xi^{(2)}_1 - i \dot \xi^{(2)}_2)  e^{i\theta_2 } \big ) Q_{\omega} E^{(2)}_1
+ \Im \big(  (\dot  \xi^{(2)}_1 - i \dot \xi^{(2)}_2)  e^{i\theta_2 } \big ) Q_{\omega} E_2^{(2)}       \big ].
\end{aligned}
\end{equation}

\section{Gluing system and derivation of the dynamics for parameters}

\medskip

In this section, we first formulate the inner-outer gluing system so that blow-up solutions with desired asymptotics can be constructed. Then we derive at leading order the dynamics that the parameter functions should satisfy.

\medskip

Since the target manifold is $\mathbb S^2$, we expect that the real solution $u$ to the harmonic map heat flow takes the form of leading profile $U_*$ plus smaller order terms such that $|u(x,t)|=1$ for all $x$ and $t$. To better evaluate the smaller order terms, we look for a solution $u$ of the form
\begin{equation}
\begin{aligned}
u=&~(1+a)U_*+\Phi-(\Phi\cdot U_*)U_*,\\
\Phi=&~\sum_{j=1}^2 \eta_R^{(j)}\Phi_{\rm in}^{(j)}(y_j,t)+\Phi_{\rm out}(x,t)+\Phi_0\\
\end{aligned}
\end{equation}
with $\Phi_0$ defined in \eqref{def-Phi0},
\begin{equation}
\begin{aligned}
&\Phi_{\rm in}^{(1)}=\varphi_{\rm in,1}^{(1)}(y_1,t) E_1^{(1)} +\varphi_{\rm in,2}^{(1)}(y_1,t) E_2^{(1)},\\
&\Phi_{\rm in}^{(2)}=\varphi_{\rm in,1}^{(2)}(y_2,t) Q_{\omega}E_1^{(2)} +\varphi_{\rm in,2}^{(2)}(y_2,t) Q_{\omega}E_2^{(2)},\\
&\eta_R^{(j)}(x,t)=\eta\left(\frac{x-\xi^{(j)}(t)}{\la^{(j)}(t)R(t)}\right),~\eta(s)=\begin{cases}1,~&\mbox{ for }s<1,\\0,~&\mbox{ for }s>2,\end{cases}
\end{aligned}
\end{equation}
where $a$ is a scalar, $\varphi^{(j)}_{\rm in,1}$, $\varphi^{(j)}_{\rm in,2}$, $\Phi_{\rm out}$ are perturbations of smaller order, and $R(t)$ will be chosen later. $\varphi^{(j)}_{\rm in,1}$, $\varphi^{(j)}_{\rm in,2}$ solve the inner problem near each bubble $U^{(j)}$, while $\Phi_{\rm out}$ handles the region away from the concentration zones.
From $|u|=1$, we see that the scalar $|a|=O(|\Phi|^2)$ is of smaller order. Notice that we only need to solve
$$
S(u)=b(x,t) U_*
$$
for some scalar $b$.
Indeed, since $|u|=1$ is enforced for all $t\in (0,T)$ and $u=U_*+\tilde w$
where the perturbation $\tilde w$ is uniformly small, we have
\begin{equation*}
b(x,t)(U_*\cdot u)= S(u)\cdot u=-\frac12 \frac{d}{dt} |u|^2+\frac12 \Delta |u|^2=0.
\end{equation*}
Thus $b\equiv 0$ follows from $U_*\cdot u\geq \delta>0$, which means $u$ is a solution to the harmonic map heat flow.

Now, we compute the error
\begin{equation*}
\begin{aligned}
S(u)=&~ -a_t U_*-(1+a)\partial_t U_*-\partial_t\Phi+(\Phi\cdot U_*)\partial_t U_*+\partial_t(\Phi\cdot U_*)U_*+\Delta a U_*\\
&~+(1+a)\Delta U_*+2\nabla a\cdot\nabla U_*+\Delta \Phi-\Delta[(\Phi\cdot U_*)U_*] \\
&~+\bigg|\nabla\big((1+a)U_*+\Phi-(\Phi\cdot U_*)U_*\big)\bigg|^2\big[(1+a)U_*+\Phi-(\Phi\cdot U_*)U_*\big]\\
=&~-\partial_t\Phi+\Delta \Phi+S(U_*)+(\Phi\cdot U_*)\partial_t U_*-(\Phi\cdot U_*)\Delta U_*-2\nabla(\Phi\cdot U_*)\cdot \nabla U_*\\
&~+\bigg|\nabla\big((1+a)U_*+\Phi-(\Phi\cdot U_*)U_*\big)\bigg|^2\big[(1+a)U_*+\Phi-(\Phi\cdot U_*)U_*\big]-|\nabla U_*|^2U_*\\
&~+2\nabla a\cdot \nabla U_*+a(\Delta U_*-\partial_t U_*)+\big[\Delta a-a_t+\partial_t (\Phi\cdot U_*)-\Delta(\Phi\cdot U_*)\big]U_*,\\
\end{aligned}
\end{equation*}
and here 
\begin{equation*}
\begin{aligned}
\Phi\cdot U_*=&~\eta_R^{(1)}\Phi_{\rm in}^{(1)}\cdot [U^{(2)}-U^{(2)*}]+\eta_R^{(2)}\Phi_{\rm in}^{(2)}\cdot [U^{(1)}-U^{(2)*}]+(\Phi_{\rm out}+\Phi_0)\cdot U_*.\\
\end{aligned}
\end{equation*}
To formulate the inner--outer gluing system, we start from $S(u)=b(x,t)U_*$ and neglect terms in $U_*$ direction due to the discussions above. One expects that the inner solution $\Phi_{\rm in}^{(j)}$ solves the linearization around the bubble $U^{(j)}$, while $\Phi_{\rm out}$ solves a non-homogenous heat equation dealing with all $\R^2_+$ including the regions away from two concentration zones. This leads to the following sufficient condition for $S(u)=b(x,t)U_*$ to hold: $\{\Phi_{\rm in}^{(j)},~\Phi_{\rm out}\}$ solve the inner--outer gluing system 
\begin{equation}\label{sys-inout}
\begin{aligned}
\partial_t\Phi_{\rm in}^{(j)}=&~\Delta \Phi_{\rm in}^{(j)}+|\nabla U^{(j)}|^2 \Phi_{\rm in}^{(j)}+2(\nabla U^{(j)}\cdot \nabla \Phi_{\rm in}^{(j)}) U^{(j)}\\
&~+\tilde L^{(j)}_{U}\Phi_{\rm out}+\mathcal K_0^{(j)}+\mathcal K_1^{(j)} ~\mbox{in}~B^{(j)}_{2\la^{(j)}R}\times (0,T),\\
\partial_t\Phi_{\rm out}=&~\Delta \Phi_{\rm out}+(1-\eta_R^{(1)}-\eta_R^{(2)})\tilde L^{(j)}_{U}\Phi_{\rm out}\\
&~+(1-\eta_R^{(1)}-\eta_R^{(2)})(\mathcal K_0^{(j)}+\mathcal K_1^{(j)})+\mathcal C_{\rm in}+\mathcal N ~\mbox{in}~\R^2_+\times(0,T),
\end{aligned}
\end{equation}
where $B_{x,2\la^{(j)}R}^{(j)}:=\{x\in\R^2_+:|x-\xi^{(j)}|\leq 2\la^{(j)}R\}$, the operator $\tilde L^{(j)}_{U}$ is defined in \eqref{def-tildeL}, and
\begin{equation}\label{def-cpls}
\begin{aligned}
\mathcal C_{in}:=&~\sum_{j=1}^2\bigg(\Phi_{\rm in}^{(j)}\Delta \eta_R^{(j)}+2\nabla \eta_R^{(j)}\cdot\nabla \Phi_{\rm in}^{(j)}-\Phi_{\rm in}^{(j)}\partial_t \eta_R^{(j)}\\
&~+\eta_R^{(j)}(\la^{(j)})^{-1}\nabla_{y_j}\Phi_{in}^{(j)}\cdot(\dot\la^{(j)}y_j+\dot\xi^{(j)})\bigg)+\eta_R^{(2)}\dot\omega J \Phi_{\rm in}^{(2)},\\
\end{aligned}
\end{equation}
\begin{equation}\label{def-nlts}
\begin{aligned}
\mathcal N:=&~(\Phi\cdot U_*)\partial_t U_*-(\Phi\cdot U_*)\Delta U_*-2\nabla(\Phi\cdot U_*)\cdot \nabla U_*\\
&~+\big|\nabla[(1+a)U_*+\Phi-(\Phi\cdot U_*)U_*]\big|^2[(1+a)U_*+\Phi-(\Phi\cdot U_*)U_*]\\
&~-|\nabla U_*|^2U_*+\Pi_{U^\perp}^{(j)} [\ttt \RR_1^{(j)}]+2\nabla a\cdot \nabla U_*+a(\Delta U_*-\partial_t U_*)+\tilde{\mathcal E}\\
&~-\sum_{j=1}^2\bigg( |\nabla U^{(j)}|^2 \Phi_{\rm in}^{(j)}+2(\nabla U^{(j)}\cdot \nabla \Phi_{\rm in}^{(j)}) U^{(j)}\\
&~+|\nabla U^{(j)}|^2\Pi_{U^{\perp}}^{(j)}\Phi_{\rm out}-2\nabla(\Phi_{\rm out}\cdot U^{(j)})\nabla U^{(j)}\bigg).\\
\end{aligned}
\end{equation}
Here
$$
J:=\begin{bmatrix}
0&-1&0\\
1&0&0\\
0&0&0\\
\end{bmatrix}.
$$
For the inner problem, we are going to write, in the complex notation
$$\Phi_{{\rm in},\mathbb C}^{(j)}=\varphi_{\rm in,1}^{(j)}  +i\varphi_{\rm in,2}^{(j)},$$
and further decompose in Fourier modes
$$
\Phi_{{\rm in},\mathbb C}^{(j)}(y,t)=\sum_{k\in\mathbb Z} e^{ik\theta_j} \varphi^{(j)}_k(\rho_j,t)
$$
in the corresponding polar coordinates. The inner problem will then be solved mode by mode.
For the outer problem, we write
$$\Phi_{\rm out}=\psi+Z^*$$
with $Z^*=(Z^*_1,Z^*_2,Z^*_3)^T:\R^2_+\times(0,\infty)\to\mathbb R^3$ satisfying
\begin{equation}\label{def-Z*}
\begin{cases}
\partial_t Z^*=\Delta Z^*~&\mbox{ in }\R^2_+\times(0,\infty),\\
\partial_{x_2}Z_1^*(\cdot,t)=0,\quad\partial_{x_2}Z_2^*(\cdot,t)=0,\quad Z_3^*(\cdot,t)=0~&\mbox{ on }\partial\R^2_+\times(0,\infty),\\
Z^*(\cdot,0)=Z^*_0~&\mbox{ in }\R^2_+.\\
\end{cases}
\end{equation}
Here $Z^*$ will be needed in the reduced problems (especially for the scaling parameters).
Then we will get a solution solving the harmonic map heat flow if $(\Phi_{\rm in}^{(1)},\Phi_{\rm in}^{(2)},\psi)$ solve the {\em inner--outer gluing system}
\begin{equation}\label{eqn-inner1}
\left\{
\begin{aligned}
&(\la^{(1)})^2\partial_t \Phi_{\rm in}^{(1)}= L^{(1)}_W[\Phi_{\rm in}^{(1)}]+(\la^{(1)})^2 \left[\tilde L^{(1)}_U[\Phi_{\rm out}]+\mathcal K_0^{(1)}+\mathcal K_1^{(1)}\right]~\mbox{ in }~\mathcal D^{(1)}_{2R}\\
&\Phi_{\rm in}^{(1)}(\cdot,0)=0~\mbox{ in }~B^{(1)}_{2R(0)}\\
&\Phi_{\rm in}^{(1)}\cdot W_1=0~\mbox{ in }~\mathcal D^{(1)}_{2R}\\
\end{aligned}
\right.
\end{equation}
\begin{equation}\label{eqn-inner2}
\left\{
\begin{aligned}
&(\la^{(2)})^2\partial_t \Phi_{\rm in}^{(2)}= L^{(2)}_W[\Phi_{\rm in}^{(2)}]+(\la^{(2)})^2 \left[\tilde L^{(2)}_U[\Phi_{\rm out}]+\mathcal K_0^{(2)}+\mathcal K_1^{(2)}\right]~\mbox{ in }~\mathcal D^{(2)}_{2R}\\
&\Phi_{\rm in}^{(2)}(\cdot,0)=0~\mbox{ in }~B^{(2)}_{2R(0)}\\
&\Phi_{\rm in}^{(2)}\cdot Q_{\omega}W_2=0~\mbox{ in }~\mathcal D^{(2)}_{2R}\\
\end{aligned}
\right.
\end{equation}
\begin{equation}\label{eqn-outer}
\left\{
\begin{aligned}
&\partial_t \psi=\Delta_x\psi+\mathcal G[\la^{(1)},p_2,\xi^{(1)},\xi^{(2)},\Phi_{\rm out},\Phi^{(1)}_{\rm in},\Phi^{(2)}_{\rm in}]~\mbox{ in }~\R^2_+\times(0,T),\\
&\partial_{x_2}\psi_1=\partial_{x_2}\psi_2=\psi_3=0 ~\mbox{ on }~\partial\R^2_+\times(0,T),
\end{aligned}
\right.
\end{equation}
where
\begin{equation*}
\begin{aligned}
\mathcal G[\la^{(1)},p_2,\xi^{(1)},\xi^{(2)},\Phi_{\rm out},\Phi^{(1)}_{\rm in},\Phi^{(2)}_{\rm in}]:=&~(1-\eta_R^{(1)}-\eta_R^{(2)})\tilde L^{(j)}_{U}\Phi_{\rm out}+\mathcal C_{\rm in}+\mathcal N\\
&~+(1-\eta_R^{(1)}-\eta_R^{(2)})(\mathcal K_0^{(j)}+\mathcal K_1^{(j)}),
\end{aligned}
\end{equation*}
 the linearization $L^{(j)}_W[\phi]$ is defined in \eqref{def-linearization}, and
$$\mathcal D^{(j)}_{2R}:=B_{2R}^{(j)}\times(0,T)=\{y_j\in\R^2_+:|y_j|\leq 2R\}\times(0,T)$$
with the radius
\begin{equation}\label{choice-R}
R=R(t)=\la_*(t)^{-\gamma_*}~\mbox{ with }~\la_*(t)=\frac{T-t}{|\log(T-t)|^2}~\mbox{ and }~\gamma_*\in(0,1/2).
\end{equation}
The reason for choosing such $R(t)$ and $\la_*(t)$ will be made clear later on. The boundary conditions in equation \eqref{eqn-outer} actually follow from the fact
\begin{equation*}
\begin{aligned}
\partial_{x_2} U_{*,1}=\partial_{x_2} U_{*,2}=U_{*,3}=0 ~\mbox{ on }~\partial\R^2_+\times(0,T),\\
\partial_{x_2} \Phi_{0,1}^{(1)}=\partial_{x_2} \Phi_{0,2}^{(1)}=\Phi_{0,3}^{(1)}=0 ~\mbox{ on }~\partial\R^2_+\times(0,T)\\
\end{aligned}
\end{equation*}
thanks to the choices of $U^{(1)}$ and the reflection $U^{(2)*}$. Here
$$U_*=(U_{*,1},U_{*,2},U_{*,3})^T,\quad \Phi_0^{(1)}=(\Phi_{0,1}^{(1)},\Phi_{0,2}^{(1)},\Phi_{0,3}^{(1)})^T.$$

\medskip

Next we derive the dynamics for the parameters $\la^{(1)}(t)$, $p_2(t)$, $\xi^{(1)}(t)$, $\xi^{(2)}(t)$ at leading order as  $t\to T$. We assume for now that the function $\Phi_{\rm out}(x,t)$ is fixed and sufficiently regular, and we regard $T$ as a parameter that will always be taken smaller if necessary. We recall that we need $\xi^{(j)}(T)=q^{(j)}$, $\la^{(1)}(T)=\la^{(2)}(T)=0$, where $q^{(1)}\in\partial \R^2_+$, $q^{(2)}\in \mathring{\R}^2_+$ are given.

\medskip

In order to find solutions to the inner problems  \eqref{eqn-inner1} and \eqref{eqn-inner2} with sufficiently fast decay in space and time, one expects certain orthogonality conditions to hold since even the stationary linearized operator (around degree one harmonic maps) has six-dimensional kernel in $L^{\infty}(\R^2)$. In fact, the linear theory that will be discussed in Section \ref{sec-lineartheories} requires the following orthogonality conditions
\begin{equation}\label{ww1}
\int_{B_{2R}^{(j)}} \mathcal H^{(j)}Z_{p,q}^{(j)}(y_j)dy_j =0\quad \mbox{ for all } t\in (0,T),
\end{equation}
where
\begin{align}
\label{def-mathcalH(1)}
    &\mathcal H^{(1)}:=(\la^{(1)})^2 \left[\tilde L^{(1)}_U[\Phi_{\rm out}]+\mathcal K_0^{(1)}+\mathcal K_1^{(1)}\right],\\
\label{def-mathcalH(2)}
    &\mathcal H^{(2)}:=(\la^{(2)})^2 Q_{\omega}^{-1}\left[\tilde L^{(2)}_U[\Phi_{\rm out}]+\mathcal K_0^{(2)}+\mathcal K_1^{(2)}\right],
\end{align}
and $Z_{p,q}^{(j)}$ is given in \eqref{kernels}. Intuitively, if $\lambda^{(j)}(t)$ has a relatively smooth vanishing as $t\to T$, then it is natural that the term $(\la^{(j)})^2 \partial_t\Phi^{(j)}_{\rm in} $ is of smaller order, and the inner problems are approximately of the form
\begin{equation}\label{linearized-elliptic}
\begin{aligned}
&L^{(1)}_ W  [\Phi_{\rm in}^{(1)}] + \mathcal H^{(1)}=0, \quad \Phi_{\rm in}^{(1)}\cdot W_1  =0  \inn B^{(1)}_{2R},\\
&L^{(2)}_ W  [\Phi_{\rm in}^{(2)}] + Q_{\omega}\mathcal H^{(2)}=0, \quad \Phi_{\rm in}^{(2)}\cdot Q_{\omega}W_2  =0  \inn B^{(2)}_{2R}.\\
\end{aligned}
\end{equation}
If there are solutions $\Phi_{\rm in}^{(j)}(y_j,t)$ to \eqref{linearized-elliptic} with sufficiently fast decay, then necessarily \eqref{ww1} hold for $p=0, 1$, $j,q=1,2$.
These orthogonality conditions in turn require the correct choices of the parameter functions so that the solution  $(\Phi_{\rm in}^{(j)},\psi)$ with appropriate asymptotics exists.

\medskip

We first derive the dynamics for the parameters $\la^{(1)}$ (mode 0), $\xi^{(1)}$ (mode 1) appearing in the boundary bubble. Write
\begin{align*}
\label{defB0j}
\mathcal B_{01}^{(1)} [\la^{(1)},\xi^{(1)}] (t) :=  \frac{\la^{(1)}}{2\pi} \int_{B_{2R}^{(1)}}   [ \KK_{0}^{(1)}[\la^{(1)},\xi^{(1)}]+ \KK_{1}^{(1)}[\la^{(1)},\xi^{(1)}]] \cdot Z_{0,1}^{(1)} (y_1)\, dy_1,~j=1,2.
\end{align*}
Combining \eqref{def-K0j}, \eqref{def-K1j} and \eqref{ww1},
 the following expression for $\mathcal B^{(1)}_{01}$ is readily obtained
 by similar computations as in \cite[Section 5] {17HMF}

\begin{align*}
\mathcal B^{(1)}_{01} [\la^{(1)},\xi^{(1)}](t)
&=
  \int_{-T} ^t  \dot \la^{(1)}(s) \,
\Gamma_1 \left ( \frac {(\la^{(1)}(t))^2}{t-s}   \right )  \,\frac{ ds}{t-s}\,  -2 \dot\la^{(1)} (t) +o(1),
\end{align*}
where $o(1)\to0$ as $t\to T$, and $\Gamma_1(\tau)$  is smooth function given by
\begin{align*}
\Gamma_1 (\tau)
&
:=  - \int_0^{\infty} \rho_1^3 w^3_{\rho_1} \left [\tilde  K ( \zeta )
+ 2 \zeta \tilde  K_\zeta (\zeta ) \frac {\rho_1^2} { 1+ \rho_1^2}
-4\cos w(\rho_1) \zeta^2 \tilde  K_{\zeta\zeta} (\zeta)
\right ]_{\zeta = \tau(1+\rho_1^2)}   \, d\rho_1,
\end{align*}
where
\[
\tilde  K(\zeta)  :=  2\frac {1- e^{-\frac{\zeta}4}} {\zeta}.
\]
Here the orthogonality with $Z_{0,2}^{(1)}$ does not contribute in the dynamics since
$$
\left[\KK_0^{(1)}+\KK_1^{(1)}\right]\cdot E_2^{(1)}=0.
$$
Using the expression of $\Gamma_1 (\tau)$, we get
\begin{equation}
\left\{
\begin{aligned}
\label{GammaNear0}
&| \Gamma_1 (\tau)- 1|
 \le   C \tau(1+  |\log\tau|)~~\hbox{ for }\tau<1 ,
\\
\nonumber
&|\Gamma_1 (\tau)|
\le  \frac C\tau~~\hbox{ for }\tau> 1.
\end{aligned}
\right.
\end{equation}
Define
\begin{align*}
\nonumber
a_{01}^{(1)}[\la^{(1)},\xi^{(1)},\Phi_{\rm out}]
&:=
- \frac{ \la^{(1)}}{2\pi} \int_{B_{2R}^{(1)}}    \tilde L^{(1)}_U  [\Phi_{\rm out}]  \cdot Z_{0,1}^{(1)} (y_1)\, dy_1.
\\
\end{align*}
Then the orthogonality condition \eqref{ww1} with $q, j=1$ and $p=0$ is reduced to
\begin{equation}\label{eqn-B01}
\mathcal B_{01}^{(1)}=a_{01}^{(1)}.
\end{equation}
We observe that
\begin{align*}
\mathcal B_{01}^{(1)}
=   \int_{-T} ^{t-(\la^{(1)})^2}    \frac{\dot \la^{(1)}(s)}{t-s}ds\, + O\big( \|\dot \la^{(1)}\|_\infty \big)+o(1) ~ \mbox{ as }~ t\to T.
\end{align*}
To get an approximation for $a_{01}^{(1)}$, we recall the operator $\tilde L^{(1)}_U$ defined in \eqref{def-tildeL(1)-1}-\eqref{def-tildeL(1)-2}. Write
$$
\Phi_{\rm out}=(\Phi_{\rm out,1},\Phi_{\rm out,2},\Phi_{\rm out,3})^T.
$$
We then get
$$
a_{01}^{(1)}[\la^{(1)},\xi^{(1)}, \Phi_{\rm out}] =  (\partial_{x_1} \Phi_{\rm out,1} + \partial_{x_2} \Phi_{\rm out,3} ) + o(1) ~ \mbox{ as }~ t\to T.
$$

 \medskip

Then the reduced problem \eqref{eqn-B01} at mode $0$ can be written in the form
\begin{equation}
\int_{-T} ^{t-(\la^{(1)})^2}     \frac{\dot \la^{(1)}(s)}{t-s}ds
 =
[\partial_{x_1} \Phi_{\rm out,1} + \partial_{x_2} \Phi_{\rm out,3}](\xi^{(1)}(t),t )  + o(1) + O(\|\dot \la^{(1)}\|_\infty),  \\
 \end{equation}
and thus neglecting lower order terms, $\la^{(1)}$ satisfies the following integro-differential equation
\begin{align}
\label{kuj}
\int_{-T} ^{t-(\la^{(1)})^2(t)} \frac{ \dot \la^{(1)}(s)}{t-s}ds   =
\partial_{x_1} \Phi_{\rm out,1}(q^{(1)},0) + \partial_{x_2} \Phi_{\rm out,3}(q^{(1)},0) =: a_0^{(1)*}.
\end{align}
At this point, we make the following assumption
\begin{align}
\label{negativeDiv}
\partial_{x_1} \Phi_{\rm out,1}(q^{(1)},0) + \partial_{x_2} \Phi_{\rm out,3}(q^{(1)},0) <0,
\end{align}
which is achieved by choosing $Z_0^*$ in \eqref{def-Z*}.
\smallskip
Then, equation \eqref{kuj} becomes
\begin{equation}\label{cccc4}
 \int_{-T} ^{t-(\la^{(1)})^2(t)} \frac{ \dot \la^{(1)}(s)}{t-s}ds   =
 - |a_0^{(1)*}| .
\end{equation}
 We claim that a good approximate solution to \eqref{cccc4} as $t\to T$  is given by
\[
\dot \la^{(1)}(t) =  -\frac {\kappa^{(1)}} {|\log(T-t)|^2}
\]
for some $\kappa^{(1)}>0$. Indeed, we have
\begin{align*}
\int_{-T}^{t-(\la^{(1)})^2(t)} \frac {\dot\la^{(1)}(s)}{t-s}\, ds \,  =& \
\int_{-T}^{t-  (T-t)  } \frac{ \dot\la^{(1)}  (s)}{t-s} \, ds +     \, \dot \la^{(1)} (t)\left [ \log (T-t)   - 2\log (\la^{(1)}(t)) \right ]\nonumber \\  &~+   \int_{t-(T-t)   } ^{ t- (\la^{(1)})^2(t)}\frac{\dot \la^{(1)}(s)-\dot \la^{(1)}(t)}{t-s} ds    \nonumber \\
  \approx & \
\int_{-T}^{t } \frac{ \dot\la^{(1)} (s)}{T-s}\, ds\,  - \, \dot \la^{(1)} (t) \log (T-t) \, := \Upsilon(t)
\end{align*}
as $t\to T$. We see that
$$
\log(T-t) \frac {d\Upsilon(t)} {dt}  =
   \frac d{dt} (|\log(T-t)|^2 \, \dot\la^{(1)}(t))= 0
$$
from the explicit form of  $\dot\la^{(1)}(t)$. Thus $\Upsilon(t)$ is a  constant. As a consequence, equation \eqref{cccc4}
is approximately satisfied if $\kappa^{(1)}$ is such that
$$
\kappa^{(1)} \int_{-T}^{T} \frac{ \dot\la^{(1)}  (s)}{T-s}\,ds\ =\ -|a_0^{(1)*}| ,
$$
which finally gives us the approximate expression
$$
\dot\la^{(1)} (t)= -  |\partial_{x_1} \Phi_{\rm out,1}(q^{(1)},0) + \partial_{x_2} \Phi_{\rm out,3}(q^{(1)},0)|\, \dot \la_* (t) ,
$$
where
\begin{equation*}
 \dot \la_* (t) = -\frac { |\log T|}{|\log(T-t)|^2}.
\end{equation*}
Naturally, imposing $\la_*(T) =0$, we then have
\begin{equation}\label{def-la*}
 \la_* (t) \sim  \frac { |\log T|}{|\log(T-t)|^2}(T-t)\, (1+ o(1)) ~\mbox{ as } ~t\to T.
\end{equation}

\medskip

Similarly, for the mode $1$ of the boundary bubble, we define
\begin{align*}
\mathcal B_{1j}^{(1)} [\la^{(1)},\xi^{(1)}] (t)  & :=  \frac{\la^{(1)}}{2\pi} \int_{B_{2R}^{(1)}}   [ \KK_{0}^{(1)}[\la^{(1)},\xi^{(1)}]+ \KK_{1}^{(1)}[\la^{(1)},\xi^{(1)}]] \cdot Z_{1,j}^{(1)} (y_1)\, dy_1, ~~ j=1,2,\\
\mathcal B^{(1)}_{1} [\la^{(1)},\xi^{(1)} ] (t) & :=   \mathcal B^{(1)}_{11}[\la^{(1)},\xi^{(1)}](t) + i \mathcal B^{(1)}_{12}[\la^{(1)},\xi^{(1)}](t).
\end{align*}
Therefore, by \eqref{def-K0j}, \eqref{def-K1j}, \eqref{kernels} and the fact that $\int_0^{\infty} \rho_1 w_{\rho_1}^2 d\rho_1\, =2$, we obtain
\[
\mathcal B_{1}^{(1)} [\la^{(1)},\xi^{(1)}](t)\,  =  \, 2\dot\xi_1^{(1)}(t)e^{i\theta_1}+o(1) ~ \mbox{ as }~ t\to T .
\]
Write
\begin{align*}
a^{(1)}_{1j} [\la^{(1)},\xi^{(1)},\Phi_{\rm out}] &:= \frac{\lambda^{(1)}}{2\pi}
\int_{B_{2R}^{(1)}}  \tilde L_U[\Phi_{\rm out}]\cdot Z^{(1)}_{1,j}(y_1)  \,dy_1, ~~ j=1,2,
\\
a_1^{(1)}[\la^{(1)},\xi^{(1)},\Phi_{\rm out}] & := - ( a^{(1)}_{11}[\la^{(1)},\xi^{(1)},\Phi_{\rm out}] + i a^{(1)}_{12} [\la^{(1)},\xi^{(1)},\Phi_{\rm out}] ) .
\end{align*}
Therefore, the orthogonality \eqref{ww1} with $j=p=q=1$ is reduced to
\begin{equation}\label{eqB1}
\mathcal B_1^{(1)}[\la^{(1)},\xi^{(1)}] =    a_1^{(1)}[\la^{(1)},\xi^{(1)},\Phi_{\rm out}].
\end{equation}
Similarly, since  $\int_0^\infty w_{\rho_1}^2 \cos w(\rho_1) \rho_1 \, d\rho_1 = 0 $, we get
\begin{align*}
a_{11}^{(1)}[\la^{(1)},\xi^{(1)},\Phi_{\rm out}] & =   2  \pd_{x_1} \Phi_{\rm out,2}(\xi^{(1)},t)  \int_{0}^\infty \cos w \,w_\rho^2 \rho \, d\rho  + O(R^{-2})
\\
&  =  o(1) ~\mbox{ as } ~t\to T
\end{align*}
by using \eqref{def-tildeL(1)-1}--\eqref{def-tildeL(1)-2}, and thus
$$
\dot\xi_1^{(1)}(t)=O(R^{-2}).
$$
This means we can have a solution
$$
\xi^{(1)}_1(t)=q^{(1)}_1+O((T-t)^{1+2\gamma_*}),
$$
where $\gamma_*$ is given in \eqref{choice-R}. Note that the imaginary part of $a_1^{(1)}$ vanishes because of the partially free boundary $\partial_{x_2}\psi_2^*(\xi^{(1)},t)=0$
for given $q_1^{(1)}\in\R$.

\medskip

For the parameters involved in the interior bubble, one can carry out a similar analysis for $\la^{(2)}(t)$ at mode 0 and for $\xi^{(2)}(t)$ at mode 1. In fact, we have
\begin{equation}
    \begin{aligned}
    &\int_{-T}^{t-(\la^{(2)})^2(t)}\frac{\dot p_2(s)}{t-s}ds=(\div \Phi_{\rm out}+i\curl \Phi_{\rm out})(\xi^{(2)},t)+O(\|\dot p_2\|_{\infty})+o(1),\\
    &\dot \xi^{(2)}(t)=O(R^{-2}) ~\mbox{ as }~t\to T,
    \end{aligned}
\end{equation}
where we have used the complex notation \eqref{def-complexnotations} writing
$$
\Phi_{\rm out}=\begin{bmatrix}
\Phi_{\rm out,1}+i\Phi_{\rm out,2}\\
\Phi_{\rm out,3}\\
\end{bmatrix}.
$$
For the expected asymptotics of the blow-up speed $p_2(t)=\la^{(2)}e^{i\omega}$ to exist, the following sign condition near the interior bubble is needed
$$
\div \Phi_{\rm out}(q^{(2)},0)<0,
$$
which determines the parameters
\begin{equation}
    \begin{aligned}
    &\int_{-T}^{t-(\la^{(2)})^2(t)}\frac{\dot \la^{(2)}(s)}{t-s}ds=-|\div \Phi_{\rm out}(q^{(2)},0)|+o(1),\\
    &\omega=\omega_0:=\arctan \frac{\curl \Phi_{\rm out}}{\div \Phi_{\rm out}}(q^{(2)},0),\\
    &\xi^{(2)}=q^{(2)}+O((T-t)^{1+2\gamma_*}),
    \end{aligned}
\end{equation}
for some given point $q^{(2)}\in\mathring{\R}^2_+$.

\medskip

In conclusion, orthogonality conditions, which are required to guarantee well-behaved solutions to the inner problems \eqref{eqn-inner1} and \eqref{eqn-inner2},  result in the following asymptotics of the modulation parameters:
\begin{equation}
\begin{aligned}
    &\la^{(1)}(t)\sim \kappa^{(1)}\frac{T-t}{|\log(T-t)|^2},\quad &\xi^{(1)}(t)\sim q^{(1)}\in\partial\R^2_+\\
    &p_2(t)\sim \kappa^{(2)}\frac{T-t}{|\log(T-t)|^2}e^{i\omega_0},\quad &\xi^{(2)}(t)\sim q^{(2)}\in\mathring{\R}^2_+\\
\end{aligned}
\end{equation}
for some $\kappa^{(1)},~\kappa^{(2)}>0$.

\bigskip

\section{Linear theories for linearized harmonic map heat flow and Stokes system}\label{sec-lineartheories}

\medskip

In this section, we give linear theories that are needed to solve the linearized harmonic map heat flow, namely the inner and outer problems \eqref{eqn-inner1}--\eqref{eqn-outer}, and also the Stokes system.

\medskip

\noindent {\bf Linear theory for the inner problems.} We first start from the model inner problem
\begin{equation*}\label{eqn-linearinner}
\begin{cases}
(\la^{(j)})^2 \partial_t \phi^{(j)} =L_W^{(j)}[\phi^{(j)}]+h^{(j)}(y_j,t)~&\mbox{ in }~\mathcal D^{(j)}_{2R},\\
\phi^{(j)}(\cdot,0)=0~&\mbox{ in }~B^{(j)}_{2R(0)},\\
\phi^{(j)}\cdot W_j=0~&\mbox{ in }~\mathcal D^{(j)}_{2R}\\
\end{cases}
\end{equation*}
for $j=1,2$, where we write
\begin{equation}\label{def-phi1phi2}
    \begin{aligned}
    &\phi^{(1)}:=\Phi_{\rm in}^{(1)},\quad h^{(1)}=\mathcal H^{(1)},\\
    &\phi^{(2)}:=Q_{\omega}^{-1}\Phi_{\rm in}^{(2)},\quad h^{(2)}=Q_{\omega}^{-1}\mathcal H^{(2)},\\
    \end{aligned}
\end{equation}
and we recall that
$$R=R(t)=\la_*(t)^{-\gamma_*}~\mbox{ with }~\la_*(t)=\frac{T-t}{|\log(T-t)|^2}~\mbox{ and }~\gamma_*\in(0,1/2).$$
For notational simplicity, in the rest of this section, we drop the superscripts and just write
\begin{equation}\label{eqn-linearinner}
\begin{cases}
\la^2 \partial_t \phi =L_W[\phi]+h(y,t)~&\mbox{ in }~\mathcal D_{2R},\\
\phi(\cdot,0)=0~&\mbox{ in }~B_{2R(0)},\\
\phi\cdot W=0~&\mbox{ in }~\mathcal D_{2R}\\
\end{cases}
\end{equation}
since the linear theory applies to both inner problems. For example, if we apply the linear theory of the model problem \eqref{eqn-linearinner} to the inner problem of the boundary bubble, then the spatial variable $y$ in \eqref{eqn-linearinner} stands for the rescaled variable $y_1$ around the boundary concentration point.

\medskip

Since the inner problem for the interior bubble does not touch the boundary $\partial\R^2_+$, one may regard it as a problem in $\R^2$ with compact support. But in the inner problem for the boundary bubble, the partially free boundary conditions play an important role. To be more precise, the partially free boundary conditions \eqref{FB} determine the symmetry of the inner problem across the boundary, and this symmetry allows us to do the Fourier expansion in modes and regard the half space problem as problem in the entire space $\R^2$, and thus the linear theory applies. Before we explain the reflection of $\Phi_{\rm in}^{(1)}$, we first introduce the Fourier modes of problem \eqref{eqn-linearinner}.

\medskip

We regard $h(y,t)$ as a function defined in $\R^2\times(0,T)$ with compact support, and assume that $h(y,t)$ has the space-time decay of the following type
$$|h(y,t)|\lesssim \frac{\la_*^{\nu}(t)}{1+|y|^a},\quad h\cdot W=0,$$
where $\nu>0$ and $a\in(2,3)$. Define the norm
\begin{equation*}
\|h\|_{\nu,a}:=\sup_{(y,t)\in\R^2\times(0,T)} \la_*^{-\nu}(t)(1+|y|^a)|h(y,t)|.
\end{equation*}
In  polar coordinates, $h(y,t)$ can be written as
\begin{equation*}
h(y,t)=h^1(\rho,\theta,t) E_1(y) +h^2(\rho,\theta,t) E_2(y),~y=\rho e^{i\theta}
\end{equation*}
since $h\cdot W=0$. We use the complex notation
$$
\tilde h(\rho,\theta,t):=h^1+i h^2
$$
and expand in Fourier series
\begin{equation}\label{def-Fourier1}
\tilde h(\rho,\theta,t)=\sum_{k=-\infty}^{\infty} \tilde h_k(\rho,t) e^{ik\theta}
\end{equation}
such that
\begin{equation}\label{def-Fourier2}
h(y,t)=\sum_{k=-\infty}^{\infty} h_k(y,t) :=h_0(y,t)+h_1(y,t)+h_{-1}(y,t)+h_{\perp}(y,t)
\end{equation}
with
\begin{equation}\label{def-Fourier3}
h_k(y,t)={\rm Re}(\tilde h_k(\rho,t)e^{ik\theta})E_1+{\rm Im}(\tilde h_k(\rho,t)e^{i k\theta}) E_2, ~~ k\in\mathbb Z.
\end{equation}
We consider the kernel functions $Z_{k,j}$ (dropping superscripts) introduced in \eqref{kernels}, and define
\begin{equation}\label{def-hbar}
\bar h_{k}(y,t)   :=       \sum_{j=1}^2   \frac {\chi Z_{k,j}(y)} { \int_{\R^2} \chi  |Z_{k,j} |^2  }  \, \int_{ \R^2 }h(z  , t )  \cdot Z_{k,j}(z)\, dz, \quad~k=0,\pm 1,~j=1,2,
\end{equation}
where
$$
\chi(y,t) = \begin{cases}
w_\rho^2(|y|)  & \hbox{ if }  |y|< 2R(t),\\
0& \hbox{ if }  |y|\ge 2R(t).
\end{cases}
$$
Similarly, we decompose the inner solution
\begin{equation*}
\phi=\sum_{k=-\infty}^{\infty} \phi_k,~~\quad \phi_k(y,t)={\rm Re} (\varphi_k(\rho,t) e^{ik\theta}) E_1 + {\rm Im} (\varphi_k(\rho,t) e^{ik\theta}) E_2.
\end{equation*}
In each mode $k$, the pair $(\phi_k,h_k)$ satisfies
\begin{equation}\label{eqn-mode}
\begin{cases}
\la^2\partial_t \phi_k = L_W[\phi_k]+h_k(y,t)~&\mbox{ in }~\mathcal D_{4R},\\
\phi_k(y,0)=0,~&\mbox{ in }~B_{4R(0)},\\
\end{cases}
\end{equation}
which is equivalent to the following problem
\begin{equation*}
\begin{cases}
\la^2\partial_t \varphi_k = \mathcal L_k [\varphi_k]+\tilde h_k(\rho,t)~&\mbox{ in }~\tilde{\mathcal D}_{4R},\\
\varphi_k(\rho,0)=0~&\mbox{ in }~(0,4R(0)),\\
\end{cases}
\end{equation*}
where $\tilde{\mathcal D}_{4R}=\left\{(\rho,t):~t\in(0,T),~\rho\in(0,4R(t))\right\}$, and
\begin{equation*}
\mathcal L_k [\varphi_k] := \partial_{\rho\rho} \varphi_k +\frac{\partial_{\rho}\varphi_k}{\rho}-(k^2+2k\cos w +\cos (2w))\frac{\varphi_k}{\rho^2}.
\end{equation*}
It is direct to see that the kernel functions for $\mathcal L_k$ such that $\mathcal L_k[Z_k]=0$ at modes $k=0,\pm 1$ are given by
\begin{equation}\label{kernels''}
Z_0(\rho)=\frac{\rho}{1+\rho^2},~~Z_1(\rho)=\frac{1}{1+\rho^2},~~Z_{-1}(\rho)=\frac{2\rho^2}{1+\rho^2}.
\end{equation}

\medskip

Now, let us go back to the Fourier expansion for the boundary inner problem. In the general setting (without symmetry across $\partial\R^2_+$), we try to look for a solution near the boundary concentration zone in the following form
\begin{equation}\label{FE}
\begin{aligned}
\phi^{(1)}(y_1,t):=&~\sum_{k\in\mathbb Z}\phi^{(1)}_k(y_1,t)\\
=&~\sum_{k\in\mathbb Z}\left[{\rm Re}(e^{ik\theta_1}\varphi^{(1)}_k(\rho_1,t))E_1^{(1)}+{\rm Im}(e^{ik\theta_1}\varphi^{(1)}_k(\rho_1,t))E_2^{(1)}\right]\\
=&~\sum_{k\in\mathbb Z}\left[(\varphi_{k,1}\cos(k\theta_1)-\varphi_{k,2}\sin(k\theta_1))E^{(1)}_1+(\varphi_{k,2}\cos(k\theta_1)+\varphi_{k,1}\sin(k\theta_1))E_2^{(1)}\right],\\
\end{aligned}
\end{equation}
where $$\varphi_k^{(1)}=\varphi_{k,1}+i\varphi_{k,2}.$$
Recall from \eqref{FB1} that the partially free boundary conditions are automatically satisfied by extending the first, second, third components of $\phi^{(1)}$ evenly, evenly, oddly in $x_2$, and we notice that the first, second and third components of $E_1^{(1)}$ are even, even, odd in $x_2$, while the first, second, third components of $E_2^{(1)}$ are odd, even/odd, even in $x_2$, respectively (see \eqref{def-E1E2}). So the terms $$\varphi_{k,1}\cos(k\theta_1)E_1^{(1)}+\varphi_{k,1}\sin(k\theta_1) E_2^{(1)}$$
have the right symmetry, but
$$
\varphi_{k,2}\cos(k\theta_1)E_2^{(1)}-\varphi_{k,2}\sin(k\theta_1)E_1^{(1)}
$$
violate the partially free boundary conditions. In other words, if
$$
\varphi_{k,2}\equiv 0,
$$
then the Fourier expansion \eqref{FE} already implies that $\phi^{(1)}$ satisfies the partially free boundary conditions. In fact, the role of the reflected bubble is to ensure that the error produced by $U_*$ is perpendicular to $E_2^{(1)}$ on $\partial\R^2_+$ so that $\varphi_{k,2}\equiv 0$, and this in turn rules out the possibility of rotations for the boundary bubble since $E_2^{(1)}$-direction corresponds exactly to the rotation around $z$-axis. This ``rigidity'' is consistent with the intuition from \eqref{FB}$_3$ that the image of the boundary under map $u$ is fixed on the equator $\Sigma$.

\medskip

Recalling the right hand side of the boundary inner problem $\mathcal H^{(1)}$ in \eqref{def-mathcalH(1)}, one needs to check
$$
\mathcal H^{(1)}\cdot E_2^{(2)}=0~\mbox{ on }~\partial\R^2_+.
$$
The most important term is in fact the coupling from the outer problem
$$
\tilde L^{(1)}_{U}[\Phi_{\rm out}]=-\frac2{\la^{(1)}} w_{\rho_1} [(\partial_{r_1}\Phi_{\rm out}\cdot U^{(1)})E_1^{(1)}-\frac1{r_1}(\partial_{\theta_1}\Phi_{\rm out}\cdot U^{(1)})E_2^{(1)}],
$$
and to avoid projection onto $E_2^{(1)}$ direction, we only need
$$\partial_{\theta_1}\Phi_{\rm out}\cdot U^{(1)}\Big|_{x_2=0}=0,$$
i.e.,
$$
(x_1\partial_2-x_2\partial_1)\Phi_{\rm out}\cdot U^{(1)}\Big|_{x_2=0}=0.
$$
This is true thanks to
$$
\partial_2 \Phi_{\rm out,1}=\partial_2 \Phi_{\rm out,2}=0 \mbox{ on } \partial\R^2_+,
$$
where we have used the explicit expression of $U^{(1)}$ in \eqref{def-U1U2}, and $\Phi_{\rm out,1}$, $\Phi_{\rm out,2}$ are the first and second components of $\Phi_{\rm out}$, respectively. Above analysis implies that for $\mathcal H^{(1)}$, there is no projection onto $E_2^{(2)}$ on the boundary, so we can regard the boundary inner problem as a problem in the entire space $\R^2$ (an interior problem).

\medskip

We now state the linear theory for both inner problems.

\begin{prop}[\cite{LCF2D}]
\label{prop-lt}
Assume that $a\in(2,3)$, $\nu>0$, $\delta\in(0,1)$ and $\|h\|_{\nu,a} <+\infty$.
Let us write $$h = h_0 + h_1 + h_{-1} + h_\perp ~\mbox{ with }~ h_{\perp} = \sum_{k\not=0,\pm 1} h_k.$$
Then there exists a solution $\phi[h]$ of problem $\eqref{eqn-linearinner}$, which defines
a linear operator of $h$, and satisfies the following estimate in $\mathcal{D}_{2R}$
\begin{align*}
&\quad
|\phi(y,t)|+(1+|y|)\left |\nabla_y \phi(y,t)\right |  +  (1+|y|)^2\left |\nabla^2_y \phi(y,t)\right |
\\
& \lesssim
   \lambda^\nu_*(t)   \,
 \min\left\{\frac{R^{\delta(5-a)}(t)}{1+|y|^3},  \frac{1}{1+|y|^{a-2}} \right\}
 \, \| h_0 -\bar h_{0} \|_{\nu,a}
+
\frac{ \lambda^\nu_*(t)  R^2(t)} {1+ |y|}  \|\bar h_0\|_{\nu,a}
\\
& \quad
+   \frac{ \lambda^\nu_*(t) }{ 1+ |y|^{a-2} }\, \left \| h_1 - \bar h_1\right  \|_{\nu,a}
+   \frac{ \lambda^\nu_*(t)  R^4(t)} {1+ |y|^2} \left \| \bar h_{1} \right  \|_{\nu,a}
\\
& \quad
+
  \lambda^\nu_*(t)
 \, \| h_{-1} -\bar h_{-1} \|_{\nu,a}
 + \lambda^\nu_*(t) \log R(t) \,   \| \bar h_{-1} \|_{\nu,a}
\\
& \quad
+
\frac{ \lambda^\nu_*(t)  }{ 1+ |y|^{a-2} }\,   \|h_\perp \|_{\nu,a} .
\end{align*}
\end{prop}

\medskip

\noindent {\bf Linear theory for the outer problem.} We then introduce the linear theory for the outer problem. For the outer problem, we will solve it componentwise since each component satisfies a nonhomogeneous heat equation. Because of the symmetry imposed on the solution, we use the Duhamel's formula for the model linear problem
\begin{equation}\label{heat-eq0}
\psi_t=\Delta_{\R^2} \psi+g
\end{equation}
with either Dirichlet or Neumann boundary condition. More precisely, for $\psi=(\psi_1,\psi_2,\psi_3)^T$ and $g=(g_1,g_2,g_3)^T$, we have
\begin{equation}
\begin{aligned}
&\psi_i(x,t)=\int_0^t\int_{\R^2_+}[\Gamma(x-y,t-s)+\Gamma(x-y^*,t-s)]g_i(y,s)dyds,\quad i=1,2,\\
&\psi_3(x,t)=\int_0^t\int_{\R^2_+}[\Gamma(x-y,t-s)-\Gamma(x-y^*,t-s)]g_3(y,s)dyds,\\
\end{aligned}
\end{equation}
where $y^*$ is the reflection of $y$. Clearly, for $j=1,2,3,$ one has
$$
|\psi_j(x,t)|\lesssim \int_0^t\int_{\R^2}\Gamma(x-y,t-s)|g_i(y,s)|dyds:= \left|\Gamma \underset{(x,t)}{*} |g|\right|
$$
since $|x-y^*|\geq |x-y|$ for $y\in\R^2_+$. So the upper bound of $\left|\Gamma \underset{(x,t)}{*} |g|\right|$ implies a weighted-$L^{\infty}$ estimate of $\psi$. We first define the weights
\begin{align}\label{weights}
\left\{
\begin{aligned}
\varrho_1^{(j)} & :=   \lambda_*^{\Theta}  (\lambda_* R)^{-1}  \chi_{ \{x\in\R^2_+: |x-q^{(j)}|\leq 3\la_*R \} },
\\
\varrho_2 & := T^{-\sigma_0} \lambda_*^{1-\sigma_0}\left( \frac{1}{|x-q^{(1)}|^2}+\frac{1}{|x-q^{(2)}|^2} \right) \chi_{\left\{ \bigcap_{j=1}^2\{x\in\R^2_+:|x-q^{(j)}| \geq  \la_*R \}\right\} },
\\
\varrho_3 & := T^{-\sigma_0},
\end{aligned}
\right.
\end{align}
where $q^{(1)}\in\partial\R^2_+$, $q^{(2)}\in\mathring{\R}^2_+$,  $\Theta>0$ and  $\sigma_0>0$ is  small.
For a function $g(x,t)$ we define the $L^\infty$-weighted norm
\begin{align}
\label{defNormRHSpsi}
\|g\|_{**} : =   \sup_{ \R^2_+ \times (0,T)}  \left(\frac{|g(x,t)|}{1+\varrho_1^{(1)}+\varrho_1^{(2)}+\varrho_2+\varrho_3}\right).
\end{align}
We define the $L^\infty$-weighted norm for $\psi$
\begin{equation}\label{normPsi}
\begin{aligned}
\| \psi\|_{\sharp, \Theta,\gamma}
&:=
\lambda^{-\Theta}_*(0)
\frac{1}{|\log T|  \lambda_*(0) R(0) }\|\psi\|_{L^\infty(\R^2_+\times (0,T))}
+ \lambda^{-\Theta}_*(0) \|\nabla_x \psi\|_{L^\infty(\R^2_+\times (0,T))}
\\
&~\quad
+
\sup_{\R^2_+\times (0,T)}   \lambda^{-\Theta-1}_*(t) R^{-1}(t)
\frac{1}{|\log(T-t)|} |\psi(x,t)-\psi(x,T)|
\\
&~\quad
+ \sup_{\R^2_+\times (0,T)} \, \lambda^{-\Theta}_*(t)
|\nabla_x \psi(x,t)-\nabla_x \psi(x,T) |+ \|\nabla_x^2 \psi\|_{L^{\infty}(\R^2_+\times(0,T))}
\\
& ~\quad
+ \sup_{}
\lambda^{-\Theta}_*(t)
(\lambda_*(t) R(t))^{2\gamma}  \frac {|\nabla_x \psi(x,t) -\nabla_x \psi(x',t') |}{ ( |x-x'|^2 + |t-t'|)^{\gamma   }} ,
\end{aligned}
\end{equation}
where $\Theta>0$, $\gamma \in (0,1/2)$, and the last supremum is taken in the region
\[
x,\, x'\in \R^2_+,\quad  t,\, t'\in (0,T), \quad |x-x'|\le 2 \la_*(t)R(t), \quad  |t-t'| < \frac 14 (T-t) .
\]

The solution $\psi$ to the model outer problem \eqref{heat-eq0} will be measured
in the norm $\|\cdot\ \|_{\sharp,\Theta,\gamma}$ defined in \eqref{normPsi} where $\gamma\in ( 0,1/2)$, and we require that $\Theta$ and $\gamma_*$ (recall that $R= \lambda_*^{-\gamma_*}$ in \eqref{choice-R}) satisfy
\begin{align}
\label{assumpPar1}
\gamma_* \in (0,1/2)  ,
\quad
\Theta \in
(0,\gamma_*).
\end{align}

\medskip

By similar computations as in \cite[Proposition A.1]{17HMF}, we have the following

\begin{prop}
Assume  \eqref{assumpPar1} holds.
For $T>0$ sufficiently small, there is a linear operator that maps a function $g:\R^2_+ \times (0,T) \to \R^3$ with  $\|g\|_{**}<\infty$ into $\psi$ which solves problem \eqref{heat-eq0}.
Moreover, the following estimate holds
\begin{align*}
\left\| \big|\Gamma \underset{(x,t)}{*} |g|\big| \right\|_{\sharp, \Theta ,\gamma}
 \leq C \|g\|_{**} ,
\end{align*}
where $\gamma\in(0,1/2)$.
\end{prop}

\medskip

\noindent {\bf Linear theory for the Stokes system.} In order to deal with the forced Navier-Stokes equation with Navier boundary conditions, we consider the following Stokes system 
\begin{equation}\label{eqn-Stokes}
\begin{cases}
\partial_t v+\nabla P=\Delta v +\nabla\cdot F~&\mbox{ in }~\R^2_+\times (0,\infty),\\
\nabla\cdot v=0~&\mbox{ in }~\R^2_+\times (0,\infty),\\
\partial_{x_2} v_1=v_2=0~&\mbox{ on }~\partial\R^2_+\times (0,\infty),\\
v(\cdot,0)=0~&\mbox{ in }~\R^2_+,\\
\end{cases}
\end{equation}
where $\nabla\cdot F$ is solenoidal (see \eqref{def-solenoidal}). Our aim is to obtain weighted $L^{\infty}$ (in space-time) control of $v$ for given forcing in divergence form. The choice of the weighted bound for the forcing $\nabla\cdot F$ is based on the behavior of inner solution near each concentration $q^{(j)}$. We have the following pointwise estimates

\medskip

\begin{prop}
Let $v$ be the solution to the Stokes system \eqref{eqn-Stokes}. If $F$ satisfies
\EQ{\label{ineqineqineq}|F(x,t)|\lesssim\frac{\la_*^{\nu-2}(t)}{1+\left|\frac{x-q}{\la_*(t)}\right|^{a+1}},\quad |\nabla_x F(x,t)|\lesssim\frac{\la_*^{\nu-3}(t)}{1+\left|\frac{x-q}{\la_*(t)}\right|^{a+2}}}
for any $q\in\overline{\R^2_+}$ fixed, $\nu>0$ and $a>1$, then solution $v$ satisfies the following pointwise bound
\begin{align*}
|v(x,t)|\lesssim \frac{\la_*^{\nu-1}}{1+\left|\frac{x-q}{\la_*(t)}\right|},\quad
|\nabla v(x,t)|\lesssim \frac{\la_*^{\nu-2}}{1+\left|\frac{x-q}{\la_*(t)}\right|}.
\end{align*}
\end{prop}

The proof of above proposition can be done either by the natural reflection or by the Green's tensor derived in Appendix \ref{sec-SO}. Indeed, similar to the discussions in Section \ref{symmetry}, there are natural reflections for $v,\,P$, and $F_D:=\nabla\cdot F:=(F_{D,1},F_{D,2})^T$ in \eqref{eqn-Stokes}
\begin{equation*}
\tilde v(x_1,x_2,t)=\begin{bmatrix}
v_1(x_1,-x_2,t)\\
-v_2(x_1,-x_2,t)\\
\end{bmatrix}, \quad \tilde P(x_1,x_2,t)=P(x_1,-x_2,t),\quad \tilde F_D=\begin{bmatrix}
F_{D,1}(x_1,-x_2,t)\\
-F_{D,2}(x_1,-x_2,t)\\
\end{bmatrix},\quad
x_2<0
\end{equation*}
such that the half space problem can be regarded as an interior problem, and thus all the precise pointwise estimates can be achieved by the Oseen tensor in $\R^2$ (see \cite[Section 3]{LCF2D}). The second method is by explicit Green's tensor for Stokes operator with Navier boundary conditions in the half space.

\bigskip

\section{Solving the partially free boundary system}

\medskip

In this section, we will solve the full system \eqref{LCF}--\eqref{FB}. First, we analyze the effect of the couplings, i.e., the transported term $v\cdot \nabla u$ in the harmonic map heat flow and the forcing $-\varepsilon_0\nabla\cdot \left(\nabla u \odot \nabla u-\frac12 |\nabla u|^2 \mathbb{I}_2\right)$ in the incompressible Navier-Stokes equation. Next, we introduce the weighted topologies for $(v,u)$ and then solve the full system using fixed point argument.

\subsection{Couplings in the full system}

\bigskip

In previous sections, we neglect the transported term and carry out all the elements for the harmonic map heat flow. Now we consider the full transported harmonic map heat flow
$$
u_t+v\cdot \nabla u=\Delta u+|\nabla u|^2u
$$
and analyze the effect of the transported term. Most importantly, we need to check $v\cdot \nabla u$ does not violate the symmetry of the boundary inner problem. In fact, by \eqref{def-E1E2}, one has
\begin{equation*}
    \begin{aligned}
(v\cdot \nabla u)\cdot E_2^{(2)}=&~(v_1\partial_{x_1}u+v_2\partial_{x_2}u)\cdot E_2^{(1)}\\
=&~-(v_1\partial_{x_1}u_1+v_2\partial_{x_2}u_1)\sin\theta_1+(v_1\partial_{x_1}u_3+v_2\partial_{x_2}u_3)\cos\theta_1\\
=&~0 ~\mbox{ on }~\partial\R^2_+
\end{aligned}
\end{equation*}
where we have used the partially free boundary conditions \eqref{FB1}.
On the other hand, it is direct to check that
\begin{equation*}
    \begin{aligned}
v\cdot\nabla u=\begin{bmatrix}
v_1\partial_1 u_1+v_2\partial_2 u_1\\
v_1\partial_1 u_2+v_2\partial_2 u_2\\
v_1\partial_1 u_3+v_2\partial_2 u_3\\
\end{bmatrix}
\end{aligned}
\end{equation*}
satisfies the partially free boundary conditions \eqref{FB1}$_1$--\eqref{FB1}$_3$, meaning that the trasported term is also compatible with the outer problem.

\medskip

For the external forcing coming from the orientation field, we notice that the orientation field $u$ only gets coupled with velocity field $v$ solved from the forcing after Helmholtz projection since the curl-free term is part of the pressure. For the Stokes system, we do Helmholtz projection first
$$
\partial_t v +\nabla P_1 = \Delta v -\mathbf P\left[ \nabla\cdot \left(\nabla  u \odot \nabla u\right)\right],
$$
where $\mathbf P$ is the Helmholtz projection such that $\mathbf P\left[ \nabla\cdot \left(\nabla  u \odot \nabla u\right)\right]$ is solenoidal
$$
\nabla\cdot \left(\nabla u \odot \nabla u-\frac12 |\nabla u|^2 \mathbb{I}_2\right)= \mathbf P\left[ \nabla\cdot \left(\nabla  u \odot \nabla u\right)\right]+\nabla \Phi_P,\quad P=P_1+\Phi_P.
$$
By our ansatz
\begin{equation*}
    \begin{aligned}
u=&~U_*+\Phi-(\Phi\cdot U_*)U_*+aU_*\\
\sim&~ U^{(1)}+U^{(2)}+\sum_{j=1}^2\eta_R^{(j)}\left(\sum_{k\in\mathbb Z} \phi_k^{(j)}(y_j,t)\right)+\Phi_{\rm out}(x,t),
\end{aligned}
\end{equation*}
we know that the forcing in the incompressible Navier-Stokes equation near $q^{(j)}$
$$\nabla\cdot\left(\nabla u\odot \nabla u -\frac12|\nabla u|^2 \mathbb I_2\right)\sim \sum\limits_{k\in\mathbb Z} \nabla\cdot\left(2\nabla U^{(j)}\odot \nabla \phi^{(j)}_k -(\nabla U^{(j)} : \nabla \phi_k^{(j)}) \mathbb I_2\right)$$
since the outer solution $\Phi_{\rm out}$ is smaller, and the leading term actually vanishes
$$\na\cdot(\na U^{(j)}\odot\na U^{(j)}-1/2\,|\na U^{(j)}|^2\,\mathbb{I}_2)=\De U^{(j)}\cdot\na U^{(j)}=-|\na U^{(j)}|^2(U^{(j)}\cdot\na U^{(j)})=0.$$
Here, $\na U^{(j)}:\na\phi_k^{(j)}=\sum\limits_{p,q}\pd_pU^{(j)}_q\pd_p(\phi_k^{(j)})_q$, where $(\phi_k^{(j)})_q$ stands for the $q$-th component of $\phi_k^{(j)}$.
In the sequel, we shall call the following term
$$
\nabla\cdot\left(2\nabla U^{(j)}\odot \nabla \phi_k^{(j)} -(\nabla U^{(j)} : \nabla \phi_k^{(j)}) \mathbb I_2\right)
$$
the forcing at mode $k$. Another useful fact observed in \cite{LCF2D} is that the forcing at mode $0$, which is of largest size on the right hand side of the Stokes system, actually enters into the pressure, so it is not involved in the loop.

\bigskip

\subsection{Inner--outer gluing system and reduced problems}

\medskip

We will get a desired solution $(v,u)$ to the partially free boundary system \eqref{LCF}--\eqref{FB} if $(v,\phi^{(1)},\phi^{(2)},\psi,Z^*,\la^{(1)},\la^{(2)},\omega,\xi^{(1)},\xi^{(2)})$ solves the following {\em inner--outer gluing system}
\begin{equation}\label{eqn-NS}
\left\{
\begin{aligned}
&\partial_t v+\nabla P=\Delta v -\varepsilon_0\nabla\cdot \mathcal F[v,\phi^{(1)},\phi^{(2)},\psi,Z^*,\la^{(1)},\la^{(2)},\omega,\xi^{(1)},\xi^{(2)}]\\
&~\qquad\qquad\qquad\qquad -\nabla\cdot(v\otimes v)~\mbox{ in }~\R^2_+\times(0,T),\\
&\nabla\cdot v=0~\mbox{ in }~\R^2_+\times(0,T),\\
&\partial_{x_2}v_1=v_2=0~\mbox{ on }~\partial\R^2_+\times(0,T),\\
&v(\cdot,0)=v_0~\mbox{ in }~\R^2_+,\\
\end{aligned}
\right.
\end{equation}
\begin{equation}\label{eqn-dinner}
\left\{
\begin{aligned}
&(\la^{(j)})^2\partial_t \phi^{(j)}= L_W^{(j)}[\phi^{(j)}]+h^{(j)}[v,\phi^{(1)},\phi^{(2)},\psi,Z^*,\la^{(1)},\la^{(2)},\omega,\xi^{(1)},\xi^{(2)}]~\mbox{ in }~\mathcal D_{2R}^{(j)},\\
&\phi^{(j)}(\cdot,0)=0~\mbox{ in }~B^{(j)}_{2R(0)},\\
&\phi^{(j)}\cdot W_j=0~\mbox{ in }~\mathcal D^{(j)}_{2R},\\
\end{aligned}
\right.
\end{equation}
\begin{equation}\label{eqn-douter}
\left\{
\begin{aligned}
&\partial_t \psi=\Delta_x\psi+\mathcal G[v,\phi^{(1)},\phi^{(2)},\psi,Z^*,\la^{(1)},\la^{(2)},\omega,\xi^{(1)},\xi^{(2)}]~\mbox{ in }~\R^2_+\times(0,T),\\
&\partial_{x_2}\psi_1=\partial_{x_2}\psi_2=\psi_3=0~\mbox{ on }~\pd\R_+^2\times(0,T),\\
&\psi(\cdot,0)=0~\mbox{ in }~\R^2_+,
\end{aligned}
\right.
\end{equation}
where
\begin{equation*}
\mathcal F[v,\phi^{(1)},\phi^{(2)},\Phi_{\rm out},\la^{(1)},\la^{(2)},\omega,\xi^{(1)},\xi^{(2)}]=\nabla u\odot \nabla u-\frac12|\nabla u|^2\mathbb I_2
\end{equation*}
with
\begin{equation*}
\begin{aligned}
u=&~(1+a)U_*+\Phi-(\Phi\cdot U_*)U_*,\\
\Phi=&~\eta_R^{(1)}\phi^{(1)}(y_1,t)+\eta_R^{(2)}Q_{\omega}\phi^{(2)}(y_2,t)+\psi+Z^*+\Phi_0,
\end{aligned}
\end{equation*}
\begin{equation*}
\begin{aligned}
&~h^{(1)}[v,\phi^{(1)},\phi^{(2)},\psi,Z^*,\la^{(1)},\la^{(2)},\omega,\xi^{(1)},\xi^{(2)}]\\
=&~
(\la^{(1)})^2 \Big[\tilde L_U^{(1)}[\Phi_{\rm out}]+\mathcal K_0^{(1)}+\mathcal K_1^{(1)}-(\la^{(1)})^{-1}\Pi^{(1)}_{U^{\perp}}(v \cdot \nabla_{y_1} u)\Big],\\
&~h^{(2)}[v,\phi^{(1)},\phi^{(2)},\psi,Z^*,\la^{(1)},\la^{(2)},\omega,\xi^{(1)},\xi^{(2)}]\\
=&~
(\la^{(2)})^2 Q^{-1}_{\omega}\Big[\tilde L^{(2)}_U[\Phi_{\rm out}]+\mathcal K_0^{(2)}+\mathcal K_1^{(2)}-(\la^{(2)})^{-1}\Pi^{(2)}_{U^{\perp}}(v \cdot \nabla_{y_2} u)\Big],\\
\end{aligned}
\end{equation*}
and
\begin{equation*}
\begin{aligned}
&~\mathcal G[v,\phi^{(1)},\phi^{(2)},\psi,Z^*,\la^{(1)},\la^{(2)},\omega,\xi^{(1)},\xi^{(2)}]\\:=&~(1-\eta_R^{(1)}-\eta_R^{(2)})\tilde L^{(j)}_{U}\Phi_{\rm out}+\mathcal C_{\rm in}+\mathcal N\\
&~+\left[v\cdot \nabla u-\eta_R^{(1)}\Pi^{(1)}_{U^{\perp}}(v\cdot \nabla u)-\eta_R^{(2)}\Pi^{(2)}_{U^{\perp}}(v\cdot \nabla u)\right].
\end{aligned}
\end{equation*}
Here the coupling $\mathcal C_{\rm in}$ and the nonlinear term $\mathcal N$ are defined in \eqref{def-cpls} and \eqref{def-nlts}, respectively.

As discussed in Section \ref{sec-lineartheories}, suitable inner solution with space-time decay can be obtained under orthogonalities by adjusting the modulation parameters $\la^{(1)}$, $\la^{(2)}$, $\omega$, $\xi^{(1)}$, $\xi^{(2)}$ at corresponding modes. To solve the inner problems \eqref{eqn-dinner} in Fourier modes, we further decompose
$$h^{(j)}= h_1^{(j)}+ h_2^{(j)}+ h_3^{(j)}+ h_4^{(j)}$$
with
\begin{equation*}
\begin{aligned}
&\qquad h^{(1)}_1[v,\phi^{(1)},\phi^{(2)},\psi,Z^*,\la^{(1)},\la^{(2)},\omega,\xi^{(1)},\xi^{(2)}]\\
&=\left(
(\lambda^{(1)})^2   \left(
[\tilde L_U]^{(1)}_0  [\Phi_{\rm out} ]
+[\tilde L_U]^{(1)}_2  [\Phi_{\rm out}] +\KK_{0}^{(1)}\right)
 +\lambda^{(1)}   [\Pi^{(1)}_{U^{\perp}}(v\cdot\nabla u)]_0\right)\chi_{\mathcal D^{(1)}_{2R} },
\\
&\qquad h^{(1)}_2[v,\phi^{(1)},\phi^{(2)},\psi,Z^*,\la^{(1)},\la^{(2)},\omega,\xi^{(1)},\xi^{(2)}]\\ &=\left(
(\lambda^{(1)})^2   \left([\tilde L_U]^{(1)}_1 [\Phi_{\rm out} ](0)
+\KK_{1}^{(1)}\right)+\lambda^{(1)}  \left([\Pi^{(1)}_{U^{\perp}}(v\cdot\nabla u)]_{1}+ [\Pi^{(1)}_{U^{\perp}}(v\cdot\nabla u)]_{\perp}\right)\right) \chi_{\mathcal D^{(1)}_{2R} },
\\
&\qquad h^{(1)}_3[v,\phi^{(1)},\phi^{(2)},\psi,Z^*,\la^{(1)},\la^{(2)},\omega,\xi^{(1)},\xi^{(2)}] \\
&= (\lambda^{(1)})^2
\left( [\tilde L_U]^{(1)}_1 [\Phi_{\rm out} ] - [\tilde L_U]^{(1)}_1  [\Phi_{\rm out}] (0)\right) \chi_{\mathcal D^{(1)}_{2R} },\\
&\qquad h^{(2)}_1[v,\phi^{(1)},\phi^{(2)},\psi,Z^*,\la^{(1)},\la^{(2)},\omega,\xi^{(1)},\xi^{(2)}]\\
&=\left(
(\lambda^{(2)})^2  Q^{-1}_{\omega} \left(
[\tilde L_U]^{(2)}_0  [\Phi_{\rm out} ]
+[\tilde L_U]^{(2)}_2  [\Phi_{\rm out}] +\KK_{0}^{(2)}\right)
 +\lambda^{(2)}  Q^{-1}_{\omega} [\Pi^{(2)}_{U^{\perp}}(v\cdot\nabla u)]_0\right)\chi_{\mathcal D^{(2)}_{2R} },
\\
&\qquad h^{(2)}_2[v,\phi^{(1)},\phi^{(2)},\psi,Z^*,\la^{(1)},\la^{(2)},\omega,\xi^{(1)},\xi^{(2)}]\\ &=\left(
(\lambda^{(2)})^2  Q^{-1}_{\omega} \left([\tilde L_U]^{(2)}_1 [\Phi_{\rm out} ](0)
+\KK_{1}^{(2)}\right)+\lambda^{(2)}  Q^{-1}_{\omega}\left([\Pi^{(2)}_{U^{\perp}}(v\cdot\nabla u)]_{1}+ [\Pi^{(2)}_{U^{\perp}}(v\cdot\nabla u)]_{\perp}\right)\right) \chi_{\mathcal D^{(2)}_{2R} },
\\
&\qquad h^{(1)}_3[v,\phi^{(1)},\phi^{(2)},\psi,Z^*,\la^{(1)},\la^{(2)},\omega,\xi^{(1)},\xi^{(2)}] \\
&= (\lambda^{(2)})^2  Q^{-1}_{\omega}
\left( [\tilde L_U]^{(2)}_1 [\Phi_{\rm out} ] - [\tilde L_U]^{(2)}_1  [\Phi_{\rm out}] (0)\right) \chi_{\mathcal D^{(2)}_{2R} },\\
\end{aligned}
\end{equation*}
where $[\Pi^{(j)}_{U^\perp}(v\cdot\nabla u)]_{0}$,  $[\Pi^{(j)}_{U^\perp}(v\cdot\nabla u)]_{1}$ and $[\Pi^{(j)}_{U^\perp}(v\cdot\nabla u)]_{\perp}$ stand respectively for the projection on modes $0$, $1$ and higher modes $|k|\geq 2$ defined in \eqref{def-Fourier1}--\eqref{def-Fourier3} (notice that there is no projection on mode $-1$ since $v\cdot\nabla u$ is in the right symmetry class), and
 \begin{equation*}
\begin{aligned}
{[\tilde L_U]}_1^{(1)}[\psi](0)& =
   -\,
  2(\la^{(1)})^{-1} w_{\rho_1}  \cos w(\rho_1) \, \big [\,(\partial_{x_1} \psi_2(\xi^{(1)}(t),t)) \cos \theta_1 +    (\partial_{x_2} \psi_2(\xi^{(1)}(t),t)) \sin  \theta_1 \, \big ]\, E_1^{(1)}
\\
& \quad + 2(\la^{(1)})^{-1}  w_{\rho_1}  \cos w(\rho_1)  \, \big [\, (\partial_{x_2} \psi_2(\xi^{(1)}(t),t)) \cos \theta_1 -    (\partial_{x_1} \psi_2(\xi^{(1)}(t),t)) \sin  \theta_1 \, \big ]\, E_2^{(1)},\\
[\tilde L_U]_1^{(2)} [\psi ](0) & =
   -\,
  2(\la^{(2)})^{-1} w_{\rho_2}  \cos w(\rho_2) \, \big [\,(\partial_{x_1} \psi_3(\xi^{(2)}(t),t)) \cos \theta_2 +    (\partial_{x_2} \psi_3(\xi^{(2)}(t),t)) \sin  \theta_2 \, \big ]\,Q_{\omega} E_1^{(2)}
\\
& \quad - 2(\la^{(2)})^{-1}  w_{\rho_2}  \cos w(\rho_2)  \, \big [\, (\partial_{x_1} \psi_3(\xi^{(2)}(t),t)) \sin \theta_2 -    (\partial_{x_2} \psi_3(\xi^{(2)}(t),t)) \cos  \theta_2 \, \big ]\, Q_{\omega}E_2^{(2)}.
\end{aligned}
\end{equation*}
Then for $j=1\,2$, by decomposing $\phi^{(j)}=\phi_1^{(j)}+\phi_2^{(j)}+\phi_3^{(j)}$ in a similar manner as $h^{(j)}_i$'s, the inner problems \eqref{eqn-dinner} become
\begin{equation}
\left\{
\begin{aligned}
(\lambda^{(j)})^2  \partial_t \phi_1^{(j)}
&= L_W^{(j)} [\phi_1^{(j)}] + h_1^{(j)}[v,\phi^{(1)},\phi^{(2)},\Phi_{\rm out},\la^{(1)},\la^{(2)},\omega,\xi^{(1)},\xi^{(2)}]\\
& \qquad
-  \sum_{  \ell =1,2} \tilde c_{0\ell}^{(j)}[ h_1^{(j)}[v,\phi^{(1)},\phi^{(2)},\Phi_{\rm out},\la^{(1)},\la^{(2)},\omega,\xi^{(1)},\xi^{(2)}]] w_{\rho_j}^2 Z^{(j)}_{0,\ell}
\\
\label{eqphi1}
& \qquad
-  \sum_{  \ell=1,2} c^{(j)}_{1\ell}[ h_1^{(j)}[v,\phi^{(1)},\phi^{(2)},\Phi_{\rm out},\la^{(1)},\la^{(2)},\omega,\xi^{(1)},\xi^{(2)}]] w_{\rho_1}^2  Z^{(j)}_{1,\ell}
\inn \mathcal D^{(j)}_{2R}
\\
\phi_1^{(j)}\cdot  W_j  &=  0   \inn \mathcal D^{(j)}_{2R}
\\
\phi_1^{(j)}(\cdot, 0) &=0 \inn B^{(j)}_{2R(0)}
\end{aligned}
\right.
\end{equation}

\begin{align}
\left\{
\begin{aligned}
\label{eqphi2}
(\lambda^{(j)})^2 \partial_t \phi_2^{(j)}
&= L_W^{(j)} [\phi _2^{(j)}]  + h_2^{(j)}[v,\phi^{(1)},\phi^{(2)},\Phi_{\rm out},\la^{(1)},\la^{(2)},\omega,\xi^{(1)},\xi^{(2)}]\\
& \qquad
 -  \sum_{  \ell=1,2} c^{(j)}_{1\ell}[  h_2^{(j)}[v,\phi^{(1)},\phi^{(2)},\Phi_{\rm out},\la^{(1)},\la^{(2)},\omega,\xi^{(1)},\xi^{(2)}]  ] w_{\rho_j}^2  Z^{(j)}_{1,\ell}    \inn \mathcal D^{(j)}_{2R}
\\
\phi_2^{(j)}\cdot  W_j  & = 0   \inn \mathcal D^{(j)}_{2R}
\\
\phi_2^{(j)}(\cdot, 0) & =0 \inn B^{(j)}_{2R(0)}
\end{aligned}
\right.
\end{align}

\begin{align}
\label{eqphi3}
\left\{
\begin{aligned}
(\lambda^{(j)})^2  \partial_t  \phi_3^{(j)} &= L_W^{(j)} [\phi_3^{(j)}] +
h_3^{(j)}[v,\phi^{(1)},\phi^{(2)},\Phi_{\rm out},\la^{(1)},\la^{(2)},\omega,\xi^{(1)},\xi^{(2)}]\\
& \qquad - \sum_{  \ell=1,2} c^{(j)}_{1\ell}[  h_3^{(j)}[v,\phi^{(1)},\phi^{(2)},\Phi_{\rm out},\la^{(1)},\la^{(2)},\omega,\xi^{(1)},\xi^{(2)}] ] w_{\rho_j}^2 Z^{(j)}_{1,\ell}
\\
& \quad
+ \sum_{\ell=1,2} c_{0\ell}^{*(j)}[v,\phi^{(1)},\phi^{(2)},\Phi_{\rm out},\la^{(1)},\la^{(2)},\omega,\xi^{(1)},\xi^{(2)}] w_{\rho_j}^2 Z^{(j)}_{0,\ell} \inn \mathcal D^{(j)}_{2R}
\\
\phi_3^{(j)}\cdot  W_j  & = 0   \inn \mathcal D^{(j)}_{2R}
\\
\phi_3^{(j)}(\cdot, 0) & =0 \inn B^{(j)}_{2R(0)}
\end{aligned}
\right.
\end{align}

\begin{align}
\label{redueqn-m0}
&c^{*(j)}_{0\ell}(t) - \tilde c_{0\ell}^{(j)} (t) = 0 ~\mbox{ for all }~ t\in (0,T), \quad (j,\ell)=\big\{(1,1),(2,1),(2,2)\big\},
\\
\label{redueqn-m1}
&c_{1\ell}^{(j)}(t)  =  0 ~\mbox{ for all }~ t\in (0,T), \quad (j,\ell)=\big\{(1,1),(2,1),(2,2)\big\}.
\end{align}

\bigskip

\subsection{Weighted topologies and fixed point argument}

\medskip

We now design the weighted topologies for inner solutions $\phi^{(1)}$, $\phi^{(2)}$, outer solution $\Phi_{\rm out}$ and parameters $(\la^{(1)},\la^{(2)},\omega,\xi^{(1)},\xi^{(2)})$.

\medskip

According to the linear theories in Section \ref{sec-lineartheories}, we shall solve the outer problem \eqref{eqn-douter}, inner problems \eqref{eqphi1}--\eqref{eqphi3} and forced Navier-Stokes equation \eqref{eqn-NS} in the norms below

\medskip

\begin{itemize}
\item We use the norm $\|\cdot\|_{**}$ defined in \eqref{defNormRHSpsi} to measure the right hand side $\mathcal G$ in the outer problem \eqref{eqn-douter}.

\item We use the norm $\| \cdot\|_{\sharp, \Theta,\gamma}$ defined in \eqref{normPsi} to measure the solution $\psi$ solving the outer problem \eqref{eqn-douter}, where $\Theta>0$ and $\gamma\in(0,1/2)$.

\item We use the norm $\|\cdot\|^{(j)}_{\nu_i,a_i}$ to measure the right hand side $h_i^{(j)}$ with $i=1,\cdots,3$, where
\begin{align}
\label{norm-h}
\|h_i^{(j)}\|^{(j)}_{\nu_i,a_i}  =
\sup_{\R_+^2 \times (0,T)} \  \frac{ |h_1^{(j)}(y_j,t)| }{ \lambda_*^{\nu_i}(t) (1+\rho_j)^{-a_i}}
\end{align}
with $\nu_i>0$, $a_i\in(2,3)$ for $i=1,\,2$, and $a_3\in (1,3)$.

\item We use the norm $\|\cdot\|^{(j)}_{*,\nu_1,a_1,\delta}$ to measure the solution $\phi_1^{(j)}$ solving \eqref{eqphi1}, where
\begin{align*}
\| \phi_1^{(j)} \|^{(j)}_{*,\nu_1,a_1,\delta}
=
\sup_{\mathcal D^{(j)}_{2R}}
\frac{| \phi_1^{(j)}(y_j,t) | + (1+\rho_j) |\nabla_{y_j} \phi_1^{(j)}(y_j,t)|+(1+\rho_j)^2 |\nabla^2_{y_j} \phi_1^{(j)}(y_j,t)|}{ \lambda_*^{\nu_1}(t) \max\left\{\frac{R^{\delta(5-a_1)}}{(1+\rho_j)^3} , \frac{1}{(1+\rho_j)^{a_1-2} }\right\}}
\end{align*}
with $\nu_1\in(0,1)$, $a_1\in(2,3)$, $\delta>0$ fixed small.

\item We use the norm $\|\cdot\|^{(j)}_{{\rm in},\nu_2,a_2-2}$ to measure the solution $\phi_2^{(j)}$ solving \eqref{eqphi2}, where
$$
\|\phi_2^{(j)}\|^{(j)}_{{\rm in},\nu_2,a_2-2}  =
\sup_{\mathcal D^{(j)}_{2R}} \  \frac{| \phi_2^{(j)}(y_j,t) | + (1+\rho_j) |\nabla_{y_j} \phi_2^{(j)}(y_j,t)|+(1+\rho_j)^2 |\nabla^2_{y_j} \phi_2^{(j)}(y_j,t)|}{ \lambda_*^{\nu_2}(t) (1+\rho_j)^{2-a_2}}
$$
with $\nu_2\in(0,1)$, $a_2\in(2,3)$.

\item We use the norm $\|\cdot\|^{(j)}_{**,\nu_3}$ to measure the solution $\phi_3^{(j)}$ solving \eqref{eqphi3}, where
\begin{align*}
\|\phi_3^{(j)}\|^{(j)}_{**,\nu_3}
= \sup_{\mathcal D^{(j)}_{2R}} \
\frac{ |\phi(y_j,t)| + (1+\rho_j)\left |\nabla_{y_j} \phi_3^{(j)}(y_j,t)\right | +(1+\rho_j)^2 |\nabla^2_{y_j} \phi_3^{(j)}(y_j,t)|}
{ \la^{\nu_3}_*(t)  R^{2}(t) ( 1+\rho_j )^{-1}  }
\end{align*}
with $\nu_3>0$.
\item We will solve the incompressible Navier--Stokes equation \eqref{eqn-NS} in the norms $$\|\cdot\|^{(j)}_{S,\nu-1,1},\,
\|\cdot\|^{(j)}_{S,\nu-2,a+1}$$
for the velocity $v$ and forcing $\mathcal F$, respectively. Here
\begin{equation}\label{def-normSSF}
\begin{aligned}
\|F\|^{(j)}_{S,\nu-2,a+1}:=&~\sup_{(x,t)\in\R^2_+\times(0,T)} \la_*^{2-\nu}(t)\left(1+\left|\frac{x-q^{(j)}}{\la_*(t)}\right|^{a+1}\right)|F(x,t)|\\
 &~+ \sup_{(x,t)\in\R^2_+\times(0,T)} \la_*^{3-\nu}(t)\left(1+\left|\frac{x-q^{(j)}}{\la_*(t)}\right|^{a+2}\right)|\nabla_x F(x,t)|,
\end{aligned}
\end{equation}
and we require $\nu\in(0,1)$ and $a\in(1,2)$. 
\end{itemize}
\medskip
The way to solve \eqref{eqn-NS} is similar to that of the transported harmonic map flow. We derive $\R^2_+$ into three parts (the regions near two concentrations $q^{(j)}$ and the intermediate region), and then analyze the behaviors of the forcing and corresponding velocity field in each region.

\medskip

We solve the inner problems in the following weighted spaces:
\begin{align*}
\tilde E_1^{(j)} &= \{ \phi^{(j)}_1 \in L^\infty(\mathcal D^{(j)}_{2R}) :
\nabla_{y_j} \phi^{(j)}_1  \in L^\infty(\mathcal D^{(j)}_{2R}), \
\|\phi_1^{(j)}\|^{(j)}_{*,\nu_1,a_1,\delta} <\infty \},
\\
\tilde E_2^{(j)} &= \{ \phi^{(j)}_2 \in L^\infty(\mathcal D^{(j)}_{2R}) :
\nabla_{y_j} \phi_2^{(j)}  \in L^\infty(\mathcal D^{(j)}_{2R}), \
\|\phi_2^{(j)}\|^{(j)}_{{\rm in},\nu_2,a_2-2}<\infty \},
\\
\tilde E_3^{(j)} &= \{ \phi^{(j)}_3 \in L^\infty(\mathcal D^{(j)}_{2R}) :
\nabla_{y_j} \phi^{(j)}_3  \in L^\infty(\mathcal D^{(j)}_{2R}), \
\|\phi_3^{(j)}\|^{(j)}_{**,\nu_3} < \infty \},
\\
\end{align*}
and denote
\begin{equation*}
E_{\phi}^{(j)} =\tilde  E_1^{(j)}\times\tilde  E_2^{(j)} \times\tilde  E_3^{(j)} ,\quad \Phi_{\rm inner}^{(j)} = ( \phi_1^{(j)},\phi_2^{(j)},\phi_3^{(j)}) \in E^{(j)}_{\phi}
\end{equation*}
$$
\|\Phi_{\rm inner}^{(j)}\|_{E^{(j)}_{\phi}} =
\|\phi^{(j)}_1\|^{(j)}_{*,\nu_1,a_1,\delta}
+\|\phi_2^{(j)}\|^{(j)}_{{\rm in},\nu_2,a_2-2}
+\|\phi_3^{(j)}\|^{(j)}_{**,\nu_3}.
$$
We define the closed ball
\[
\mathcal B^{(j)} = \{  \Phi^{(j)}_{\rm inner} \in E_{\phi}^{(j)} : \| \Phi^{(j)}_{\rm inner}\|_{E^{(j)}_{\phi}} \leq 1 \} .
\]
For the outer problem, we introduce
\begin{equation*}
E_{\psi}=\left\{\psi\in L^{\infty}(\R^2_+\times(0,T)):\|\psi\|_{\sharp,\Theta,\gamma}<\infty\right\}.
\end{equation*}
For the incompressible Navier--Stokes equation, we shall solve the velocity field $v$ in the space
\begin{equation}\label{class-v}
E_v=\left\{v\in L^2(\R^2_+;\R^2): \nabla\cdot v=0,~\|v\|_{S,\nu-1,1}<M \varepsilon_0\right\},
\end{equation}
where $0<\varepsilon_0\ll 1$ is the universal number in \eqref{LCF}, and $M>0$ is some fixed number.

\medskip

For the parameters $$\la^{(1)}(t),\,  p_2(t)=\la^{(2)}(t)e^{i\omega(t)},\, \xi^{(1)}(t)=(\xi^{(1)}_1(t),0),\,\xi^{(2)}(t)=(\xi^{(2)}_1(t),\,\xi^{(2)}_2(t)),$$
we set
\begin{equation*}
\begin{aligned}
& X_{\la^{(1)}} := \{  \la^{(1)} \in C([-T,T;\R]) \cap C^1([-T,T;\R]) \ : \ \la^{(1)}(T) = 0 , \ \big\|\la^{(1)}-\kappa^{(1)}\la_*(t)\big\|_{*,3-\sigma}<\infty \} ,\\
& X_{p_2} := \{  p_2 \in C([-T,T;\mathbb C]) \cap C^1([-T,T;\mathbb C]) \ : \ p_2(T) = 0 , \ \big\|p_2-\kappa^{(2)}e^{i\omega_0}\la_*(t)\big\|_{*,3-\sigma}<\infty \} ,\\
& X_{\xi^{(1)}}=\left\{\xi_1^{(1)}\in C^1((0,T);\R):\dot\xi_1^{(1)}(T)=0,\|\xi^{(1)}\|_{X_{\xi^{(1)}}}<\infty\right\},\\
& X_{\xi^{(2)}}=\left\{\xi^{(2)}\in C^1((0,T);\R^2):\dot\xi^{(2)}(T)=0,\|\xi^{(2)}\|_{X_{\xi^{(2)}}}<\infty\right\},\\
 \end{aligned}
\end{equation*}
where $\| \cdot \|_{*,3-\sigma}$ is defined by
\begin{align}
\nonumber
\|g\|_{*,3-\sigma} := \sup_{t\in [-T,T]}  |\log(T-t)|^{3-\sigma} |\dot g(t)|,
\end{align}
and
\begin{equation*}
\|\xi^{(j)}\|_{X_{\xi^{(j)}}}=\|\xi^{(j)}\|_{L^{\infty}(0,T)}+\sup_{t\in(0,T)} \la_*^{-\sigma}(t)|\dot\xi^{(j)}(t)|
\end{equation*}
for some $\sigma\in(0,1)$.

\medskip

By similar computations as those in \cite[Section 4]{LCF2D}, under following restrictions on the constants in the norms introduced above
\begin{equation}\label{choice-outer}
\left\{
\begin{aligned}
&0<\Theta<\min\left\{\gamma_*,\frac12-\gamma_*,\nu_1-1+\gamma_*(a_1-1),\nu_2-1+\gamma_*(a_2-1),\nu_3-1\right\},\\
&\Theta<\min\left\{\nu_1-\delta\gamma_*(5-a_1)-\gamma_*,\nu_2-\gamma_*,\nu_3-3\gamma_*,\right\},\\
&\delta\ll 1,\quad \nu>\frac12,\\
\end{aligned}
\right.
\end{equation}
we have

\medskip

\begin{prop}\label{prop-outerop}
Assume \eqref{choice-outer} hold true. If $T>0$ is sufficiently small, then there exists a solution $\psi=\Psi(v,\Phi^{(1)}_{\rm inner},\Phi^{(2)}_{\rm inner},\la^{(1)},p_2,\xi^{(1)},\xi^{(2)})$ to the outer problem \eqref{eqn-douter} with
\begin{align*}
&\qquad\|\Psi(v,\Phi^{(1)}_{\rm inner},\Phi^{(2)}_{\rm inner},\la^{(1)},p_2,\xi^{(1)},\xi^{(2)})\|_{\sharp,\Theta,\gamma}\\
&\lesssim T^{\epsilon}\sum_{j=1}^2 \bigg(\|v\|^{(j)}_{S,\nu-1,1}+\|\Phi^{(j)}_{\rm inner}\|_{E_{\phi}^{(j)}}+\|\la^{(1)}\|_{X_{\la^{(1)}}}+\|p_2\|_{X_{p_2}}+\|\xi^{(j)}\|_{X_{\xi^{(j)}}}+1\bigg),
\end{align*}
for some $\epsilon>0$.
\end{prop}
We define $\mathcal T_{\psi}$ by the operator which solves $\psi$ in Proposition \ref{prop-outerop}. For the inner problems \eqref{eqphi1}--\eqref{eqphi3}, we then take $\Phi^{(j)}_{\rm inner}\in E_{\phi}^{(j)}$ and substitute
$$\Phi_{\rm out}(v,\Phi^{(1)}_{\rm inner},\Phi^{(2)}_{\rm inner},\la^{(1)},p_2,\xi^{(1)},\xi^{(2)}) = Z^*+ \Psi(v,\Phi^{(1)}_{\rm inner},\Phi^{(2)}_{\rm inner},\la^{(1)},p_2,\xi^{(1)},\xi^{(2)})$$
into inner problems \eqref{eqn-dinner}. We can then write equations \eqref{eqn-dinner}
as the fixed point problems
\begin{align}
\label{ptofijo}
\Phi_{\rm inner}^{(j)} = \mathcal A^{(j)} \big(\Phi_{\rm inner}^{(j)}\big),\quad j=1,\,2,
\end{align}
where
\begin{align*}
\mathcal A^{(j)} (\Phi^{(j)}_{\rm inner}) = ( \mathcal A^{(j)}_1(\Phi^{(j)}_{\rm inner}) ,  \mathcal A^{(j)}_2(\Phi_{\rm inner}^{(j)}) ,  \mathcal A_3^{(j)}(\Phi_{\rm inner}^{(j)}) ) , \quad \mathcal A^{(j)} : \bar{\mathcal B}^{(j)}_1\subset E^{(j)}_{\phi} \to E^{(j)}_{\phi}
\end{align*}
with
\begin{align*}
\mathcal A^{(j)}_1(\Phi^{(j)}_{\rm inner}) &=   \mathcal T^{(j)}_{1}  \bigg(
h_1[v,\Phi_{\rm out}(v,\Phi^{(1)}_{\rm inner},\Phi^{(2)}_{\rm inner},\la^{(1)},p_2,\xi^{(1)},\xi^{(2)}),\la^{(1)},p_2,\xi^{(1)},\xi^{(2)}  ] \bigg),
\\
\mathcal A^{(j)}_2(\Phi^{(j)}_{\rm inner}) &=   \mathcal T^{(j)}_{2}  \bigg(
h_2[v,\Phi_{\rm out}(v,\Phi^{(1)}_{\rm inner},\Phi^{(2)}_{\rm inner},\la^{(1)},p_2,\xi^{(1)},\xi^{(2)}),\la^{(1)},p_2,\xi^{(1)},\xi^{(2)}  ] \bigg),
\\
\mathcal A^{(j)}_3(\Phi^{(j)}_{\rm inner}) &=   \mathcal T^{(j)}_{3 }
\Bigg(
h_3[v,\Phi_{\rm out}(v,\Phi^{(1)}_{\rm inner},\Phi^{(2)}_{\rm inner},\la^{(1)},p_2,\xi^{(1)},\xi^{(2)}),\la^{(1)},p_2,\xi^{(1)},\xi^{(2)}  ]\\
&~\qquad+
\sum_{\ell=1}^2 c_{0\ell}^{*(j)}[v,\Phi_{\rm out}(v,\Phi^{(1)}_{\rm inner},\Phi^{(2)}_{\rm inner},\la^{(1)},p_2,\xi^{(1)},\xi^{(2)}),\la^{(1)},p_2,\xi^{(1)},\xi^{(2)}  ] w_{\rho_j}^2 Z_{0,\ell}^{(j)}
\Bigg).
\\
\end{align*}
Here $\mathcal T^{(j)}_{1}(\cdot)$, $\mathcal T^{(j)}_{2}(\cdot)$, $\mathcal T^{(j)}_{3}(\cdot)$ stand for the operators that solve the inner problems \eqref{eqphi1}, \eqref{eqphi2}, \eqref{eqphi3},  respectively.

\medskip

By the linear theories in Section \ref{sec-lineartheories}, it is direct to check (similar to \cite[Section 4]{LCF2D}) that the inner problems can be solved provided
\begin{equation}\label{choice-inner}
\left\{
\begin{aligned}
&\nu=\nu_1=\nu_2<\min\left\{1,1-\gamma_*(a_2-2)\right\},\\
&\nu_3<\min\left\{1+\Theta+2\gamma_*\gamma,\nu_1+\frac12 \delta\gamma_*(a_1-2)\right\},\\
&1<a<2,\\
&0<\varepsilon_0\ll 1.\\
\end{aligned}
\right.
\end{equation}
Here $\varepsilon_0\ll 1$ is required to ensure the implementation of the loop. More precisely, we have

\begin{prop}\label{prop-innerop}
Assume \eqref{choice-inner} are satisfied. If $T>0$ and $\varepsilon_0>0$ are sufficiently small,   then the system of equations \eqref{ptofijo} for $\Phi^{(j)}_{\rm inner}=(\phi^{(j)}_1,\phi^{(j)}_2,\phi^{(j)}_3)$ has a solution $\Phi^{(j)}_{\rm inner}\in  E_{\phi}^{(j)}$ for $j=1,\,2$.
\end{prop}

Above propositions together with the compactness from H\"older regularity complete the proof of Theorem \ref{thm1} by Schauder fixed point theorem.

\bigskip

\appendix

\section{Derivation of the partially free boundary system}\label{sec-deriv}

\medskip

We derive in this appendix the energy law and compatibility for the partially free boundary system \eqref{LCF}--\eqref{FB}.

\medskip

We first derive the energy law. Multiplying \eqref{LCF}$_1$ by $v$ and integrating over $\Omega\subset\R^d$ $(d\leq 3)$, we get
\begin{equation*}
\frac12 \frac{d}{dt}\int_{\Omega} |v|^2 + \int_{\Omega} (v\cdot \nabla v)\cdot v + \int_{\Omega} \nabla P\cdot v=-\int_{\Omega} |\nabla v|^2- \int_{\Omega} (\Delta u\cdot \nabla u)\cdot v,
\end{equation*}
where we have used
$$\nabla\cdot \left(\nabla u \odot \nabla u-\frac12 |\nabla u|^2 \mathbb{I}_d\right)=\Delta u\cdot \nabla u.$$
By \eqref{LCF}$_2$ and \eqref{FB}$_1$,
$$\int_{\Omega} (v\cdot \nabla v)\cdot v = \int_{\Omega} \nabla P\cdot v=0.$$
So we have
\begin{equation}\label{el-1}
\frac12 \frac{d}{dt}\int_{\Omega} |v|^2 =-\int_{\Omega} |\nabla v|^2- \int_{\Omega} (\Delta u\cdot \nabla u)\cdot v.
\end{equation}
Next we multiply \eqref{LCF}$_3$ with $\Delta u +|\nabla u|^2 u$ and integrate over $\Omega$
\begin{equation*}
-\frac12 \frac{d}{dt} \int_{\Omega} |\nabla u|^2 +\int_{\Omega} (v\cdot \nabla u)\cdot (\Delta u +|\nabla u|^2 u)=\int_{\Omega}\left|\Delta u +|\nabla u|^2 u\right|^2.
\end{equation*}
Since
$$
\int_{\Omega} (v\cdot \nabla u)\cdot (|\nabla u|^2 u)=\int_{\Omega} |\nabla u|^2 v\cdot \frac{\nabla(|u|^2)}{2}=0,
$$
we obtain
\begin{equation}\label{el-2}
-\frac12 \frac{d}{dt} \int_{\Omega} |\nabla u|^2 +\int_{\Omega} (\Delta u\cdot \nabla u)\cdot v=\int_{\Omega}\left|\Delta u +|\nabla u|^2 u\right|^2.
\end{equation}
Combining \eqref{el-1} and \eqref{el-2}, we get
\begin{equation}\label{energylaw}
\frac12 \frac{d}{dt}\left(\int_{\Omega}|v|^2+|\nabla u|^2\right)=-\int_{\Omega} |\nabla v|^2-\int_{\Omega}\left|\Delta u +|\nabla u|^2 u\right|^2
\end{equation}
which is called {\it the basic energy law} (see \cite{LL1995CPAM}). The energy law \eqref{energylaw} reflects the energy dissipation property of the flow of liquid crystals.

\medskip

On the other hand, the physical compatibility condition should be satisfied
\begin{equation}
\left<\left(\frac12(\nabla v+(\nabla v)^T)-P\mathbb I_d-\nabla u\odot \nabla u\right)\nu,\tau\right>=0\quad \mbox{ on }~\partial \Omega,
\end{equation}
where
$$
\nabla\cdot \left(\frac12(\nabla v+(\nabla v)^T)-P\mathbb I_d-\nabla u\odot \nabla u\right)
$$
is the {\it stress tensor}. It is easy to see that $<P\mathbb I_d \nu,\tau>=0$ as $<\nu,\tau>=0$. Also,
$$
\left<\frac12(\nabla v+(\nabla v)^T)\nu,\tau\right>=0
$$
is the {\it Navier boundary condition} \eqref{FB}$_2$, and
$$
0=\left<(\nabla u\odot \nabla u)\nu,\tau\right>=\left<\nabla_{\nu} u, \nabla_{\tau} u\right>
$$
implies the partially free boundary condition \eqref{FB}$_4$
$$
\frac{\partial u}{\partial \nu}\perp T_u \Sigma ~\mbox{ on }~\partial \Omega\times (0,T).
$$

\medskip

In conclusion, with the partially free boundary conditions \eqref{FB}, the system \eqref{LCF} is physically meaningful.

\bigskip

\section{Multiple bubbles: analysis of the interactions}\label{k-bubbles}

\medskip

In fact, we have three different cases for multiple bubbles, and we can take the following ansatz for each case
\begin{itemize}
\item multiple bubbles all placed in the interior
$$
u_*=\sum_{j=1}^k Q_{\omega_j}\left[W_1\left(\frac{x-\xi^{(j)}}{\la^{(j)}}\right)+W_1\left(\frac{x-(\xi^{(j)})^*}{\la^{(j)}}\right)\right]-(2k-1)W_1(\infty),
$$
where $\xi^{(j)}=(\xi^{(j)}_1,\xi^{(j)}_2)$ and $(\xi^{(j)})^*=(\xi^{(j)}_1,-\xi^{(j)}_2)$.
\medskip
\item multiple bubbles all placed on the boundary
$$
u_*=\sum_{j=1}^k W_1\left(\frac{x-\xi^{(j)}}{\la^{(j)}}\right)-(k-1)W_1(\infty) ~\mbox{ with }~\xi^{(j)}=(\xi^{(j)}_1,0).
$$
\medskip
\item mixed case: finite linear combination of interior and boundary bubbles
\begin{equation*}
\begin{aligned}
u_*=&~\sum_{j=1}^{k_{\mathcal B}} W_1\left(\frac{x-\xi_{\mathcal B}^{(j)}}{\la_{\mathcal B}^{(j)}}\right)-(k_{\mathcal B}-1)W_1(\infty)\\
&~+\sum_{j=1}^{k_{\mathcal I}} Q_{\omega_{\mathcal I,j}}\left[\left(W_1\left(\frac{x-\xi_{\mathcal I }^{(j)}}{\la_{\mathcal I }^{(j)}}\right)-W_1(\infty)\right)+\left(W_1\left(\frac{x-(\xi^{(j)}_{\mathcal I })^*}{\la_{\mathcal I }^{(j)}}\right)-W_1(\infty)\right)\right],
\end{aligned}
\end{equation*}
where $\xi_{\mathcal I }^{(j)}=(\xi^{(j)}_{\mathcal I,1},\xi^{(j)}_{\mathcal I,2})$, $(\xi_{\mathcal I }^{(j)})^*=(\xi^{(j)}_{\mathcal I,1},-\xi^{(j)}_{\mathcal I,2})$, and $\xi_{\mathcal B}^{(j)}=(\xi_{\mathcal B,1}^{(j)},0)$.
\medskip
\end{itemize}
Here $W_1$ is defined in \eqref{def-W1}. The purpose of the reflection terms in the first and third case is to enforce the symmetry \eqref{symmetryclass} for the full system. Observe that the first case is ``localized'' as the core inner region does not touch the boundary, while the second case is automatically in the symmetry class \eqref{symmetryclass}. So the most interesting case is the third case of mixed bubbles. Here we give some heuristic discussions about the reason for placing the reflected bubble in the mix case.

\medskip

\begin{remark}
Notice that above general ansatz is slightly different from what we take for the two-bubble case in Section \ref{Sec-approx}, and the choices are of course not unique. 
\end{remark}

\medskip

The key here is the interaction between the boundary bubble and the pair of the interior bubble and its reflection. More precisely, we need to analyze the error produced by the interior bubble and its associated reflection entering into the tangent plane of the boundary bubble.  Write
\begin{align*}
&~Q_{\omega}\left[W_1\left(\frac{x-\xi^{(2)}}{\la^{(2)}}\right)+W_1\left(\frac{x-(\xi^{(2)})^*}{\la^{(2)}}\right)\right]\\
=&~\begin{bmatrix}
       \cos\omega & -\sin\omega & 0 \\[0.3em]
       \sin \omega & \cos\omega  & 0 \\[0.3em]
       0 & 0 & 1\\
     \end{bmatrix}\begin{bmatrix}
     \frac{2\la^{(2)}(x_1-\xi^{(2)}_1)}{(x_1-\xi^{(2)}_1)^2+(x_2-\xi^{(2)}_2)^2+(\la^{(2)})^2}+\frac{2\la^{(2)}(x_1-\xi^{(2)}_1)}{(x_1-\xi^{(2)}_1)^2+(x_2+\xi^{(2)}_2)^2+(\la^{(2)})^2}\\
     \frac{(x_1-\xi^{(2)}_1)^2+(x_2-\xi^{(2)}_2)^2-(\la^{(2)})^2}{(x_1-\xi^{(2)}_1)^2+(x_2-\xi^{(2)}_2)^2+(\la^{(2)})^2}+\frac{(x_1-\xi^{(2)}_1)^2+(x_2+\xi^{(2)}_2)^2-(\la^{(2)})^2}{(x_1-\xi^{(2)}_1)^2+(x_2+\xi^{(2)}_2)^2+(\la^{(2)})^2}\\
     \frac{2\la^{(2)}(x_2-\xi^{(2)}_2)}{(x_1-\xi^{(2)}_1)^2+(x_2-\xi^{(2)}_2)^2+(\la^{(2)})^2}+\frac{2\la^{(2)}(x_2+\xi^{(2)}_2)}{(x_1-\xi^{(2)}_1)^2+(x_2+\xi^{(2)}_2)^2+(\la^{(2)})^2}\\
     \end{bmatrix}\\
     :=&~\tilde W_1+\tilde W_2,
\end{align*}
where $\xi^{(2)}=(\xi^{(2)}_1,\xi^{(2)}_2)$ with $\xi^{(2)}_2>0$.
Then the error
$$
S(\tilde W_1+\tilde W_2)=-\partial_t \tilde W_1-\partial_t \tilde W_2+|\nabla \tilde W_1|^2\tilde W_2+|\nabla \tilde W_2|^2\tilde W_1+2(\nabla \tilde W_1\cdot\nabla \tilde W_2)(\tilde W_1+\tilde W_2).
$$
In the proof of Theorem \ref{thm1}, a crucial observation is that the projection of the error onto $Q_{\omega}E_2$-direction of the tangent plane for the boundary linearization will in fact destroy the symmetry of the boundary inner problem across the $\partial\R^2_+$. Recall
$$
E_2=\begin{bmatrix} -\frac{x_2}{|x-\xi^{(1)}|}\\ 0\\ \frac{x_1-\xi^{(1)}_1}{|x-\xi^{(1)}|}\\ \end{bmatrix}.
$$
Here for notational simplicity, we have dropped the superscripts, and $\xi^{(1)}=(\xi^{(1)}_1,0)$ is the concentration point on $\partial\R^2_+$. Since $\tilde W_1$ and $\tilde W_2$ are symmetric about the boundary $\partial \R^2_+$, it is clear that
$$
(-\partial_t \tilde W_1-\partial_t \tilde W_2)\cdot Q_{\omega} E_2\Big|_{x_2=0}=0,
$$
$$
(\tilde W_1+\tilde W_2)\cdot Q_{\omega} E_2\Big|_{x_2=0}=0,
$$
and thus
$$
(2(\nabla\tilde  W_1\cdot\nabla \tilde W_2)(\tilde W_1+\tilde W_2))\cdot Q_{\omega} E_2\Big|_{x_2=0}=0.
$$
For the term $|\nabla \tilde W_1|^2 \tilde W_2+|\nabla \tilde W_2|^2\tilde W_1$, straightforward computations imply that on $\partial \R_+^2$
$$
|\nabla\tilde  W_1|^2=|\nabla \tilde W_2|^2,
$$
and thus
$$
(|\nabla \tilde W_1|^2\tilde W_2+|\nabla\tilde  W_2|^2\tilde W_1)\cdot Q_{\omega} E_2\Big|_{x_2=0}=(|\nabla \tilde W_1|^2(\tilde W_1+\tilde W_2))\cdot Q_{\omega} E_2\Big|_{x_2=0}=0.
$$
Therefore, we obtain
$$
S(\tilde W_1+\tilde W_2)\cdot Q_{\omega} E_2\Big|_{x_2=0}=0
$$
as desired.

\medskip

Certainly, one may deal with the inner problem touching the boundary without using the reflection, but this will rely on careful analysis on the non-degeneracy results and boundary linear theory, which might also be very interesting. Above formal analysis implies that the presence of reflected bubble can simplify the analysis of the inner problem near boundary concentration.

\bigskip

\section{Analysis of the Stokes operator}\label{sec-SO}

\medskip

As mentioned in Section \ref{sec-lineartheories}, another straightforward way to capture the precise pointwise estimates for the velocity field without using reflection is by explicit Green's tensor, which we now derive.

\medskip

We consider the following Stokes system with Navier boundary conditions
\begin{equation}\label{Stokes-Navier}
\left\{
\begin{aligned}
&\partial_t v  +\nabla P = \Delta v +F~&\mbox{ in }~\R^2_+\times (0,\infty),\\
&\nabla\cdot v =0~&\mbox{ in }~\R^2_+\times (0,\infty),\\
&\partial_{x_2} v_1=v_2=0~&\mbox{ on }~\partial\R^2_+\times (0,\infty),\\
&v\big|_{t=0}=0,\\
\end{aligned}
\right.
\end{equation}
where $F=(F_1,F_2)^T$ is solenoidal:
\begin{equation}\label{def-solenoidal}
\nabla\cdot F=0,\quad F_2\big|_{x_2=0}=0.
\end{equation}
Our aim is to construct Green's tensor and its associated pressure tensor to \eqref{Stokes-Navier}. To this end, we first consider the homogeneous system with $F=0$ and non-zero boundary condition on $x_2=0$
\begin{equation}\label{Stokes-Golovkin}
\left\{
\begin{aligned}
&\partial_t u  +\nabla P = \Delta u~&\mbox{ in }~\R^2_+\times (0,\infty),\\
&\nabla\cdot u =0~&\mbox{ in }~\R^2_+\times (0,\infty),\\
&\partial_{x_2} u_1\Big|_{x_2=0}=a_1(x_1,t),\quad u_2\Big|_{x_2=0}=a_2(x_1,t),\\
&u\big|_{t=0}=0,\\
&u\to 0 ~\mbox{ as }~|x|\to+\infty,\\
\end{aligned}
\right.
\end{equation}
and then use Duhamel's principle to get a solution for the non-homogeneous system \eqref{Stokes-Navier}. We will use Fourier transform in $x_1$ and Laplace transform in $t$
$$
\widetilde u(\xi_1,x_2,s)=\mathcal F_{x_1}\mathcal L_{t}[u]:=\int_{\R^1}\int_0^{\infty} e^{-ix_1\xi_1-st}u(x_1,x_2,t)dt dx_1.
$$
Then taking Fourier-Laplace transform on \eqref{Stokes-Golovkin}, we get
\begin{equation*}
\left\{
\begin{aligned}
&s\widetilde u_1+\xi_1^2 \widetilde u_1-\frac{d^2 \widetilde u_1}{dx_2^2}+i\xi_1 \widetilde P=0\\
&s\widetilde u_2+\xi_1^2 \widetilde u_2-\frac{d^2 \widetilde u_2}{dx_2^2}+\frac{d\widetilde P}{dx_2}=0\\
&i\xi_1\widetilde u_1+\frac{d\widetilde u_2}{dx_2}=0\\
&\frac{d\widetilde u_1}{d x_2}\Bigg|_{x_2=0}=\widetilde a_1(\xi_1,s), \quad\widetilde u_2\Bigg|_{x_2=0}=\widetilde a_2(\xi_1,s)\\
&\widetilde u\to 0~\mbox{ as }~x_2\to \infty.\\
\end{aligned}
\right.
\end{equation*}
We look for solution of the form
$$
u=u'+\nabla\varphi,
$$
where $u'$ is solenoidal which solves the heat equation and $\varphi$ is harmonic. Since $\widetilde u\to 0$ as $x_2\to \infty$, we have
$$
\widetilde{u'}=\theta_1(\xi_1,s)e^{-\sqrt{\xi_1^2+s}x_2},\quad \widetilde\varphi=\theta_2(\xi_1,s)e^{-|\xi_1|x_2}.
$$
Then
$$
\widetilde u=\phi(\xi_1,s)e^{-\sqrt{\xi_1^2+s}x_2}+\psi(\xi_1,s)e^{-|\xi_1|x_2},
$$
where
\begin{equation*}
\phi=\left(\phi_1,\frac{i\xi_1\phi_1}{\sqrt{\xi_1^2+s}}\right),\quad \psi=\left(i\xi_1,-|\xi_1|\right)\varphi(\xi_1,s).\\
\end{equation*}
The functions $\phi_1$ and $\varphi$ are determined by the boundary condition:
\begin{equation*}
\left\{
\begin{aligned}
&\phi_1\sqrt{\xi_1^2+s}+i\xi_1|\xi_1|\varphi=-\widetilde a_1,\\
&\frac{i\xi_1\phi_1}{\sqrt{\xi_1^2+s}}-|\xi_1|\varphi=\widetilde a_2.\\
\end{aligned}
\right.
\end{equation*}
Then
\begin{equation*}
\left\{
\begin{aligned}
&\phi_1=\frac{(i\xi_1\widetilde a_2-\widetilde a_1)\sqrt{\xi_1^2+s}}{s},\\
&\varphi=-\frac{i\xi_1\widetilde a_1+(\xi_1^2+s)\widetilde a_2}{s|\xi_1|},\\
\end{aligned}
\right.
\end{equation*}
and thus
\begin{equation}\label{FL-uP-N}
\begin{aligned}
\widetilde u_1(\xi_1,x_2,s)=&~\frac{(i\xi_1\widetilde a_2-\widetilde a_1)\sqrt{\xi_1^2+s}}{s}e^{-\sqrt{\xi_1^2+s}x_2}+\frac{\xi_1^2\widetilde a_1-i\xi_1(\xi_1^2+s)\widetilde a_2}{s|\xi_1|}e^{-|\xi_1|x_2},\\
\widetilde u_2(\xi_1,x_2,s)=&~-\frac{\xi_1^2\widetilde a_2+i\xi_1\widetilde a_1}{s}e^{-\sqrt{\xi_1^2+s}x_2}+\frac{i\xi_1 \widetilde a_1+(\xi_1^2+s)\widetilde a_2}{s}e^{-|\xi_1|x_2},\\
\widetilde P(\xi_1,x_2,s)=&~-s\varphi e^{-|\xi_1|x_2}=-\frac{i\xi_1\widetilde a_1+(\xi_1^2+s)\widetilde a_2}{|\xi_1|} e^{-|\xi_1|x_2},
\end{aligned}
\end{equation}
where $\widetilde a=(\widetilde a_1,\widetilde a_2)$.
Then we look for solution of the non-homogeneous problem \eqref{Stokes-Navier} with zero Navier boundary condition
$$
v=u+w,
$$
where
\begin{equation}\label{rep-www}
\begin{aligned}
&w_1=\int_0^t\int_{\R^2_+}\left[\Gamma(x-y,t-\tau)-\Gamma(x-y^*,t-\tau)\right]F_1(y,\tau)dyd\tau\\
&w_2=\int_0^t\int_{\R^2_+}\left[\Gamma(x-y,t-\tau)+\Gamma(x-y^*,t-\tau)\right]F_2(y,\tau)dyd\tau\\
\end{aligned}
\end{equation}
with $y^*=(y_1,-y_2)$ being the reflection of $y=(y_1,y_2)$. Then it is direct to see from the fact that $F$ is solenoidal that
\begin{equation}\label{solenoidal-heat-2}
\begin{cases}
\partial_t w  = \Delta w+F\\
\nabla\cdot w =0\\
\partial_{x_2} w_1\Big|_{x_2=0}=-b_1(x_1,t),\quad w_2\Big|_{x_2=0}=-b_2(x_1,t),\\
\end{cases}
\end{equation}
where
\begin{equation*}
\begin{aligned}
&b_1(x_1,t)=-\int_0^t\int_{\R^2_+}\partial_{x_2}\left[\Gamma(x-y,t-\tau)-\Gamma(x-y^*,t-\tau)\right]F_1(y,\tau)dyd\tau\Big|_{x_2=0},\\
&b_2(x_1,t)=0.\\
\end{aligned}
\end{equation*}
Therefore, from \eqref{solenoidal-heat-2}, $u$ solves the homogeneous problem \eqref{Stokes-Golovkin} with $a(x_1,t)=b(x_1,t)$, and
$$\widetilde a_1=\int_0^{\infty}e^{-\sqrt{\xi_1^2+s}y_2}\widetilde F_1(\xi_1,y_2,s)dy_2.$$
By \eqref{FL-uP-N}, one has
\begin{equation*}
\begin{aligned}
\widetilde u_1(\xi_1,x_2,s)=&~\frac{\widetilde a_1}{s}\frac{\partial}{\partial x_2}\left(e^{-\sqrt{\xi_1^2+s}x_2}-e^{-|\xi_1|x_2}\right),\\
\widetilde u_2(\xi_1,x_2,s)=&~-i\xi_1\frac{\widetilde a_1}{s}\left(e^{-\sqrt{\xi_1^2+s}x_2}-e^{-|\xi_1|x_2}\right),\\
\widetilde P(\xi_1,x_2,s)=&~-i\xi_1\widetilde a_1 \frac{e^{-|\xi_1|x_2}}{|\xi_1|}.\\
\end{aligned}
\end{equation*}
Taking inverse Fourier-Laplace transform, we obtain
\begin{equation}\label{u-P-Navier}
\begin{aligned}
u_1=&~\frac{\partial}{\partial x_2}\Bigg[-2\int_0^{t}ds\int_{\R^2_+}\frac{\partial \Gamma(x-y^*,s-\tau)}{\partial x_2} F_1(y,\tau)dy d\tau\\
&~-4\int_{\R}\frac{\partial E(x-y)}{\partial x_2}dy_1\int_0^t ds\int_0^s\int_{\R^2_+}\frac{\partial \Gamma(y-z,s-\tau)}{\partial z_2}F_1(z,\tau)dzd\tau\Big|_{y_2=0}\Bigg],\\
u_2=&~\frac{\partial}{\partial x_1}\Bigg[-2\int_0^{t}ds\int_{\R^2_+}\frac{\partial \Gamma(x-y^*,s-\tau)}{\partial x_2} F_1(y,\tau)dy d\tau\\
&~-4\int_{\R}\frac{\partial E(x-y)}{\partial x_2}dy_1\int_0^t ds\int_0^s\int_{\R^2_+}\frac{\partial \Gamma(y-z,s-\tau)}{\partial z_2}F_1(z,\tau)dzd\tau\Big|_{y_2=0}\Bigg],\\
P(x,t)=&~4\frac{\partial}{\partial x_1}\Bigg[\int_{\R^1}E(x-y)dy_1\int_0^t d\tau \int_{\R^2_+}\frac{\partial \Gamma(y-z,t-\tau)}{\partial z_2}F_1(z,\tau)dz\Big|_{y_2=0}\Bigg],
\end{aligned}
\end{equation}
where we have used the following inverse Fourier-Laplace transforms
\begin{equation}\label{F-L-table}
\begin{aligned}
&(\mathcal F_{x_1}\mathcal L_{t})^{-1}\left[e^{-|\xi_1|x_2}\right]=2\delta(t)\frac{\partial E(x)}{\partial x_2},\quad (\mathcal F_{x_1}\mathcal L_{t})^{-1}\left[\frac{e^{-|\xi_1|x_2}}{|\xi_1|}\right]=-2\delta(t)E(x)\\
&(\mathcal F_{x_1}\mathcal L_{t})^{-1}\left[e^{-\sqrt{\xi_1^2+s}x_2}\right]=-2\frac{\partial \Gamma(x,t)}{\partial x_2},\quad (\mathcal F_{x_1}\mathcal L_{t})^{-1}\left[\frac{e^{-\sqrt{\xi_1^2+s}x_2}}{\sqrt{\xi_1^2+s}}\right]=2\Gamma(x,t)\\
&(\mathcal F_{x_1}\mathcal L_{t})^{-1}\left[\frac1s\right]=\delta(x_1)\\
\end{aligned}
\end{equation}
with
 \begin{equation}
 E(z)=\frac{1}{2\pi}\log|z|, \qquad
 \Gamma(x,t)=\left\{
 \begin{aligned}
 &\frac{1}{4\pi t} e^{-\frac{|x|^2}{4t}},~&~t>0,\\
 &0,~&~ t<0.\\
 \end{aligned}
 \right.
 \end{equation}
Then from \eqref{rep-www}, the representation formula for Navier slip boundary is the following
\begin{equation}\label{v-P-Navier}
\begin{aligned}
v_1(x,t)=&~\int_0^t\int_{\R^2_+}\left[\Gamma(x-y,t-\tau)-\Gamma(x-y^*,t-\tau)\right]F_1(y,\tau)dyd\tau\\
&~+\frac{\partial}{\partial x_2}\Bigg[-2\int_0^{t}ds\int_{\R^2_+}\frac{\partial \Gamma(x-y^*,s-\tau)}{\partial x_2} F_1(y,\tau)dy d\tau\\
&~-4\int_{\R}\frac{\partial E(x-y)}{\partial x_2}dy_1\int_0^t ds\int_0^s\int_{\R^2_+}\frac{\partial \Gamma(y-z,s-\tau)}{\partial z_2}F_1(z,\tau)dzd\tau\Big|_{y_2=0}\Bigg],\\
v_2(x,t)=&~\frac{\partial}{\partial x_1}\Bigg[-2\int_0^{t}ds\int_{\R^2_+}\frac{\partial \Gamma(x-y^*,s-\tau)}{\partial x_2} F_1(y,\tau)dy d\tau\\
&~-4\int_{\R}\frac{\partial E(x-y)}{\partial x_2}dy_1\int_0^t ds\int_0^s\int_{\R^2_+}\frac{\partial \Gamma(y-z,s-\tau)}{\partial z_2}F_1(z,\tau)dzd\tau\Big|_{y_2=0}\Bigg],\\
P(x,t)=&~4\frac{\partial}{\partial x_1}\Bigg[\int_{\R}E(x-y)dy_1\int_0^t d\tau \int_{\R^2_+}\frac{\partial \Gamma(y-z,t-\tau)}{\partial z_2}F_1(z,\tau)dz\Big|_{y_2=0}\Bigg].
\end{aligned}
\end{equation}
Therefore, the Green's tensor and its associated pressure tensor for half space with Navier boundary conditions have been constructed:
\begin{prop}\label{prop-rep-N}
The solution to \eqref{Stokes-Navier} with solenoidal forcing can be expressed in the form
\begin{equation*}
\begin{aligned}
&v(x,t)=\int_0^t \int_{\R^2_+} \mathcal G^0(x,y,t-\tau)F(y,\tau)  dyd\tau+\int_0^t\int_{\R^2_+} \mathcal G^*(x,y,t-\tau)\int_0^{\tau}F(y,s)ds dyd\tau,\\
&P(x,t)=\int_0^t \int_{\R^2_+} \mathcal P(x,y,t-\tau)\cdot  F(y,\tau) dyd\tau\\
\end{aligned}
\end{equation*}
with $\mathcal G^0=(G_{ij}^0)_{i,j=1,2}$, $\mathcal G^*=(G_{ij}^*)_{i,j=1,2}$, $\mathcal P=(P_j)_{j=1,2}$, and
\begin{equation*}
\begin{aligned}
G^0_{ij}(x,y,t)=&~\delta_{ij}(\Gamma(x-y,t)-\Gamma(x-y^*,t)),\\
G^*_{ij}(x,y,t)=&~(1-\delta_{ij})\frac{\partial}{\partial x_1}\Bigg[-2\frac{\partial\Gamma(x-y^*,t)}{\partial x_2}-4\int_{\R}\frac{\partial E(x_1-z_1,x_2)}{\partial x_2}\frac{\partial\Gamma(z_1-y_1,y_2,t)}{\partial y_2}dz_1\Bigg]\\
&~+\delta_{ij}\frac{\partial}{\partial x_2}\Bigg[-2\frac{\partial\Gamma(x-y^*,t)}{\partial x_2}-4\int_{\R}\frac{\partial E(x_1-z_1,x_2)}{\partial x_2}\frac{\partial\Gamma(z_1-y_1,y_2,t)}{\partial y_2}dz_1\Bigg],\\
P_j(x,y,t)= &~4(1-\delta_{j2}) \frac{\partial}{\partial x_j}\Bigg[\int_{\R} E(x_1-z_1,x_2)\frac{\partial \Gamma(z_1-y_1,y_2,t)}{\partial y_2}dz_1\Bigg].
\end{aligned}
\end{equation*}
\end{prop}

To derive the pointwise estimates of the Green's tensor and pressure tensor as in \eqref{v-P-Navier}, we have the following lemma in the general case $\R^n_+$ $(n\geq 2)$, whose proof is similar to \cite{Solonnikov1976} in the no-slip boundary case (see also \cite[Proposition 2.3]{solonnikov2003estimates}).

\medskip

\begin{lemma}\label{lem-scaling}
Let $M(x,t)$ be a function defined for $x\in\R^n_+$ $(n\geq 2)$ and $t>0$ with the properties
\begin{equation*}
\begin{aligned}
&M(\la x,\la^2 t)=\la^q M(x,t) \quad \forall \la>0\\
&|D_x^k D_t^s M(x,t)|\lesssim t^{\frac{q-|k|-2s}{2}}e^{-c\frac{|x|^2}{t}}.\\
\end{aligned}
\end{equation*}
Then the integral
$$
J_i(x,y_n,t)=\int_{\R^{n-1}} \partial_{y_i} E(y) M(x'-y',x_n,t) dy'
$$
satisfies
\begin{equation}\label{est-modelJ1}
J_i(\la x, \la y_n, \la^2 t)=\la^{q}J_i(x,y_n,t),
\end{equation}
and
\begin{equation}\label{est-modelJ2}
|D^k_x D_{y_n}^l \partial_t^s J_i(x,y_n,t)|\lesssim t^{\frac{q+n-1-2s-k_n}{2}}[|x'|^2+(x_n+y_n)^2+t]^{-\frac{|k'|+l+n-1}{2}} e^{-\frac{cx_n^2}{t}},
\end{equation}
where $x'=(x_1,\dots,x_{n-1})$ and $k'=(k_1,\dots,k_{n-1})$.
\end{lemma}

As a consequence of the above proposition, we have the following pointwise estimates for the Green's tensor and pressure tensor.

\medskip

\begin{prop}
The pressure tensor and Green's tensor in Proposition \ref{prop-rep-N} have the following pointwise upper bounds
\begin{equation*}
\begin{aligned}
&|\partial_t^sD_x^kD_y^mP_j(x,y,t)|\lesssim t^{-1-s-\frac{m_2}{2}}(|x-y^*|^2+t)^{-\frac{1+|k|+|m'|}{2}}e^{-\frac{cy_2^2}{t}},\\
&|\partial_t^sD_x^kD_y^mG^*_{ij}(x,y,t)|\lesssim t^{-1-s-\frac{m_2}{2}}(|x-y^*|^2+t)^{-\frac{2+|k|+|m'|}{2}}e^{-\frac{cy_2^2}{t}}.
\end{aligned}
\end{equation*}
\end{prop}

\medskip

Let us now consider 
\begin{equation}\label{model1}
\begin{cases}
\partial_t v  +\nabla P = \Delta v +F~&\mbox{ in }~\R^2_+\times (0,T),\\
\nabla\cdot v =0~&\mbox{ in }~\R^2_+\times (0,T),\\
\partial_{x_2} v_1\big|_{x_2=0}=0,\quad v_2\big|_{x_2=0}=0,\quad v\big|_{t=0}=v_0,\\
\end{cases}
\end{equation}
where the forcing is not solenoidal. By the Helmholtz decomposition
$$
F=\mathbb P F +\mathbb Q F,
$$
where $\mathbb P F$ is a potential $\nabla \Phi_P$, and $\mathbb Q F$ is divergence-free, one can write
$$
\Phi_P(x,t)=-\int_{\R^2_+} \nabla_y N(x,y)\cdot F(y,t) dy,
$$
where
$$
N(x,y)=E(x-y)+E(x-y^*)
$$
 is the Green function of the Neumann problem for the Laplace operator in the half-space with $E$ defined as
 $$
  E(z)=\frac{1}{2\pi}\log|z|,
 $$
and $y^*=(y_1,-y_2)$. Thus,
\begin{equation}\label{Helmholtz}
\mathbb Q F=F+\nabla \int_{\R^2_+} \nabla_y N(x,y)\cdot F(y,t) dy
\end{equation}
is solenoidal, i.e.,
$$
\div (\mathbb Q F)=0,\quad \mathbb Q F\big|_{x_2=0}=0.
$$
Then, a solution to the model problem \eqref{model1} is defined by the following representation formulae
\begin{equation}\label{GreenNavier-v}
\begin{aligned}
v(x,t)=&~\int_{\R^2_+} \mathcal G^0(x,y,t) v_0(y) dy +\int_0^t \int_{\R^2_+} \mathcal G^*(x,y,t-s) v_0(y) dyds\\
&~+\int_0^t\int_{\R^2_+} \mathcal G^0(x,y,t-s) \mathbb Q F(y,s) dy ds+\int_0^t \int_{\R^2_+} \mathcal G^*(x,y,t-\tau) \int_0^{\tau} \mathbb Q F(y,s)ds dyd\tau,\\
\end{aligned}
\end{equation}
\begin{equation}\label{GreenNavier-P}
P(x,t)=\int_{\R^2_+} \mathcal P(x,y,t)\cdot v_0(y) dy+\Phi_P(x,t) +\int_0^t\int_{\R^2_+} \mathcal P(x,y,t-s)\cdot \mathbb Q F(y,s) dy ds,
\end{equation}
where $\mathcal G^0$, $\mathcal G^*$ and $\mathcal P$ are given in Proposition \ref{prop-rep-N}.

\bigskip

\section*{Acknowledgements}

F. Lin is partially supported by the NSF grant DMS-1955249. J. Wei is partially supported by NSERC of Canada.
Y. Sire is partially supported by the Simons foundation.

\bigskip

\bibliographystyle{plain}






\end{document}